\documentclass[12pt]{article}
\usepackage[psamsfonts]{amssymb}
\usepackage{amsfonts}
\usepackage[mathcal]{eucal}
\newcommand{\g}{\mathfrak{g}}
\renewcommand{\a}{\mathfrak{a}}
\renewcommand{\b}{\mathfrak{b}}
\newcommand{\h}{\mathfrak{h}}
\renewcommand{\k}{\mathfrak{k}}
\newcommand{\q}{\mathfrak{q}}
\newcommand{\p}{\mathfrak{p}}
\newcommand{\n}{\mathfrak{n}}
\newcommand{\m}{\mathfrak{m}}
\renewcommand{\l}{\mathfrak{l}}
\renewcommand{\j}{\mathfrak{j}}
\newcommand{\s}{\mathfrak{s}}
\renewcommand{\t}{\mathfrak{t}}
\newcommand{\z}{\mathfrak{z}}
\newcommand{\N}{\mathbb{N}}
\newcommand{\Q}{\mathbb{Q}}

\newcommand{\R}{\mathbb{R}}
\newcommand{\C}{\mathbb{C}}

\newcommand{\qed}{{\null\hfill\ \raise3pt\hbox{\framebox[0.1in]{}}}
\break\null}
\newtheorem{theo}{Th\'eor\`eme}
\newtheorem{prop}{Proposition}
\newtheorem{lem}{Lemme}
\newtheorem{cor}{Corollaire}
\newtheorem{rem}{Remarque}
\newtheorem{defi}{D\'efinition}
\newcommand{\ste}{\hfill\break}

\begin{document}
\def\b{\mathfrak{b}}
\def\c{\mathfrak{c}}
\def\d{\mathfrak{d}}
\def\e{\mathfrak{e}}
\def\f{\mathfrak{f}}
\def\h{\mathfrak{h}}
\def\i{\mathfrak{i}}
\def\j{\mathfrak{j}}
\def\r{\mathfrak{r}}
\def\k{\mathfrak{k}}
\def\q{\mathfrak{q}}
\def\p{\mathfrak{p}}
\def\n{\mathfrak{n}}
\def\m{\mathfrak{m}}
\def\l{\mathfrak{l}}
\def\j{\mathfrak{j}}
\def\s{\mathfrak{s}}
\def\t{\mathfrak{t}}

\def\u{\mathfrak{u}}
\def\v{\mathfrak{v}}
\def\w{\mathfrak{w}}
\def\x{\mathfrak{x}}
\def\K{\mathbb{K}}

\def\dem{ {\em D\'emonstration : \ste }}
\def\sac{ sous-alg\`ebre de Cartan }
\def\beq{\begin{equation}}
\def\eeq{\end{equation}}
\def\z{\mathfrak{z}}
\textheight 659pt\textwidth  444pt
\oddsidemargin  -1mm
\evensidemargin -1mm
\topmargin      -8mm
\pagestyle{plain}

\title{\bf Classification des triples de Manin pour les alg\`ebres
 de Lie
r\'eductives complexes
}
\author{Patrick Delorme}

\maketitle
\setcounter{section}{-1}
\section{Introduction}

Let $\g$ be a finite dimensional, complex,  reductive Lie algebra.
One says that a  symmetric, $\g$-invariant,
 $\R$(resp. $\C$)-bilinear form on $\g$ is a {\bf Manin form} if and
only if its signature is
$(dim_{\C}\g, dim_{\C}\g)$ (resp. is non degenerate).  \ste
Recall that a {\bf Manin-triple} in $\g$ is a triple $(B,\i, \i')$,
where
$B$ is a real (resp. complex) Manin form and  where
$\i,
\i'$ are real (resp. complex) Lie subalgebras of $\g$, isotropic for
$B$, and such that $\i+ \i'=\g$. Then this is a direct sum and $\i,
\i'$ are of real dimension equal to the complex dimension of $\g$.
\ste    Our goal is to classify all Manin-triples of $\g$, up to
conjugacy under the action on
$\g$ of the connected, simply connected Lie group $G$ with Lie
algebra $\g$, by induction on the rank of the derived algebra of
$\g$.\ste
One calls {\bf af-involution (resp. f-involution)} of a complex
semi-simple Lie algebra
$\m$, any $\R$(resp. $\C$)-linear involutive automorphism of $\m$,
$\sigma$, such that
 there exists :\ste 1) an ideal ${\tilde  \m_{0}}$ of $\m$, which
is reduced to zero for f-involutions,  and a real form
${\tilde \h}_{0}$ of ${\tilde \m_{0}},$\ste   2) simple  ideals of
$\m$, $\m'_{j}$,
$\m''_{j}$, $j=1,\dots,p$,
\ste
3)  an  isomorphism of   complex Lie algebras, $\tau_{j}$,
between
$\m'_{j}$ and
$\m''_{j}$, $j=1,\dots,p$, \ste
such that, denoting by  $ \h_{j}:=\{(X, \tau_{j}(X))\vert X
\in
\m'_{j}\}$, and by  $\h$  the fixed point set  of $\sigma$,
one has  :

$$\m={\tilde \m_{0}}\oplus (\oplus_{ j=1,\dots,p} (\m'_{j}\oplus
\m''_{j}))$$
$$\h= {\tilde \h}_{0}\oplus (\oplus_{ j=1,\dots,p} \h_{j}) $$
Notice that an $\R$-linear involutive automorphism of $\m$ is
determined by its fixed point set, as the set of antiinvariant
points is just the orhogonal of the fixed point set,  for the
Killing form of
$\m$, viewed as a real Lie algebra.
The following Theorem generalizes previous results of E.Karolinsky
(cf [K2], Theorem 1 (i) and [K1] for the proof, see also [K3]
Proposition 3.1), as we do not make any restriction on the Manin
form. If we are dealing with $\R$(resp. $\C$)-bilinear Manin form,
Manin triple will mean real (resp. complex) Manin triple
 \ste{\bf Theorem 1}
\ste{\em Let
$B$ be a Manin form and let
$\i$ be a {\bf Lagrangian subalgebra of $\g$ for $B$}, i.e. a real
(resp. complex)  Lie subalgebra of $\g$ whose  real dimension is
equal to the complex dimension of
$\g$ and which is  isotropic for $B$ . Then :
\ste a) If we denote by $\p$ the normalizer in $\g$ of the nilpotent
radical,
$\n$, of $\i$, then $\p$ is a parabolic subalgebra of $\g$, with
nilpotent radical equals to $\n$.\ste
b) Let $\l$ be a Levi subalgebra of $\p$ (i.e. $\l$ is a reductive
Lie subalgebra of $\p$ with $\p=\l\oplus\n$), and denote  by $\m$
its
derived ideal and $\a$ its center, then  the intersection, $\h$, of
$\i$ and $ \m$ is the fixed point set of an af-involution (resp.
f-involution) of
$\m$, which is  isotropic for
$B$. \ste
c) The intersection,
$\i_{\a}$, of
$\a$ and
$\i$, is isotropic for
$B$, and its real dimension equals  the complex dimension of
$\a$.\ste  d) One has
$\i=\h\oplus\i_{\a}\oplus\n$.\ste  Reciprocally, any real (resp.
complex)   Lie subalgebra, $\i$, of $\g$, which is of this form is
Lagrangian  for
$B$. \ste Then,  one
says that {\bf
$\i$ is under $\p$ }}
\ste
One chooses a Cartan subalgebra $\j_{0}$ of $\g$, and a Borel
subalgebra of $\g$ containing $\j_{0}$, $\b_{0}$. Then, from
Theorem 1 and the Bruhat decomposition, one sees (cf. Proposition
1) that  every Manin triple  is conjugated, under
$G$, to a Manin
triple
$(B,\i, \i')$ such that  $\i$ is under $\p$ and  $\i'$ is under
$\p'$, with $\p$ containing $\b_{0}$ and $\p'$ containing the
opposite Borel subalgebra to
$\b_{0}$, with respect to $\j_{0}$. A Manin
triple  satisfying these
conditions will be called {\bf standard}, under $(\p, \p')$.
\ste
If $\r$ is a real subalgebra of $\g$, we denote by $R$ the
analytic subgroup of
$G$ with
 with Lie algebra $\r$ .\ste
If $\e$ is an abelian real subalgebra of $\g$, $\r$ is  an
$\e$-invariant subspace of
$\g$, let
$\Delta(\r,\e)$ the set of weights of $\e$ in $\r$, which is the
subset of $Hom_{\R}(\e,\C)$ whose elements are joint eigenvalues of
elements in $\e$ acting on $\r$. The weight space of
$\alpha$ in
$\Delta(\r,\e)$ is denoted by $\r^{\alpha}$. \ste
Let $\p$, $\p'$ be
given as above,  and let
$B$ be a Manin form on $\g$.\ste \ste
{\bf Theorem 2}  {\em \ste  If there exists a standard
Manin
triple
$(B,\i,
\i')$  under $(\p, \p')$,
 then $\p$ or
$\p'$ is different from
$\g$ }.\ste\ste {\bf Theorem 3} {\em \ste Let  $(B,\i,
\i')$ be  a real (resp. complex) standard  Manin triple under $(\p,
\p')$). Let  $p^{\n'}$  be the projection of $\g$ on the
$\j_0$-invariant supplementary subspace of the nilpotent radical
$\n'$ of $\p'$ , with kernel $\n'$ . Let $\l\oplus \n$  the
Langlands decomposition of   $\p$, such that $\l$ contains
$\j_{0}$. Set
$\i_{1}=p^{\n'}(\tilde{\h}\cap \p')$, where $\tilde{\h}=\i\cap \l$.
One defines similarly $\i'_{1}$. \ste
Then $\i_{1},\i'_{1}$ are contained in $\l\cap \l'$. Moreover, if
$B_{1}$ denotes the restriction of $B$ to $\l\cap \l'$,
$(B_{1},\i_{1},\i'_{1})$ is a real (resp. complex)  Manin triple for
$ \l\cap \l'$. We set $\g_{1}:=\l\cap\l '$.\ste We will use freely
the notation of Theorem 1 for $(B_{1},\i_{1},\i'_{1})$, which is
called the {\bf predecessor  of the standard Manin triple }$(B,\i,
\i')$ .\ste} \ste
{\bf Theorem 4}\ste
{\em Every real (resp. complex )  Manin triple  under
$(\p,\p')$ is conjugate, by an element of  $P\cap P'$, to a
real
 (resp. complex) Manin triple under  $(\p,\p')$, $(B,\i,\i')$,
 whose all successive predecessors,  $(B,\i_{1},\i'_{1}),
(B,\i_{2},\i'_{2}), \dots$, are standard Manin triples
 in  $\g_{1}=\l\cap \l', \g_{2}, \dots$, with
respect to
the intersection of $\b_{0}, \b'_{0}$, with
$\g_{1}, \g_{2}, \dots$, and such that the intersection
$\f_{0}$ ( resp. $\f'_{0}$) of $\j_{0}$ with $\i$ (resp.
$\i'$) is a  fundamental Cartan subalgebra of  $\i$ (resp.
$\i'$), contained in $\i_{1}, \i_{2}, \dots $ (resp.
$\i'_{1}, \i'_{2}, \dots) $ . \ste Such a Manin triple will
be called {\bf strongly standard}. The smallest integer,
 $k$, such that $\g_{k}=\j_{0}$, is called the height of the
strongly standard Manin triple.}
\ste

Now, we assume that $B$ is
$\C$-bilinear.  One defines : $\C^{+}:= \{\lambda \in
\C^{*}\vert Re\> \lambda <0, \>\> or\>\>  Re\> \lambda=0\>\>
et \>\> Im\> \lambda >0\}, \>\> \C^{-}= \C^{*}\setminus
\C^{+}$.  If  $B $  is a complex Manin form  on  $\g$, one denotes
by
 $\g_{+}$ (resp. $\g_{-}$) the sum of the ideals of  $\g$,
$\g_{i}$ , for which the restriction of  $B$ to
$\g_{i}$ is equal to $\lambda_{i}K_{\g_{i}}$, where  $\lambda
_{i}\in \C^{+}$ (resp. $\C^{-}$), and $K_{\g_{i}}$ is the Killing
form
of $\g_{i}$.\ste
Set  $\j_{+}=\j_{0}\cap \g_{+}$, $\j_{-}=\j_{0}\cap \g_{-}$
. The restriction from
$\j_{0}$ to  $\j_{+}$ identifies the roots from  $\j_{0}$ in
$\g_{+}$ to those from  $\j_{+}$ in $\g_{+}$. Let us denote by $
{\tilde
R}_{+}$ the set of these roots and by  $\Sigma_{+}$, the set of
simple roots of the set of positive roots,  ${\tilde R}_{+}^{+}$,
of  ${\tilde R}_{+}$, whose elements are the non zero  weights  of
$\j^{+}$ in
$\b_{0}\cap \g^{+}$. One defines similarly  ${\tilde
R}_{-}$.\ste
 One defines also
$\Sigma_{-}$, the set of simple roots of the set of positive roots,
${\tilde R}_{-}^{+}$, of  ${\tilde R}_{-}$, whose elements are the
non zero  weights  of  $\j_{-}$ in $\b'_{0}\cap
\g^{-}$.\ste
For $\alpha\in{\tilde R}={\tilde R}_{+}\cup {\tilde R}_{-}$, let
$H_{\alpha}\in
\j_{0}$ be the coroot of
$\alpha$. Let ${\cal W}=(H_{\alpha},X_{\alpha},
Y_{\alpha})_{\alpha\in
\Sigma}$, be a Weyl system of generator of $[\g,\g]$, where
$\Sigma=\Sigma_{+}\cup \Sigma_{-}$. More precisely, for all  $
\alpha,
\beta \in \Sigma$, one  has  :
$$ [X_{\alpha}, Y_{\beta}]=\delta_{\alpha\beta} H_{\beta} $$
$$ [H_{\alpha}, X_{\beta}] = N_{\alpha\beta }X_{\beta} $$
$$ [H_{\alpha}, Y_{\beta}]=-N_{\alpha\beta }Y_{\beta}$$
where : $$ N_{\alpha\beta}=\beta(H_{\alpha})=2K_{\g}(H_{\alpha},
H_{\beta})/K_{\g}(H_{\alpha},H_{\alpha})$$
{\bf Definition}\ste
{\em One calls $(A,A',
\i_{\a},\i_{\a'})$ generalized  Belavin-Drinfeld data with respect
to $B$, if :
\ste 1)
 $A$
is a  bijection from a subset  $\Gamma_{+}$ of  $\Sigma_{+}$ on a
subset  $\Gamma_{-}$ of  $\Sigma_{-}$, such that  :  $$
B(H_{A\alpha}, H_{A\beta)}=-B(H_{\alpha},H_{ \beta}), \>\> \alpha,
\beta
\in
\Gamma_{+}
$$
2) $A'$ is a  bijection from a subset $\Gamma'_{+}$ of  $\Sigma_{+}
$
on a
subset $\Gamma'_{-}$ of  $\Sigma_{-}$, such that :
 $$
B(H_{A'\alpha}, H_{A'\beta})=-B(H_{\alpha},H_{ \beta}), \>\> \alpha,
\beta
\in
\Gamma'_{+}
$$
3) Let    $C=''A^{-1}A' \>''$ be the map defined  on $dom
\>C=\{\alpha\in
\Gamma'_{+}\vert A'\alpha\in \Gamma_{-}\}$ by $C\alpha=
A^{-1}A'\alpha$, $\alpha\in dom \>C$. Then
$C$ satisfies  :\ste  For all $\alpha \in
dom\> C$, there exists $n\in \N^{*}$ such that $\alpha, \dots,
C^{n-1}\alpha \in dom\> C\>\>$ and  $\>\> C^{n}\alpha\notin dom\> C
$.
\ste
4)  $\i_{\a}$ (resp.$\i_{\a'}$) is a complex  vector subspace  of
$\j_{0}$, i, included and Lagrangian in the orthogonal, $\a$
(resp.
$\a'$) for the  Killing form of  $\g$ (or for  $B$), to the suspace
generated by $H_{\alpha}$, $\alpha \in
\Gamma:=\Gamma_{+}\cup\Gamma_{-}$ (resp.
$\Gamma':=\Gamma'_{+}\cup\Gamma'_{-}$.\ste
 5) Let
$\f$ be the subspace of
$\j_{0}$ generated by the family
$H_{\alpha}+H_{A\alpha}$, $\alpha\in \Gamma_{+}$. One defines
similarly
 $\f'$. Then  : $$ (\f\oplus\i_{\a} )\cap (\f'\oplus
\i_{\a'})=\{0\}$$  We will denote by $R_{+}$ the sub-system of roots
of
${\tilde R}$ whose elements those of  ${\tilde R}$ which
are linear combination of elements of  $\Gamma_{+}$. One defines
similarly $R_{-}$, $R'_{+}$, $R'_{-}$. We will denote also by
$A$ (resp. $A'$) the  $\R$-linear extension of  $A$ (resp.
$A'$), which defines a bijection from  $R_{+}$ on  $R_{-}$ (resp.
$R'_{+}$ on  $R'_{-}$). }\ste
If  $A$ satisfies the  condition 1) above, there exists a unique
isomorphism $\tau $ from the subalgebra $\m_{+}$ of  $\g$, generated by
 $X_{\alpha},H_{\alpha}, Y_{\alpha}$, $\alpha \in \Gamma_{+}$,
on the subalgebra  $\m_{-}$ of $\g$, generated by
$X_{\alpha},H_{\alpha}, Y_{\alpha}$, $\alpha \in \Gamma_{-}$, such
that :
$$\tau (H_{\alpha})=H_{A\alpha}, \tau (X_{\alpha})=X_{A\alpha},
\tau (Y_{\alpha})=Y_{A\alpha}, \alpha \in
\Gamma^{+}$$
  \ste {\bf Theorem } (cf. Proposition 8 et Th\'eor\`eme 5)\ste
{\em (i) Let  ${\cal BD} =(A,A', \i_{\a}, \i_{\a'})$  be
generalized
Belavin-Drinfeld data, with respect to $B$. Let  $\p$ be
the  parabolic subalgebra  of  $\g$, containing
$\b_{0}$ and  $\m$. Its Langlands decomposition $\p=\l\oplus
\n$, where  $\l$ contains  $\j_{0}$, satisfies $\l=\m\oplus \a$.
Let  $\i$ be equal to $\h\oplus \i_{\a}\oplus \n$, where $
\h:=\{X+\tau
(X)\vert X\in
\m_{+}\}$. One  defines similarly $\i'$. \ste
Then $(B,\i,\i')$ is a strongly standard Manin triple. \ste (ii)
Every Manin triple  is conjugate by an element of $G$ to a
unique  Manin triple of this type. If the original triple is
moreover
strongly standard, the element of $G$ can be taken in $J_{0}$, or,
in
other words, every strongly standard is of the precceeding type, if
one allows to change ${\cal W}$}.
\ste \ste
One shows easily that the preceeding Theorem implies the
classification of certain $R$-matrix given by Belavin and Drinfeld
([BD], Theorem 6.1).\ste One proves also results for
real Manin triples. One retrieves a result of A. Panov [P1] which
classifies certain Lie bialgebras structures on a real simple Lie
algebra.

\ste{\bf Aknowledgment} : I thank very much C. Klimcik for
suggesting me this work, and for many interesting discussions. I
thank  J.L. Brylinski for pointing out to me the work of E.
Karolinsky. I thank also B. Enriquez and Y. Kosmann-Scwarzbach for
telling to me the relations of my earlier work [De] with the work
of A. Belavin and G. Drinfeld [BD], and A. Panov [P1].

\section{Sous-alg\`ebres de Lie  Lagrangiennes}
Dans tout l'article, alg\`ebre de Lie voudra dire alg\`ebre de Lie
de
dimension finie.\ste
 Si $\g$ est une alg\`ebre de Lie on notera souvent  $\g^{der}$ son
id\'eal d\'eriv\'e.\ste  Soit $\a$ est une alg\`ebre de Lie
ab\'elienne sur
$\K=\R$ ou $\C$, $V$ un $\a$-module complexe. Pour $\lambda\in
Hom_{\K}(\a,\C)$, on note
$V^{\lambda}:=\{v\in V\vert Xv=\lambda(X) v, \>\> X\in \a\}$, qui
est appel\'e le sous-espace de poids $\lambda$ de $V$. On dit que
$\lambda $ est un poids de $\a$ dans $V$ si $V^{\lambda}$ est non
nul et on note $\Delta (V,\a)$ l'ensemble des poids non nuls de
$\a$ dans $V$.\ste
Si $G$ est un groupe de Lie, on notera $G^{0}$ sa composante neutre.

\begin{lem}\ste(i)Soit $\g$ une alg\`ebre de Lie semi-simple
complexe,
$\g_{1},\dots, \g_{n}$  ses id\'eaux simples. Toute forme
$\R$ (resp.$\C$)-bilin\'eaire $\g$-invariante sur $\g$ est du type
$B_{\lambda}$ ou
$B_{\lambda}^{\g}$ (resp. $K_{\lambda}$ ou $K_{\lambda}^{\g}$), o\`u
$\lambda=(\lambda_{1},\dots,
\lambda_{n})\in\C^{n}$ et :
$$B_{\lambda}(X_{1}+\dots+X_{n},Y_{1}+\dots+Y_{n})=
\sum_{i=1,\dots,n}Im(\lambda_{i}K_{\g_{i}}(X_{i},Y_{i}))$$
$$K_{\lambda}(X_{1}+\dots+X_{n},Y_{1}+\dots+Y_{n}=
\sum_{i=1,\dots,n}\lambda_{i}K_{\g_{i}}(X_{i},Y_{i}),$$
 si les $X_{i},Y_{i}$ sont des \'el\'ements de $\g_{i}$. Ici
$K_{\g_{i}}$ d\'esigne la forme de Killing de $\g_{i}$.\ste
En particulier, une telle forme est sym\'etrique et les id\'eaux
simples sont deux \`a deux orthogonaux pour une telle forme.\ste
(ii) La forme $B_{\lambda}$ (resp. $K_{\lambda}$)  est non
d\'eg\'en\'er\'ee si et seulement si chacun des $\lambda_{i}$ est
non nul. La forme $B_{\lambda}$ est alors de signature
$(dim_{\C}\g, dim_{\C}\g)$.\ste
(iii) La restriction de $B_{\lambda}^{\g}$ \`a une sous-alg\`ebre
complexe et simple, $\s$, de $\g$, est de la forme $B^{\s}_{\mu}$,
o\`u $\mu=\sum_{i=1,\dots,n}q_{i}\lambda_{i}$ et pour chaque $i$,
$q_{i}$ est un nombre rationnel positif. De plus $q_{i}$ est non
nul si et seulement si $\s$ a un crochet non nul avec $\g_{i}$.
\end{lem}{\em D\'emonstration : }On traite le cas des formes
$\R$-bilin\'eaires, celui des formes $\C$-bilin\'eaires \'etant
semblable. Traitons d'abord le cas o\`u
$\g$ est simple, auquel cas $n=1$. La donn\'ee d'une forme
$\R$-bilin\'eaire $\g$-invariante sur $\g$, \'equivaut \`a celle
d'une
application $\R$-lin\'eaire entre $\g$ et $Hom_{\R}(\g, \R)$, qui
commute \`a l'action de $\g$, regard\'ee comme alg\`ebre de Lie
r\'eelle.
Mais la partie imaginaire de la forme de Killing de $\g$
(regard\'ee
comme complexe) d\'etermine un isomorphisme de $\g$-modules entre
$\g$ et $Hom_{\R}(\g, \R)$. Finalement, la donn\'ee de notre forme
\'equivaut \`a la donn\'ee d'un endomorphisme, $T$, $\R$-lin\'eaire
 de
$\g$, commutant \`a l'action de $\g$.
\ste Soit $\k$ une forme r\'elle compacte de $\g$. On \'ecrit :$$
T(X)=Re\>T(X)+i Im \>T(X),\>\>X\in \k,$$ o\`u $Re\>T(X),Im\>T(X)\in
\k$. Alors
$Re\>T$, $Im \>T$ sont des \'el\'ements de $Hom_{\k}(\k,\k)$. Comme
$\k$ est simple, il r\'esulte du lemme de Schur que cet espace est
\'egal \`a $\R \> Id_{\k}$. Donc, il existe $\lambda\in \C$ tel que
:$$T(X)=\lambda X,\>\> X\in \k$$ Maintenant, si $X,Y\in \k$, on a :
$$T(i[X,Y])=T([iX,Y])=[iX,TY]$$
la derni\`ere \'egalit\'e r\'esultant du fait que $T$ commute \`a
l'action
de $\g$. Joint \`a ce qui pr\'ec\`ede, cela donne :
$$T(i[X,Y])=\lambda i[X,Y], \>\>X,Y\in \k$$
Comme $[\k,\k]=\k$, on conclut que $T$ est la multiplication par
$\lambda$. D'o\`u (i) dans le cas o\`u $\g$ est simple.\ste
Supposons
maintenant que $\g$ soit la somme directe de deux id\'eaux
$\g'$, $\g''$. Soit $B$ une  forme
$\R$-bilin\'eaire $\g$-invariante sur $\g$. On a :
$$B([X',Y'],X'')=-B(Y',[X',X''])=0, \>\>X',\>\> Y'\in\g'
,\>\>X''\in \g'', $$ la derni\`ere \'egalit\'e r\'esultant du fait
que $\g'$,
$\g''$ commutent entre eux. Comme $[\g,\g]=\g$, on a aussi
$[\g',\g']=\g'$. Finalement $\g'$ et $\g''$ sont orthogonaux.
Alors,  on d\'eduit (i) pour $\g$ semi-simple du cas o\`u $\g$ est
simple.\ste (ii) L'assertion sur la non nullit\'e des $\lambda_{i}$
est
claire. Pour l'\'etude de la signature, on se ram\`ene au cas o\`u
$\g$
est simple. Supposons $\lambda\in\C$, non nul. Soit $\k$ une forme
r\'eelle compacte de $\g$. On fixe une base $X_{1},\dots, X_{l}$
de
$\k$. On choisit une racine carr\'ee $\mu$ de
$i\lambda^{-1}$.On pose $Y_{i}=\mu X_{i}$, $Z_{i}=i\mu X_{i}$.
Alors
 $B_{\lambda}(Y_{i},Z_{j})=0$,
$B_{\lambda}(Y_{i},Y_{j})=\delta_{i,j}$,
$B_{\lambda}(Z_{i},Z_{j})=-\delta_{i,j}$. D'o\`u l'on d\'eduit
l'assertion voulue sur la signature.
\ste (iii)
On utilisera le fait suivant :\ste
{\em Si $\rho$ est une repr\'esentation complexe d'une
alg\`ebre de Lie simple complexe $\s$ dans un espace de dimension
finie $V$, on a :$$tr(\rho(X)\rho(Y))=qK_{\s}(X,Y)$$
o\`u $q$ est un nombre rationel positif. De plus $q$ est nul si et
seulement si $\rho$ est triviale.}\ste
L'existence d'un coefficient de proportionnalit\'e $q$ est claire ,
car
la forme de Killing est, \`a un scalaire mutiplicatif pr\`es, la
seule
forme
$\C$-bilin\'eaire invariante sur
$\s$. On se ram\`ene, pour l'\'etude
de $q$, au cas o\`u  $\rho$ est simple. On consid\`ere, sur $V$, un
produit scalaire invariant par une forme r\'eelle compacte, $\k$,
de
$\s$. Si
$X\in
\k$, $\rho(X)$ est antihermitien et :
$$tr(\rho(X) \rho(X))=-tr(\rho(X) \rho(X)^{*})\leq 0$$
cette trace \'etant nulle seulement si $\rho(X)$ est nul. On en
d\'eduit que $q>0$ si $\rho$ est non triviale. Puis on prend un
\'el\'ement
non nul d'une sous-alg\`ebre de Cartan $\j$ de $\s$, sur lequel tous
les poids entiers de $\j$   sont entiers. On en d\'eduit que
$K_{\g}(X,X)$ et
$tr(\rho(X)
\rho(X))$ sont des entiers, le premier nombre \'etant non nul car
\'egal \`a la somme sur toutes les racines, $\alpha$, de
$(\alpha(X))^{2}$. La rationalit\'e de $q$ en r\'esulte.\qed
\begin{defi}\ste Si $\g$ est une alg\`ebre de Lie r\'eductive
complexe,
une forme
$\R$(resp. $\C$)-bilin\'eaire\ste sym\'etrique sur $\g$ et
invariante
par
$\g$ est dite forme de Manin si et seulement si elle est de
signature
$(dim_{\C}\g,dim_{\C}\g)$ (resp. si et seulement si elle est non
d\'eg\'en\'er\'ee). Une forme de Manin est dite forme  sp\'eciale
si sa
restriction
\`a toute sous-alg\`ebre de Lie complexe semi-simple est non
d\'eg\'en\'er\'ee.
\end{defi}

\begin{lem}
\ste (i) Une
forme
$\R$(resp. $\C$)-bilin\'eaire sym\'etrique $\g$-invariante sur $\g$
est \ste sp\'eciale
si et seulement si sa restriction \`a $\g^{der}$ est sp\'eciale et
si
sa restriction au centre, $\z$, de $\g$ est de signature
$(dim_{\C}\z,dim_{\C}\z)$ (resp. est non d\'eg\'en\'er\'ee).\ste
(ii) La restriction d'une forme sp\'eciale \`a une sous-alg\`ebre de Lie
semi-\ste simple complexe de $\g$ est sp\'eciale.\ste
 (iii) La restriction
d'une forme sp\'eciale au centralisateur d'un \'el\'ement
semi-simple de
$\g$,  dont l'image par la repr\'esentation adjointe de $\g$ n'a que
des  valeurs propres  r\'eelles,
 est sp\'eciale.\ste (iv) Si $\g$ est semi-simple, et
$B=B_{\lambda}^{\g}$ (resp. $K=K_{\lambda}^{\g}$), o\`u
$\lambda=(\lambda_{1},\dots,
\lambda_{n})\in\C^{n}$ v\'erifie :
\begin{equation}Si\>\> \sum_{i=1,\dots, n}q_{i}\lambda_{i}=0\>\>
avec\>\> q_{i}\in
\Q^{+},\>\> alors \>\>les \>\>q_{i}\>\> sont\>\> tous\>\>
nuls\end{equation}
alors $B$ (resp. $K$) est sp\'eciale.
On note que (1.1) est satisfait d\`es que les $\lambda_{i}$ sont
ind\'ependants sur $\Q$, ou bien tous strictement positifs. \ste
(v) Si $\g$ est simple, toute forme
$\R$ (resp. $\C$)-bilin\'eaire sym\'etrique $\g$-inva-\ste riante
sur $\g$
est sp\'eciale.
\end{lem}
{\em D\'emonstration :} On traite le cas des formes
$\R$-bilin\'eaires, celui des formes $\C$-bilin\'eaires \'etant
semblable. (i) r\'esulte du fait
que toute sous-alg\`ebre semi-simple de $\g$ est
contenue dans $\g^{der}$ et que, pour toute forme
$\R$-bilin\'eaire sym\'etrique $\g$-invariante sur
$\g$, le centre de $\g$ et $\g^{der}$ sont
orthogonaux.\ste
(ii) est clair. \ste
Montrons  (iii). Comme le centralisateur d'un \'el\'ement de $\g$
est la somme de son intersection avec $\g^{der}$
et de celle avec le centre, on se r\'eduit ais\'ement,
gr\^{a}ce \`a (i) au cas o\`u $\g$ est semi-simple, ce que
l'on suppose dans la suite. Soit $X$ un \'el\'ement semi-simple de
$\g$ tel que $ad\>X$ n'a que des valeurs propres r\'eelles,
soit $\l$ son centralisateur et $\c$ le centre de $\l$. D'apr\`es
(i) et (ii), il suffit de voir que la restriction d'une forme
sp\'eciale \`a $\c$ est de signature $(dim_{\C}\c,dim_{\C}\c)$. Soit
$\j$ une sous-alg\`ebre de Cartan   de $\g$ contenant $X$. Alors
$\c$ est \'egal \`a l'intersection des noyaux des racines de $\j$
dans $\g$ s'annulant sur $X$. Cela montre que $\c$ est la
somme de ses intersections avec les id\'eaux simples de $\g$.
Il suffit alors de prouver notre assertion sur la signature
dans le cas o\`u $\g$ est simple. Alors $B=B_{\lambda}$,avec
$\lambda$ non nul. Soit $\j_{\R}$ l'espace form\'e des \'el\'ements
 de
$\j$ sur lesquels toutes les racines de $\j$ dans $\g$ sont
r\'eelles, qui est une forme r\'eelle de $\j$. Il est clair que
$\c$ est la somme directe de $\c_{\R}:=\c\cap\j_{\R}$ avec
$i\c_{\R}$. On fixe
 une
base orthonorm\'ee,  $X_{1},\dots, X_{l}$,  de $\c_{\R}$, pour la
forme de Killing de $\g$. Celle-ci existe car la forme de Killing
est d\'efinie  positive sur
$\j_{\R}$. On choisit une racine carr\'ee $\mu$ de $i\lambda^{-1}$.
On pose
$Y_{i}=\mu X_{i}$, $Z_{i}=i\mu X_{i}$. Alors
 $B_{\lambda}(Y_{i},Z_{j})=0$,
$B_{\lambda}(Y_{i},Y_{j})=\delta_{i,j}$,
$B_{\lambda}(Z_{i},Z_{j})=-\delta_{i,j}$. D'o\`u l'on d\'eduit
l'assertion voulue sur la signature, ce qui prouve (iii).\ste
(iv) est une cons\'equence imm\'ediate du Lemme 1 et (v) est un cas
particulier de (iv)\qed

\begin{cor}
(i) \ste La restriction  d'une forme
$\R$ (resp.$\C$)-bilin\'eaire, sym\'etrique, $\g$-invariante sur
$\g$, et non d\'eg\'en\'er\'ee,
$B$ , au
 centralisateur, $\l$ d'un \'el\'ement semi-simple de
$\g$,  dont l'image par la repr\'esentation adjointe de $\g$ n'a que
des  valeurs propres  r\'eelles, est non d\'eg\'en\'er\'ee. Il en va de
m\^{e}me de sa restriction \`a $\l^{der}$ et au centre $\a$ de $\l$
.\ste (ii) Si
$B$ est une forme de Manin sur $\g$, sa restriction \`a
$\l$ (resp. $\l^{der}$, $\a$) est une forme de Manin sur $\l$
(resp. $\l^{der}$, $\a$).\ste
(iii) Une forme bilin\'eaire sym\'etrique, $\g$-invariante est une
forme de Manin si et seulement si sa restriction \`a $\g^{der}$ et
sa
restriction au centre de $\g$ sont des formes de Manin.
\end{cor}
\dem On traite le cas des formes
$\R$-bilin\'eaires, celui des formes $\C$-bilin\'eaires \'etant
semblable. Montrons
(i). L'alg\`ebre de Lie $\l$ est la somme du centre $\z$ de $\g$
avec la
somme de ses intersections avec les id\'eaux simples de $\g$. Ces
sous-alg\`ebres sont deux \`a deux orthogonales pour $B$, d'apr\`es
 le
Lemme 1. La restriction de $B$ \`a $\z$ est non d\'eg\'en\'er\'ee,
car
$\g^{der}$ et $\z$ sont orthogonaux. De plus la restriction de
$B$ \`a l'intersection de $\l$ avec un id\'eal simple de $\g$ est
non
d\'eg\'en\'er\'ee, d'apr\`es le Lemme 2 (iii), appliqu\'e \`a cet
id\'eal
simple. D'o\`u l'on d\'eduit que la restriction de $B$ \`a $\l$ est
 non
d\'eg\'en\'er\'ee, ce qui implique le m\^{e}me fait pour sa
restriction \`a
$\l^{der}$ et $\a$.
\ste Si $B$ est de signature $(dim_{\C}\g, dim_{\C}\g) $, sa
restriction \`a $\z$ est de signature $(dim_{\C}\z, dim_{\C}\z) $.
Alors, les assertions sur la signature se d\'emontre comme ci-dessus,
gr\^{a}ce au Lemme 2 (ii), (v), appliqu\'e aux id\'eaux simples de
$\g$. D'o\`u (ii). \ste
La partie si de (iii) est claire. La partie seulement si r\'esulte
de (ii) \qed \ste
On rappelle que le radical  d'une alg\`ebre de Lie,
$\g$, est son plus grand id\'eal r\'esoluble,
et que son radical nilpotent,
  est le plus grand id\'eal, dont les \'el\'ements sont
repr\'esent\'es
par des endomorphismes nilpotents dans chaque repr\'esentation de
dimension finie de $\g$.
 Suivant Bourbaki, on
appelle sous-alg\`ebre de Levi d'une alg\`ebre de Lie, toute
sous-alg\`ebre semi-simple suppl\'e-\ste mentaire du radical. Si
$\g$
est une alg\`ebre de Lie semi-simple complexe on appelle
d\'ecomposition
de Langlands d'une sous-alg\`ebre parabolique $\p$ de $\g$ une
d\'ecomposition de la forme $\p=\l\oplus\n$, o\`u $\n$ est le
radical nilpotent et $\l$ est une sous-alg\`ebre de Lie complexe de
$\g$, r\'eductive dans
$\g$. \ste
Rassemblons dans un Lemme quelques propri\'et\'es des
d\'ecompositions de Langlands d'une sous-alg\`ebre parabolique.
\ste

\begin{lem}
Soit  $\p$ une sous-alg\`ebre parabolique  de $\g$, $\n$ son
radical nilpotent.\ste (i) Si $\j$ est une sous-alg\`ebre de Cartan
de $\p$, c'est une sous-alg\`ebre de Cartan
de $\g$  , et il existe une seule d\'ecomposition de
Langlands  de $\p$,
$\p=\l\oplus\n$, telle que $\j$ soit contenue dans $\l$.
\ste (ii) Si  $\j$, $\j'$  sont  deux sous-alg\`ebres de Cartan de
$\g$, contenues dans $\p$, elles sont conjugu\'ees par un
\'el\'ement de
$P$, qui conjugue les alg\`ebres $\l$ et $\l'$
correspondantes.\ste  (iii) Si $\p=\l\oplus\n$ est une
d\'ecomposition de Langlands de
$\p$, toute sous-alg\`ebre de Cartan
de $\p$ est une  sous-alg\`ebre de Cartan
de $\g$.
\end{lem}
\dem Revenant \`a la d\'efinition des sous-alg\`ebres paraboliques
(cf.
[Bou], Ch. VIII, Paragraphe 3.4, D\'efinition 2, par exemple), on
voit qu'il existe   une sous-alg\`ebre de Cartan
de $\g$, $\j_{1}$, et une d\'ecomposition de Langlands de $\p$,
$\p=\l_{1}\oplus \n$, avec $\j_{1}$ contenue dans $\l_{1}$, telle
que $\l_{1}$ soit la somme des sous-espaces poids de  $\j_{1}$
dans $\p$, qui ne rencontrent pas $\n$. En particulier les poids de
$\j_{1}$ dans $\l_{1}\approx \p/\n$ sont distincts de ceux dans
$\n$. Si
$\p=\l'_{1}\oplus \n$ est une d\'ecomposition de Langlands de $\p$,
avec $\j_{1} $ contenue dans $\l'_{1}$, $\l'_{1}$ est somme des
sous-espaces poids de $\j_{1}$ dans $\l'_{1}\approx \g/\n$. D'o\`u
l'\'egalit\'e de  $\l_{1}$ et $\l'_{1}$. Ceci assure l'unicit\'e de
$\l$ pour $\j=\j_{1}$. Maintenant toutes les sous-alg\`ebres de
Cartan de
$\p$ sont conjugu\'ees \`a $\j_{1}$, par un \'el\'ement de $P$ (cf.
 [Bour], Chapitre VII,
Paragraphe 3.3, Th\'eor\`eme 1). On en d\'eduit (i) par transport de
structure, et
(ii) r\'esulte de la preuve de (i).\ste
Montrons (iii). Si $\p=\l\oplus\n$ est une d\'ecomposition de
Langlands de
$\p$, l'action du centre de $\l$ sur $\g$ est semi-simple, puisque
$\l$ est r\'eductive dans $\g$. Cela implique que, si $\j$ une une
sous-alg\`ebre de Cartan de $\l$, $\j$ est ab\'elienne et  form\'ee
d'\'el\'ements semi-simples de $\g$. Mais $\l$ est isomorphe
\`a $\p/\n$,
qui est une alg\`ebre r\'eductive de m\^{e}me rang que $\g$.
Pour des
raisons de dimension, on voit que $\j$ est une sous-alg\`ebre de
Cartan de $\g$.\qed
 \begin{defi} On appelle af-involution (resp. f-involution),
o\`u a
vaut  pour antilin\'eaire  et f pour flip, d'une alg\`ebre de Lie
semi-simple complexe $\m$, tout automorphisme involutif,
$\R$-lin\'eaire (resp. $\C$-lin\'eaire),
$\sigma$, de $\m$, pour lequel  il existe :\ste
1) un id\'eal ${\tilde  \m_{0}}$ de  $\m$, qui est en outre est
 r\'eduit \`a
z\'ero pour les f-involutions, et une forme r\'eelle,
${\tilde \h}_{0}$, de  ${\tilde  \m_{0}}$ .\ste   2) des id\'eaux
simples de
$\m$, $\m'_{j}$,
$\m''_{j}$, $j=1,\dots,p$.\ste
3) un  isomorphisme d'alg\`ebres de Lie complexes, $\tau_{j}$,
entre
$\m'_{j}$ et
$\m''_{j}$, $j=1,\dots,p$, \ste
tel que, notant  $ \h_{j}:=\{(X, \tau_{j}(X))\vert X
\in
\m'_{j}\}$, et notant  $\h$, l'ensemble des points fixes de
$\sigma$, on ait :

$$\m={\tilde \m_{0}}\oplus (\oplus_{ j=1,\dots,p} (\m'_{j}\oplus
\m''_{j}))$$
$$\h= {\tilde \h}_{0}\oplus (\oplus_{ j=1,\dots,p} \h_{j}) $$
\end{defi}
Il est bon de remarquer qu'un automorphisme involutif de $\m$ est
caract\'eris\'e par son espace de points fixes, car l'espace des
\'el\'ements anti-invariants est juste l'orthogonal de celui-ci,
pour
la forme de Killing de $\m$ regard\'ee comme r\'eelle.\ste
On remarque qu'une f-involution est en particulier une
af-involution.\ste D\'ebutons par quelques propri\'et\'es
\'el\'ementaires.
\begin{lem}
Soit $\h$ une forme r\'eelle simple d'une alg\`ebre de Lie
semi-simple
complexe $\s$. \ste (i) L'alg\`ebre $\s$ n'est pas simple, si et
seulement si
$\h$ admet une structure complexe \ste  (ii) Dans ce cas, $\s$ est
le produit   de deux id\'eaux simples, $\s_{1}$, $\s_{2}$,
isomorphes
\`a
$\h$. \ste (iii) Toujours dans ce cas, il existe un isomorphisme
antilin\'eaire,
$\tau$, entre les alg\`ebres de Lie
$\s_{1}$, $\s_{2}$, regard\'ees comme r\'eelles, tel que
: $$\h=\{(X,\tau (X))\vert X\in \s_{1}\}$$
\end{lem}
\dem Les  points (i) et (ii) sont  bien connus. Montons (iii). Comme
$\h$ est une forme r\'eelle de $\s_{1}\oplus \s_{2}$, la projection
de
$\h$ sur chacun des deux facteurs est non nulle, donc induit un
isomorphisme $\R$-lin\'eaire  de $\h$ avec chacun de ces facteurs.
Il en r\'esulte que $\h$ a la forme indiqu\'ee, mais on sait
seulement
que $\tau $ est $\R$-lin\'eaire. Mais alors $\h$ appara\^{i}t comme
l'ensemble des points fixes de l'automorphisme involutif
$\R$-lin\'eaire de
$\s$ d\'efini par : $$(X,Y)\mapsto (\tau ^{-1}(Y),\tau (X)), X\in
\s_{1}, Y\in \s_{2}
$$
D'apr\`es la remarque qui pr\'ec\`ede le Lemme, cette involution
doit
\^{e}tre \'egale \`a la conjugaison par rapport \`a  $\h$, donc elle est
antilin\'eaire. Ceci implique l'antilin\'earit\'e de $\tau$. \qed

\begin{lem}
On se donne une af-involution, $\sigma$, d'une alg\`ebre de Lie
semi-simple de $\m$. Les id\'eaux simples de $\m$ sont permut\'es
par
$\sigma$. On note
$\m_{j}$,
$j= 1,
\dots,r$, les id\'eaux simples de $\m$. On d\'efinit une involution
$\theta$ de
$\{1,
\dots,r\}$ caract\'eris\'ee par  : $\sigma (\m_{j})=\m_{\theta (j)},
\>\>j =1,
\dots,r$. Pour $j =1,
\dots,r$, l'une des propri\'et\'es suivantes est v\'erifi\'ee :\ste
1)  $\theta (j)=j$ et la restriction de $\sigma $ \`a $\m_{j}$ est
un automorphisme antilin\'eaire de $\m_{j}$. \ste
2) $\theta (j)\not =j$ et la restriction de $\sigma $ \`a $\m_{j}$
est
un isomorphisme antilin\'eaire de $\m_{j}$ sur $\m_{\theta
(j)}$. \ste   3) $\theta (j)\not =j$ et la restriction de $\sigma $
\`a
$\m_{j}$ est un isomorphisme $\C$-lin\'eaire de $\m_{j}$ sur
$\m_{\theta (j)}$.\ste
Si on est dans le cas 1) ou 2), $\m_{j}$ est
contenu dans l'id\'eal ${\tilde \m}_{0}$ de la d\'efinition des
af-involutions. En particulier si $\sigma$ est une f-involution, on
est toujours dans le cas 3).
\end{lem}
\dem
En effet, soit $\h_{p+l}$, $l=1, \dots , q$, les id\'eaux simples de
${\tilde \h_{0}}$. Comme ${\tilde \h_{0}}$ est une forme r\'eelle de
${\tilde \m_{0}}$, ${\tilde \m_{0}}$ est la somme directe des
$\h_{l}+i\h_{l}$, qui sont en outre des id\'eaux. Si
$\h_{p+l}+i\h_{p+l}$ est simple, c'est un id\'eal simple de $\m$ et
on
est dans le cas 1). Sinon $\h_{l}+i\h_{l}$ est le produit de deux
id\'eaux simples et l'on est dans le cas
2), d'apr\`es le Lemme pr\'ec\'edent. \ste On traite de m\^{e}me le
cas o\`u
$\m_{l}$ est
\'egal
\`a l'un des
$\m'_{j}$,
$\m''_{j}$,
$j=1, \dots \>\>, p$, en remarquant que $\h_{j}$ \ste est l'ensemble
des points fixes de l'involution $\C$-lin\'eaire de $(\m'_{j},
\m''_{j})$  donn\'ee par :
\beq(X,Y)\mapsto (\tau ^{-1}(Y),\tau (X)), X\in
\m'_{j}, Y\in \m'' _{j}
\eeq
\qed
\begin{lem}
Tout isomorphisme $\R$-lin\'eaire entre deux alg\`ebres de Lie
simples complexes est soit $\C$-lin\'eaire, soit antilin\'eaire.
\end{lem}
\dem
Deux alg\`ebres de Lie semi-simples complexes qui sont isomorphes
comme alg\`ebres r\'eelles, sont isomorphes comme alg\`ebres de Lie
complexes. En effet, leurs syst\`emes de racines restreintes sont
alors isomorphes. Mais chacun de  ceux-ci est  isomorphe au
syst\`eme
de racine  associ\'e \`a une \sac. D'o\`u l'assertion. Ceci implique
que l'on peut se limiter, pour prouver le Lemme, aux automorphismes
$\R$-lin\'eaire  d'une alg\`ebre de Lie simple complexe, $\g$.
Consid\'erant la conjugaison par rapport \`a une forme r\'eelle de
$\g$, $X\mapsto {\overline X}$, l'alg\`ebre de Lie
$\g':=\{(X,{\overline X)}\vert X\in \g\}$ est une forme r\'eelle de
$\g\times \g$ isomorphe \`a $\g$. Alors tout automorphisme, $\sigma$,
$\R$-lin\'eaire de
$\g$, d\'efinit,  par transport de  structure, un automorphisme de
$\g'$, $\sigma'$, qui poss\`ede un unique prolongement
$\C$-lin\'eaire
\`a
$\g\times \g$, $\sigma''$ . Il existe deux automorphismes
$\C$-lin\'eaires de
$\g$, $\tau $ et $\sigma$, tels que $\sigma''$ v\'erifie :
$$\sigma''(X,Y)=(\tau(X), \theta(Y)) \>\>ou \>\>bien\>\> (\tau(Y),
\theta(X)), (X,Y)\in \g\times \g$$
Ecrivant la d\'efinition de $\sigma'$, la stabilit\'e de
$\g'$, par
$\sigma'' $ implique que $\sigma$ est $\C$-lin\'eaire dans le
premier
cas et antilin\'eaire dans le second.\qed \ste
{\bf Corollaire }
{\em Une involution $\R$ (resp. $\C$)-lin\'eaire d'une alg\`ebre
de Lie semi-simple complexe est
une af-involution (resp. f-involution) si et seulement si sa
restriction \`a tout id\'eal
simple qu'elle laisse invariant est antilin\'eaire (resp . si
 elle ne laisse  aucun id\'eal simple invariant)
}\ste
\dem En effet, d'apr\`es le Lemme 5, il suffit de voir que tout
automorphisme du produit de deux alg\`ebres de Lie simples complexes
$\s_{1}$, $\s_{2}$, permutant les facteurs, est soit $\C$-lin\'eaire,
soit antilin\'eaire. Mais un tel automorphisme est de la forme :
$$(X,Y)\mapsto (\tau ^{-1}(Y),\tau (X)), X\in
\s_{1}, Y\in \s_{2},
$$
o\`u $\tau$ est un isomorphisme $\R$-lin\'eaire entre $\s_{1}$,
$\s_{2}$
. On conclut gr\^{a}ce au Lemme pr\'ec\'edent.\qed
Soit $E$ un espace
vectoriel complexe de dimension finie muni d'une forme
$\R$ (resp. $\C$)-bilin\'eaire sym\'etrique non d\'eg\'en\'er\'ee.
Tout
sous-espace vectoriel r\'eel (resp.  complexe) isotrope est de
dimension r\'eelle inf\'erieure ou \'egale \`a la dimension complexe de $E$
Un sous-espace vectoriel r\'eel (resp. complexe) de $E$, muni d'une
d'une forme
$\R$ (resp. $\C$)-bilin\'eaire sym\'etrique non d\'eg\'en\'er\'ee
 est dit
Lagrangien s'il est isotrope et de dimension r\'eelle \'egale \`a la
dimension complexe de $E$. Un tel espace existe si et seulement la
forme est de signature
$(dim_{\C}E,dim_{\C}E)$ (resp. si $E$ est de dimension complexe
paire).\ste Comme on l'a indiqu\'e dans l'introduction, le
Th\'eor\`eme suivant g\'en\'eralise des r\'esultats d'E. Karolinsky
(cf. [K1], Th\'eor\`eme 3 (i) et [K3] Proposition 3.1)
\begin{theo} Soit $\g$ une alg\`ebre de Lie r\'eductive complexe et
$B$ une forme de Manin $\R$ (resp. $\C$)-bilin\'eaire. Soit $\i$ une
sous-alg\`ebre de Lie r\'eelle (resp. complexe) de
$\g$, Lagrangienne pour
$B$.\ste On a les propri\'et\'es suivantes  :\ste (i) Si l'on note
$\p$
le normalisateur dans $\g$ du radical nilpotent, $\n$, de $\i$,
$\p$
est une sous-alg\`ebre parabolique de $\g$, contenant $\i$,
de radical
nilpotent
$\n$.\ste (ii) Soit $\p=\l\oplus \n$ une d\'ecomposition de
Langlands de $\p$, $\a$ le centre de $\l$ et  $\m$ son id\'eal
d\'eriv\'e. On note   $\h$ l'intersection de $\i$ et $\m$. Elle est
isotrope  pour
$B$.
De plus $\h$ est l'espace des points fixes d'une af-involution
(resp. f-involution), $\sigma$, de
$\m$. \ste Si $B$ est r\'eelle et sp\'eciale, celle-ci est
 antilin\'eaire et $\h$
est une forme r\'eelle de $\m$.
\ste(iii) L'intersection $\i_{\a}$ de $\a$ et $\i$ est
Lagrangienne  pour la restriction de $B$ \`a $\a$.
\ste (iv) On a $\i=\h\oplus\i_{\a}\oplus \n$.\ste\ste
 R\'eciproquement si une sous-alg\`ebre de Lie r\'eelle, $\i$, de
$\g$ est de la forme ci-dessus, elle est Lagrangienne pour $B$.
 On dit alors que {\bf  $\i$ est sous $\p$}.
\end{theo}

{\em D\'ebut de la d\'emonstration du Th\'eor\`eme 1 :}
 Soit $\i$  une sous-alg\`ebre de
Lie r\'eelle (resp. complexe) de $\g$ Lagrangienne pour $B$. On
note
$\r$ son radical et on pose :  \beq
\n:=\{ X\in \r\cap\g^{der}\vert ad_{\g}(X)\>\> est\>\>
nilpotent\}\eeq
Soit  $\h$ une sous-alg\`ebre
de Levi de $\i$.
\begin{lem}
L'ensemble $\n$ est un id\'eal de $\i$ et $[\i,\r]$ est contenu
dans $\n$.
\end{lem}
{\em D\'emonstration :} Montrons que $\n$ est un id\'eal de $\r$
contenant $[\r, \r]$. En effet, comme $\r$ est r\'esoluble, dans
une base, sur $\C$,  bien choisie de $\g$, les $ad_{\g}(X)$,
$X\in \r$ s'\'ecrivent sous forme de matrices triangulaires
sup\'erieures. Pour $X\in \r$, les entr\'ees de la diagonale de
cette matrice sont not\'ees $\lambda_{1}(X),\dots,\lambda_{p}(X)
$, o\`u les $\lambda _{i}$ sont des caract\`eres de $\r$. Alors
$\n$ est l'intersection des noyaux de ces caract\`eres avec
$\g^{der}$. Donc  $\n$ est  un id\'eal de $\r$ contenant $[\r,
\r]$.\ste Si $\f$ est une sous-alg\`ebre de Cartan de $\h$,
$\f\oplus\r$ est encore une alg\`ebre de Lie r\'esoluble car
$[\f,\f]=\{0\}$ et $[\f, \r]\subset \r$. Un argument similaire \`a
celui ci-dessus montre que $[\f, \r]$ est contenu dans $\n$. La
r\'eunion de toutes les sous-alg\`ebres de Cartan de $\h$ est
dense dans $\h$, d'apr\`es la densit\'e des \'el\'ements r\'eguliers
(cf.[Bou], Ch. VII, Paragraphe 2.2 et Paragraphe 2.3,
Th\'eor\`eme 1).
 Par continuit\'e
et densit\'e, on en d\'eduit que $[\h,\r]\subset \n$.\qed

\begin{lem}Soit $k$ un entier compris entre 0 et la dimension
r\'eelle (resp. complexe) de $\r / \n$. Il existe un sous-espace
r\'eel (resp. complexe), ab\'elien,
$\a_{k}$,  de $\r$, de dimension $k$, tel que :\ste
(i) $\a_{k}\cap\n=\{ 0\}.$
\ste(ii) $\a_{k} $ est form\'e d'\'el\'ements semi-simples de $\g$.
\ste (iii) $\a_{k} $ et $\h$ commutent. \end{lem}  {\em
D\'emonstration : } On proc\`ede par r\'ecurrence sur $k$. Si $k=0$,
le Lemme est clair. Supposons le d\'emontr\'e pour
$k<dim_{\K}(\r/\n)$ (o\`u $\K=\R$, resp. $\C$) et montrons le au
rang
$k+1$. Alors
$\h\oplus
\a_{k}$ est une alg\`ebre de Lie r\'eductive dans $\g$, regard\'ee
comme r\'eelle (resp. complexe). D'autre part, comme $\n$ contient
$[\i,
\r]$  d'apr\`es
le Lemme pr\'ec\'edent, on voit que
$\i$  et donc  $\h\oplus \a_{k}$ agit trivialement sur $\r/\n$.
Ceci
implique que : $$[\h\oplus \a_{k},\a_{k}\oplus\n]\subset \n$$
Donc, $ \a_{k}\oplus\n$ est un $(\h\oplus \a_{k})$-sous-module
de $\r$, qui admet un suppl\'ementaire dans $\r$ commutant \`a
$\h\oplus \a_{k}$, puisque
$\h\oplus \a_{k}$ est r\'eductive dans $\g$ et que le quotient
$\r/\a_{k}\oplus\n$ est un $\h\oplus \a_{k}$-module trivial. \ste
On choisit un \'el\'ement non nul de ce suppl\'ementaire, $X$. Alors :
\begin{equation}X\in \r, \>\>X\notin \h\oplus \a_{k}, \>\> et
\>\>[X,\h\oplus \a_{k}]=\{0\} \end{equation} On \'ecrit
$X=X_{s}+X_{n}$, o\`u $X_{n}$ est un \'el\'ement de $\g^{der}$,
$X_{s}$ est un \'el\'ement de $\g$  commutant \`a $X_{n}$
 tels que
$ad_{\g}X_{s}$ est semi-simple et $ad_{\g}X_{n}$est nilpotent.
On sait qu'alors $ad_{\g}X_{s}$, $ad_{\g}X_{n }$ sont des
polyn\^omes en $ad_{\g}X$. Joint \`a (1.2), cela implique :
\begin{equation}[X_{s},\h\oplus
\a_{k}]=\{0\},\>\>[X_{n},\h\oplus
\a_{k}]=\{0\}
\end{equation}
Montrons que $X_{n}$ appartient \`a $\i$. Soit $\j$ une
sous-alg\`ebre
de Cartan de $\h$. Alors $\j\oplus \r$ est r\'esoluble. On peut donc
choisir une base de $\g$ dans laquelle les
$ad_{\g}Y$, $Y\in \j\oplus \r$, sont repr\'esent\'es par des
matrices triangulaires sup\'erieures. On peut choisir cette base
de sorte qu'elle soit la r\'eunion de bases des id\'eaux simples de
$\g$ avec une base du centre de $\g$, ce que l'on fait dans la
suite. Comme  $ad_{\g}X_{n }$ est un polyn\^ome en $ad_{\g}X$, et
que $X\in \r$, il est repr\'esent\'e dans cette base par une
matrice triangulaire sup\'erieure. Comme cet endomorphisme est
nilpotent, sa diagonale est nulle. On en d\'eduit que, pour tout
$Y\in \j\oplus\r$,  les composantes de $X_{n}$ et $Y$ dans les
id\'eaux simples  de $\g$ sont deux \`a deux
orthogonales pour la forme de Killing de $\g$.
Alors, il r\'esulte de l'orthogonalit\'e, pour $B$,  du centre de
$\g$  \`a $\g^{der}$ et du Lemme 1 (i), que :
$$B(X_{n},Y)=0, \>\> Y\in \j\oplus \r$$
En utilisant la densit\'e dans $\h$ de la r\'eunion de ses
sous-alg\`ebres de Cartan, on en d\'eduit que $X_{n}$ est
orthogonal \`a $\i$ pour $B$. Mais $\i$ est un  sous-espace
isotrope pour $B$ de dimension maximale. Donc $X_{n}$ est
\'el\'ement
de $\i$ comme d\'esir\'e.\ste
Ecrivons $X_{n}=H+R$ avec  $H\in \h$, $R\in  \r$. Comme
$X_{n}$ commute \`a $\h$, d'apr\`es (1.3) et que $[\h, \r]\subset
\r$, on voit que $H$ commute \`a $\h$. Donc  $H$ est nul puisque
$\h$ est semi-simple. Finalement $X_{n}\in \r$, et en fait
$X_{n}\in\n$, d'apr\`es la d\'efinition de $\n$. Comme $X$
appartient \`a un suppl\'ementaire de $\a_{k}+\n$ dans $\r$ et que
$X=X_{s}+X_{n}$, on a :
$$X_{s}\in \r, X_{s}\notin \a_{k}+\n$$
On pose $\a_{k+1}=\a_{k}+\K X_{s}$. D'apr\`es (1.3) et la
semi-simplicit\'e de $ad_{\g}X_{s}$, $\a_{k+1}$ v\'erifie les
propri\'et\'es voulues.\qed
\ste{\em Suite de la d\'emonstration du Th\'eor\`eme 1:}\ste On pose
$\i_{\a}:= \a_{p}$, avec $p=dim_{\K}\r/\n$, de sorte que
$\i=\h\oplus \i_{\a}\oplus \n$, o\`u $\i_{\a}$ est form\'e
d'\'el\'ements
semi-simples de $\g$ avec : $$[\h,
\i_{\a}]=\{0\},\>\>\r=\i_{\a}\oplus \n$$ Comme $\i_{\a}\oplus \n$
 est
r\'esoluble, il existe une sous-alg\`ebre de Borel, $\b$, de $\g$,
contenant $\i_{\a}\oplus \n$. A noter que $\n$ est contenue dans le
radical nilpotent, $\v$, de $\b$, d'apr\`es la d\'efinitionde $\n$
et les propri\'et\'es du radical nilpotent d'une sous-alg\`ebre de
Borel. Montrons que
$\i_{\a}$ est contenue dans une sous-alg\`ebre de Cartan de $\g$,
contenue dans
$\b$. En effet, d'apr\`es  [Bor], Proposition 11.15,
la sous-alg\`ebre
de Borel $\b$ contenant $\i_{\a}$, elle contient une sous-alg\`ebre
de Borel du centralisateur $\l$ de $\i_{\a}$ dans
$\g$. Celle-ci contient une sous-alg\`ebre de Cartan $\j$ de
$\l$. Celle-ci est aussi une sous-alg\`ebre de Cartan
de $\g$ contenant $\i_{\a}$ (cf. [Bou], Ch. VII, Paragraphe 2.3,
Proposition 10). \ste Soit $\u$ la somme  des sous
espaces poids de $\i_{\a}$, dans $\v$, , pour des poids non nuls.
Alors
$\p:=\l\oplus\u$ est une sous alg\`ebre parabolique de $\g$,
contenant $\b$. Comme $\l$ est r\'eductive,
$\m:=\l^{der}$ est semi-simple et le radical de $\p$ est \'egal
\`a la
somme du centre $\a$ de $\l$ avec $\u$. La d\'efinition de $\u$
montre que $[\p,\p]=\m\oplus \u$, donc le radical nilpotent de $\p$
est \'egal \`a $\u$ (cf. [Bou], Ch. I, Paragraphe 5.3,
 Th\'eor\`eme 1).
 Comme
 $\i=\h\oplus \i_{\a}\oplus \n$, que $\i_{\a}\oplus \n$ est
contenu dans $\b$ et  que $\h$ est contenu dans $\l$, on a :
$$\i\subset\p$$
Or $\p$ (resp. $\u$) est la somme de ses intersections
$\p_{i}$ (resp. $\u_{i}$) avec les id\'eaux simples $\g_{i}$ de
$\g$. Comme $\p_{i}$ est orthogonal \`a $\u_{i}$ pour la forme de
Killing de $\g_{i}$, on en d\'eduit que $\u$ est orthogonal \`a $\p$
pour $B$ (cf. Lemme 1 (i)). Comme $\i$ est un sous-espace
isotrope pour $B$, de dimension maximale et contenu dans $\p$,
$\u$ est inclus dans $\i$. Il r\'esulte alors de la d\'efinition de
$\n$, que $\u$ est contenu dans $\n$. Par suite, on a  :
\begin{equation}
\i=\u\oplus
(\i\cap\l),\>\>\n=\u\oplus(\n\cap\l)
\end{equation}
 Remarquons que $\a$  contient
$\i_{\a}$.
 On a $\v=\u\oplus
(\v\cap\m)$. Comme $\n\subset \v$ et $\u\subset\n$, on en d\'eduit
que $\n=\u\oplus (\n \cap\m)$. On d\'eduit alors de (1.6) que :
$\n\cap \l=\n\cap\m$. Finalement, on a :
$$\i=\h\oplus \i_{\a}\oplus (\n\cap\m)\oplus \u$$
Alors, posant  :$$\i':=\i\cap\m ,$$ on a :
$$\i'=\h\oplus(\n\cap\m)$$ C' est une sous-alg\`ebre isotrope de
$\m$ pour la restriction de $B$
\`a
$\m$, donc, d'apr\`es le Corollaire du Lemme 2 (ii), de dimension
r\'eelle inf\'erieure ou \'egale \`a la dimension complexe de $\m$.
\ste  De m\^{e}me, $\i_{\a}$ est un sous espace isotrope de $\a$
pour la restriction de $B$ \`a $\a$. D'apr\`es le Corollaire du
Lemme 2, la restriction de $B$ \`a $\a$ est de signature
$(dim_{\C}\a, dim_{\C}\a)$ (resp. est non d\'eg\'en\'er\'ee).
Il en r\'esulte
que la dimension r\'eelle de $\i_{\a}$ est inf\'erieure ou \'egale
\`a
$dim_{\C}\a$.  Mais $dim_{\C}\g=dim_{\C}\m+dim_{\C}\a+dim_{\R}\u$.
Comme
$dim_{\R}\i=dim_{\C}\g$, on d\'eduit de ce qui pr\'ec\`ede que l'on
a
: \begin{equation}dim_{\R} \i'=dim_{\C}\m,\>\>
dim_{\R}\i_{\a}=dim_{\C}\a
\end{equation}
\begin{lem}L'alg\`ebre de Lie  $\n':=\n\cap\m$ est r\'eduite
\`a z\'ero
et $\h$ a la forme indiqu\'ee dans le Th\'eor\`eme.
\end{lem}{\em D\'emonstration :}
Si
$\f$ est une sous-alg\`ebre de Cartan de
$\h$,
$\f+\n'+\i\n'$ est une sous-alg\`ebre de Lie r\'eelle et r\'esoluble de
$\m$. On peut donc choisir une base de $\m$, r\'eunion de bases
des id\'eaux simples de $\m$ , telle que, pour
tout $X\in\f+\n'+\i\n'$, $ad_{\m}X$ soit repr\'esent\'e, dans cette
base, par une matrice triangulaire sup\'erieure. De plus, si $X$
est \'el\'ement de $\n'+\i\n'$, les \'el\'ements diagonaux de cette
matrice sont nulles. On voit, gr\^{a}ce au Lemme 1 (i), que
$\n'+i\n'$ est orthogonal \`a $\f+\n'$, pour la restriction, $B'$,
de
$B$ \`a $\m$. Ceci \'etant vrai pour tout $\f$, $\n'+i\n'$
est  orthogonal \`a
$\i'$ ($=\h+\n'$), pour $B'$. Mais $B'$ est non d\'eg\'en\'er\'ee,
d'apr\`es
le Corollaire du Lemme 2, donc, d'apr\`es le Lemme 1 (ii) et (1.7),
$\i'$ est un sous-espace isotrope de $\m$, pour
$B'$, de dimension maximale. Il en r\'esulte que
$\n'+i\n'$ est contenu dans $\i'$. Mais $\n'+i\n'$ est aussi
  contenu dans $ \v\subset \g^{der}$. Finalement $\n'+i\n'$ est
 contenudans l'intersection  de $\i'$ avec  $\n$, d'apr\`es
la d\'efinition de celui-ci. Mais, comme $\i'$ est contenu dans $\m$,
$\i'\cap\n =\n'$. Alors on a :
$\n'+i\n'\subset
\n'$, c'est \`a dire que $\n'$ est un sous-espace vectoriel
complexe de $\g$, ce qui bien sur \'evident dans le cas complexe.
\ste

Soit
$\h_{j}$,
$j=1,\dots,r$, les id\'eaux simples de $\h$. Comme  $\h_{j}\cap i
\h_{j}$ est un id\'eal de l'alg\`ebre de Lie simple r\'eelle
$\h_{j}$, il y a deux possiblit\'es pour
$\h_{j}$. Ou bien $\h_{j}\cap i \h_{j}=\{0\}$, et alors $\h_{j}+
i \h_{j}$ est une alg\`ebre de Lie semi-simple complexe dont
$\h_{j}$ est une forme r\'eelle. Ou bien $\h_{j}\cap i
\h_{j}=\h_{j}$ et $\h_{j}$ est une sous-alg\`ebre simple complexe
de $\g$. On remarquera que cette deuxi\`eme possibilit\'e est
exclue, si $B$ est sp\'eciale, puique $\h_{j}$ serait alors
semi-simple complexe et isotrope pour $B$.
\ste On suppose que, pour $j=1,\dots , p$,  $\h_{j}\cap i
\h_{j}=\{0\}$, et que pour $j=p+1,\dots,r$, $\h_{j}\cap i
\h_{j}=\h_{j}$. Si $j=1,\dots,p$,  on note $\k_{j}= \h_{j}\oplus i
\h_{j}$. Dans le cas complexe,
i.e. si $B$ est $\C$-bilin\'eaire, $\h_{j}$ est toujours complexe et
$p=0$. Si
$j=p+1,\dots,r$,  on note $\k_{j}$ la somme des projections de
$\h_{j}$ dans les id\'eaux simples de $\m$. On note aussi $\k'_{j}=
\h_{j}+i
\h_{j}$, pour $j=1, \dots, r$. On note
$\k=\sum _{j=1, \dots, r}\k_{j}$ et $\k'=\sum_{j=1,\dots,
r}
\k'_{j}= \h+\i\h$ , qui est contenu dans $\k$.  Par ailleurs,
deux \'el\'ements,
$X$ et
$Y$, de
$\m$ commutent si et seulement $X$ et $iY$ commutent (resp. $X$
commute \`a chacune des projections de $Y$ dans les id\'eaux simples de
$\m$). Il en r\'esulte que, pour $j\not=l$,
$\k_{j}$ commute \`a $\h_{l}$, donc $\k_{j}\cap (\sum_{l\not=
j}\k_{l})$ est contenu dans le centre de $\k_{j}$
, qui est semi-simple complexe.   Il en
r\'esulte que :
\begin{equation}\k=\oplus_{j=1, \dots,r}
\k_{j},\>\> \k'=\oplus_{j=1,
\dots,r}\k'_{j}
\end{equation}
Alors $\k$, $\k'$  sont des   sous-alg\`ebres de Lie semi-simples
complexes de $\m$. \ste
Montrons que : \beq\k\cap \n'=\{0\}\eeq
En effet $\n'$ est un id\'eal dans $\i'=\h+\n'$, puisque
$\i'=\i\cap\m$ et $\n'$ est l'intersection de l'id\'eal $\n$ de $\i$
avec $\m$. C'est donc un $\h$-module, et aussi un $\k'$-module
puisque $\n'$ est un espace vectoriel complexe. Donc
$\k\cap\n'$ est  un sous-$\k'$-module, et aussi une
sous-alg\`ebre r\'esoluble  de $ \k$. Il est clair que les $\k_{j}$
sont des sous-$\k'$-modules de  $\k$, qui n'ont aucun sous-quotient
simple en commun. En effet, d'une part  $\k'_{l}$ agit trivialement
sur
$\k_{j}$, si $j\not=l$. D'autre part, d'apr\`es les d\'efinitions, on
voit que  les sous-quotients simples du
$\k'_{j}$-module $\k_{j}$ sont isomorphes \`a des sous-quotients  de
$\k'_{j}$, dont aucun n'est trivial, puique $\k'_{j}$ est
une alg\`ebre de Lie semi-simple. Donc, si
$\k\cap
\n'$ est non nul, il a une intersection non nulle,
$\k''$, avec l'un des $\k_{j}$, qui est un $\k'_{j}$-sous-module.
 Comme $\k_{j}$ est une sous-alg\`ebre de Lie de $\m$, il en va de
m\^{e}me de $\k''$, qui est de plus r\'esoluble, puisque c'est le cas
de

$\n'$.\ste  Si
$j=1,\dots ,p$,
$\k_{j}=\k'_{j}$ et un $\k'_{j}$-sous-module de $\k_{j}$ est
isomorphe un id\'eal de $\k_{j}$. Alors $\k\cap
\n'$ est \`a la fois semi-simple et r\'esoluble. Une contradiction
qui montre (1.9) dans ce cas.\ste Si $j=p+1,\dots ,r$,
$\k'_{j}=\h_{j}$ est simple, donc l'une des projections de $\k''$
sur un id\'eal simple de $\m$ est isomorphe  \`a $\h_{j}$. Cette
projection
\'etant un morphisme d'alg\`ebres de Lie, il en r\'esulte que l'
alg\`ebre
de Lie r\'esoluble $\k'' $, admet un quotient semi-simple. Une
contradiction qui ach\`eve de prouver (1.9). \ste
Pour $j=p+1, \dots, r$, $\h_{j}$ ne peut \^{e}tre contenu dans un
id\'eal simple de $\m$. En effet, d'apr\`es le Corollaire du Lemme  2,
la restriction de $B$ \`a $\m$ est une forme de Manin. D'apr\`es le
Lemme 1 (ii) et le Lemme 2 (v), la restriction de $B$ \`a un id\'eal
simple de $\m$ est sp\'eciale, et notre assertion en r\'esulte, car
pour $j=p+1, \dots,r$, $\h_{j}$ est isotrope et semi-simple
complexe . Pour $j=p+1, \dots,r$, on notera
$n_{j}$, le nombre d'id\'eaux simples de
$\m$ dans lesquels
$\h_{j}$ a une projection non nulle, et pour $j=1,\dots , p$, on
pose $n_{j}=1$. On vient de voir que :
\beq n_{j}\geq 2, \>\> j=p+1, \dots, r\eeq
Montrons que
$\n'=\{0\}$.  \ste
On a \'evidemment :
$$dim_{\R}\i'=(\sum_{j=1,\dots,r}dim_{\R}\h_{j})+dim_{\R}\n'$$
Alors,  en posant $p_{j}=1$, pour $ j=1,\dots,p$ et $p_{j}=2$, pour
$ j=p+1,\dots,r$, on a :
\beq dim_{\R}\i'=(\sum_{j=1,\dots,r}p_{j}dim_{\C}\k'_{j})
+dim_{\R}\n'
\eeq
Soit $\j_{\k}$ une \sac de $\k$, $\b_{\k}$ une
sous-alg\`ebre de Borel de $\k$, contenant $\j_{\k}$, de radical
nilpotent $\n_{\k}$. Alors $\b_{\k}\oplus \n'$ est une alg\`ebre de
Lie r\'esoluble, contenue dans $\m$, donc contenue dans une
sous-alg\`ebre de Borel, $\b_{\m}$, de $\m$. Alors $\j_{\k}$ est
contenue dans une sous alg\`ebre de Cartan, $\j_{\m}$, contenue dans
$\b_{\m}$ (voir avant (1.6)). De plus $\n_{k}$ est contenu dans le
radical nilpotent de $\b_{\m}$, $\n_{\m}$, qui v\'erifie :
$$\n_{\m}=\{X\in \b_{\m}\vert ad_{\m}X\>\> est \>\>nilpotent\}$$
 En effet, les \'el\'ements de $\n_{\k}$ sont repr\'esent\'es, dans
toute
repr\'esentation de dimension finie de $\n_{\k}$, et donc de
$\b_{\m}$, par des op\'erateurs nilpotents.\ste
De m\^eme, $\n'$ est contenudans $\m$, car pour tout $X\in \n'$,
$ad_{\g}X$, et donc $ad_{\m}X$ est nilpotent.\ste
On note $\j''$ (resp. $\n''$) un suppl\'ementaire de $\j_{\k}$ dans
$\j_{\m}$, resp. $\n_{k}\oplus \n'$ dans $\n_{\m}$. Un calcul
imm\'ediat montre :
\beq
dim_{\C}\m=dim_{\C}\k+2 dim_{\C}\n'+dim_{\C}\j''+2 dim_{\C}\n''
\eeq
En posant $n_{j}=1$ pour $j=1,\dots,p$, on a imm\'ediatement :
\beq
dim_{\C}(\k)=\sum_{j=1,\dots,r}n_{j}dim_{\C}\k'_{j}
\eeq
Alors (1.11), joint \`a la premi\`ere \'egalit\'e de (1.7), et \`a
(1.12)
,
(1.13), implique :
$$
2 dim_{\C}\n'+\sum_{j=1,\dots,r}p_{j}dim_{\C}\k'_{j} =2
dim_{\C}\n'+dim_{\C}\j'' +2dim_{\C}
\n''+\sum_{j=1,\dots,r}n_{j}dim_{\C}\k'_{j}
$$
Comme $n_{j}$ est sup\'erieur ou \'egal \`a $p_{j}$, d'apr\`es (1.10)
et la
d\'efinition des $n_{j}$, $p_{j}$, on en d\'eduit :
$$p_{j}=n_{j}, j=1, \dots, r, \>\> et \>\> \n''=\j''=\{0\}$$
Alors $\k\oplus \n'$ contient la sous-alg\`ebre de Borel, $\b_{\m}$,
de $\m$. C'est une sous-alg\`ebre parabolique dont le radical est
\'egal \`a $\n'$, donc est nilpotent. Elle est donc \'egale \`a $\m$ et
son
radical nilpotent $\n'$ est r\'eduit \`a z\'ero, comme d\'esir\'e.
En outre
$\m=\k$,
 donc les $\k_{j}$ sont des id\'eaux
de $\m$. On pose
${\tilde \m_{0}}=\oplus_{j=1,\dots,p}\k_{j}$,
${\tilde \h_{0}}=\oplus_{j=1,\dots,p}\h_{j}$ . On pose $q= r-p$.
Pour
$l=1,
\dots, q$,
$\k_{p+l}$ est somme de deux id\'eaux simples, $\m'_{l}$, $\m''_{l}$,
car
$n_{p+l}=2$. La projection de $\h_{l+p}$ sur chacun de ces id\'eaux
est bijective, sa surjectivit\'e r\'esultant de la d\'efinition de
$\k_{l}$, son injectivit\'e r\'esultant de la simplicit\'e de
$\h_{l+p}$ et de la non nullit\'e de ce morphisme d'alg\`ebres de
Lie. Donc
$\h_{p+l}:=\{(X,
\tau_{l}(X))\vert X \in \m'_{l}\}$, o\`u $ \tau_{l}$  est un
isomorphisme  $\C$-lin\'eaire de l'alg\`ebre de Lie  $\m'_{l}$ sur
$\m''_{l}$. Donc $\h$  a la forme voulue \qed

\begin {lem}
Aucun poids non nul de $\a$  dans
$\g$ n'est nul sur $\i_{\a}$.\ste

\end{lem}
{\em D\'emonstration :} Raisonnons par l'absurde et supposons
qu'il existe un poids non nul $\alpha$ de $\a$ dans $\g$, nul sur
$\i_{\a}$. Soit
$H_{\alpha}\in \a$ tel que :
\begin{equation}
K_{\g}(H_{\alpha}, X)=\alpha (X),\>\> X\in \a
\end{equation}
Alors
$H_{\alpha}$ appartient \`a l'un des id\'eaux simples de $\g$. En
effet, soit $\j$ une sous-alg\`ebre de Cartan de $\g$, contenant
$\a$. Alors $\alpha $ est la restriction \`a $\a$ d'une racine $\beta$
de $\j$ dans $\g$ et l'on a :
\begin{equation}
K_{\g}( H_{\alpha}, H_{\alpha})>0
\end{equation}
 Soit
$H_{\beta}\in \j$ tel que :
$$K_{\g}(H_{\beta}, X)=\beta (X),\>\> X\in \j$$
Alors $\j$ (resp. $\a$) est la somme directe de ses intersections
avec les id\'eaux simples de $\g$, et
$H_{\beta}$ (resp. $H_{\alpha}$) appartient \`a l'une de celles-ci. On
d\'eduit alors du Lemme 1 (i), qu' il existe $\mu\in \C$ , non nul
car $B$ est non d\'eg\'en\'er\'ee, tel que :
\begin{equation}
B(\lambda H_{\alpha}, X)=Im ( K_{\g}(\lambda \mu
H_{\alpha}, X)),
\>\> \lambda
\in \C \>\> X\in
\g
\end{equation}
Comme $\alpha$ est nulle sur $\i_{\a}$,   il r\'esulte
de (1.14) et (1.16) que $\C H_{\alpha}$ est orthogonale \`a $\i_{a}$,
pour $B$. Comme $B$ est une forme de Manin,  la restriction de $B$ \`a
$\a$ est de signature $(dim_{\C}\a, dim_{\C}\a)$. Tenant compte de
(1.7), on voit que $\i_{\a}$ est un sous-espace de $\a$, isotrope
pour $B$, de dimension maximale. Alors, ce qui pr\'ec\`ede montre que
$\C H_{\alpha}$ est contenu dans $\i_{\a}$. Par ailleurs,  si
$\lambda$ est une racine carr\'ee de $i\mu^{-1}$, $B( \mu H_{\alpha},
\mu H_{\alpha})$ est non nul d'apr\`es (1.15) et (1.16). Une
contradiction avec le fait que $\i_{\a}$ est isotrope qui ach\`eve de
prouver le Lemme. \qed
{\em Fin
de la d\'emonstration du Th\'eor\`eme 1: }\ste
Montrons la propri\'et\'e suivante :
\begin{eqnarray} &&Toute\>\> sous-alg\grave{e}bre\>\>
parabolique\>\>
\q \>\>de\>\>
\g\>\> est\>\> \acute{e}gale \>\>au \>\>normalisateur \nonumber \\&&
dans\>\>
\g \>\>de \>\>son\>\> radical \>\>nilpotent \>\>\w\>\>
\end{eqnarray}   D'abord $\q$ normalise $\w$. Donc,
le normalisateur $\r$ de $\w$ dans $\g$  contient $\q$. C'est
donc une sous-alg\`ebre parabolique de $\g$, qui contient $\w$
comme id\'eal. Compte tenu de [Bou], Ch. I, Paragraphe 5.3, Remarque 2,
$\w$
est contenu dans le radical nilpotent, $\x$, de $\r$. On a alors
: $\p\subset \r,\>\>\w\subset \x$, d'o\`u l'on d\'eduit facilement
que $\x=\w$ et $ \r=\q$, comme d\'esir\'e. Donc $\q$ est bien le
normalisateur de $\n$. Ceci ach\`eve de prouver (1.17).\ste
Montrons que $\n$ est le radical nilpotent de $\i$. En effet, comme
$\n$ est somme de sous-espaces poids sous $\a$ et qu'aucun de ces
poids n'est nul sur $\i_{\a}$, d'apr\`es le Lemme pr\'ec\'edent, on a :
$$[\i_{\a}, \n]= \n$$ Donc $[\i,\i]=\h\oplus \n$ et l'intersection
de $[\i,\i] $ avec le radical $\r=\i_{\a}+\n$ de $\i$ est \'egal \`a
$\n$. Donc, d'apr\`es [Bou] Ch. I, Paragraphe 5.3, Th\'eor\`eme 1, $\n$
est bien le
radical nilpotent de $\i$. On a donc montr\'e que
$\i$ s'\'ecrit de la mani\`ere voulue, pour une d\'ecomposition de
Langlands particuli\`ere  du normalisateur, $\p$, de $\n$ .  Si
$\p=\l'\oplus
\n$ est une autre d\'ecomposition de Langlands de $\p$, $\l$ et $\l'$
sont isomorphes, puisqu'elles sont toutes les deux isomorphes \`a
$\g/\n$. Comme $\i$ contient $\n$, les intersections de $\i$
avec $\l$ et $\l'$ se correspondent dans cet isomorphisme, et la
d\'ecomposition de $\i$ qu'on en d\'eduit, relativement \`a cette
nouvelle d\'ecomposition de Langlands de $\p$, a les propri\'et\'es
voulues.\ste
Etudions la partie r\'eciproque du Th\'eor\`eme. Une sous-alg\`ebre
parabolique de $\g$ est la somme de ses intersections avec
les id\'eaux simples de $\g$. En outre, elle est orthogonale \`a son
radical nilpotent, pour la forme de Killing de $\g$. On conclut
que si $\i$ a une d\'ecomposition comme dans l'\'enonc\'e, elle est
isotrope pour $B$, et de dimension r\'eelle \'egale \`a  la dimension
complexe de $\g$. \qed

\begin{defi}On rappelle qu'une \sac   d'une alg\`ebre
de Lie semi-simple r\'eelle est une \sac fondamentale si et seulement
si elle contient des \'el\'ements r\'eguliers dont l'image par la
repr\'esen-\ste tation adjointe n'a que des valeurs propres imaginaires
pures. Cela \'equivaut au fait qu'aucune racine de cette \sac
n'est r\'eelle.\ste  Une \sac d'une alg\`ebre de Lie r\'eelle est dite
fondamentale si sa projection dans une sous-alg\`ebre de Levi,
parall\`element au radical, est une \sac fondamentale de cette
alg\`ebre de Lie semi-simple r\'eelle.
\end{defi}
Comme toutes les sous-alg\`ebres de Cartan  fondamentales d'une
alg\`ebre de Lie r\'eelle semi-simple sont conjugu\'ees entre elles par
des automorphismes int\'erieurs, il en va de m\^{e}me pour les
sous-alg\`ebres de Cartan fondamentales d'une alg\`ebre de Lie r\'eelle.
En effet il suffit d'adapter la preuve du fait que
 (ii) implique (i), dans [Bou], Ch. VII, Paragraphe 3.5,
 Proposition 5, en
remarquant pour cela que tout automorphisme int\'erieur d'une alg\`ebre
de Levi d'une alg\`ebre de Lie r\'eelle, s'\'etend en un automorphisme
int\'erieur de l'alg\`ebre de Lie.

\begin{lem}
On conserve les hypoth\`eses et notations du Th\'eor\`eme 1. \ste
(i) Si $\f$ est une sous-alg\`ebre de Cartan de $\h$, il existe
des \'el\'ements r\'eguliers de $\g$ contenus dans $\f\oplus
\i_{\a}$. Le centralisateur dans $\g$, $\j$, de ${\tilde
\f}:=\f\oplus
\i_{\a}$,
est une sous-alg\`ebre de Cartan
de $\g$, contenue dans $\l$, v\'erifiant $\j=(\j\cap\m)\oplus
\a$.\ste (ii) Si
$\f$ est une sous-alg\`ebre de Cartan de $\h$,
$ \f\oplus
\i_{\a}$ est une sous-alg\`ebre de Cartan  de $\i$.
\ste (iii) Toute
sous-alg\`ebre de Cartan de $\i$ (resp. sous-alg\`ebre de
Cartan de $\i$ contenue dans $\h+\i_{\a}$) est conjugu\'ee, par un
automorphisme int\'erieur de
$\i$,  (resp. \'egale)   \`a une alg\`ebre de ce type.\ste
(iv) Soit ${\tilde \f'}$ une \sac de $\i$. Il existe une unique
d\'ecomposi-\ste tion de Langlands $\l'+\n$ de $\p$, telle que $\l'$
contienne ${\tilde \f'}$. Alors, notant  $\f':={\tilde
\f'}\cap\l'^{der}$, on a ${\tilde \f'}=\f'+(\i\cap\a')$, o\`u $\a'$
est le centre de $\l'$. De plus ${\tilde \f'}$ est une \sac
fondamentale de $\i$, si et seulement si $\f'$ est une \sac
fondamentale de
$\h':=\i\cap\l'^{der}$
\end{lem} \dem
Montrons (i). On raisonne par l'absurde. On note $\j'$ une \sac de
$\g$, qui contient $\f\oplus \a$. Celle-ci  existe puisque les
\'el\'ements de  $\f\oplus \a$ sont semi-simples. Supposons qu'aucun
\'el\'ement  de
${\tilde \f}:=\f\oplus
\i_{\a}$ ne soit r\'egulier dans
$\g$. Alors, pour tout $X\in {\tilde \f}$, il existe une racine
$\alpha_{X} $ de $\j'$ dans $\g$, nulle sur $X$. Pour une racine
donn\'ee, l'intersection de son noyau avec ${\tilde \f}$ est un ferm\'e
de ${\tilde \f}$. Notre hypoth\`ese montre que ${\tilde \f}$ est la
r\'eunion de ces ferm\'es. Il en r\'esulte que l'un de ces sous-espaces
vectoriels est  d'int\'erieur non vide, donc \'egal \`a ${\tilde \f}$.
Cela signifie qu'une racine de $\j'$ s'annule sur ${\tilde \f}$.
Alors, d'apr\`es le Lemme 7, dont on v\'erifie ais\'ement qu'il  est
valable pour toute d\'ecomposition de Langlands de $\p$ , celle-ci
doit \^{e}tre nulle sur
$\a$. C'est donc une racine de
$\j'$ dans $\m$, qui ne peut \^{e}tre nulle sur $\f$. Une contradiction
qui prouve la premi\`ere partie de (i). Le centralisateur $\j$ de
${\tilde \f}$ est donc une sous-alg\`ebre de Cartan. Par ailleurs
$\j$ contient  $\a$. On en d\'eduit la deuxi\`eme partie de (i).
\ste \ D'apr\`es le Lemme 7,  le nilespace de ${\tilde \f}$ dans $\n$
est r\'eduit \`a z\'ero. Comme $\f$ est une sous-alg\`ebre de
Cartan de $\h$, le nilespace de ${\tilde \f}$ dans $\h$ est
\'egal \`a $\f$. Finalement le nilespace de ${\tilde \f}$ dans
$\i$ est \'egal \`a ${\tilde \f}$. Alors (ii) r\'esulte de [Bou],
Ch. VII,
Paragraphe 2.1, Proposition 3. \ste Montrons (iii). Soit ${\tilde
\f}'$ une
autre sous-alg\`ebre de Cartan de $\i$. La projection, ${\tilde
\f''}$, de
${\tilde
\f}'$ sur
${\tilde \h}:=\h+\i_{\a}$, parall\`element \`a $\n$, est une
sous-alg\`ebre de Cartan de ${\tilde \h}$ (cf. [Bou], Ch. VII,
 Paragraphe 2.1,
Corollaire 2 de la Proposition 4), donc de la forme $\f'+\i_{\a}$,
o\`u $\f'$ est une sous-alg\`ebre de
Cartan de $\h$. Alors ${\tilde \f}'$ et ${\tilde \f}''$ sont deux
sous alg\`ebres de Cartan de $\i$, ayant la m\^{e}me projection sur
${\tilde \h}$, parall\`element \`a $\n$, donc conjugu\'ees par un
automorphisme int\'erieur de $\i$, d'apr\`es [Bou], Ch. VII,
Paragraphe 3.5,  Proposition 5  (voir aussi
apr\`es la D\'efinition 3). Si de plus ${\tilde \f}'$ est contenue dans
$ {\tilde
\h}$, le raisonnement ci-dessus montre qu'elle a la forme indiqu\'ee.
 Ce
qui prouve (iii). \ste
Prouvons (iv). Gr\^{a}ce \`a (iii), on se ram\`ene, par conjugaison,
au  cas o\`u
${\tilde
\f'}$ est contenue dans $\l$ et comme dans (i). Si $\l'+\n$ est une
d\'ecomposition de Langlands  de $\p$, o\`u $\l'$ contient ${\tilde
\f'}$, $\l'$ contient un
\'el\'ement r\'egulier de $\g$, contenu dans ${\tilde
\f'}$, dont le centralisateur dans $\g$ est une \sac de $\g$,
contenue dans $\l'$. Celle-ci est \'egale au centralisateur dans $\g$
de ${\tilde
\f'}$. D'o\`u l'unicit\'e de $\l'$, gr\^{a}ce aux propri\'et\'es des
d\'ecompositions de Langlands (cf. Lemme 3 (i)). L'assertion
sur les \sac fondamentales est claire car $\h'$ est une
sous-alg\`ebre de Levi de $\i$, d'apr\`es le Th\'eor\`eme 1.\qed

\section{Triples de Manin pour une alg\`ebre de Lie r\'eductive
complexe : Descente}
\setcounter{equation}{0}
Dans toute la suite $\g$ d\'esignera une alg\`ebre de Lie
r\'eductive complexe. On fixe,  $\j_{0}$, une sous alg\`ebre de
Cartan de
$\g$, $\b_{0}$ une sous-alg\`ebre de  Borel de $\g$, contenant
$\j_{0}$. On note
$\b'_{0}$ la sous-alg\`ebre de Borel oppos\'ee \`a $\b_{0}$,
relativement \`a $\j_{0}$.
\begin{defi} Un triple de Manin pour $\g$ est un triplet
$(B,
\i,\i')$, o\`u $B$ est une forme de Manin  sur $\g$, $\i$ et $\i'$
sont des sous-alg\`ebres de Lie r\'eelles de $\g$,  isotropes pour $B$,
telles que $\g=:\i\oplus \i' $. La signature de $B$ \'etant \'egale \`a
$(dim_{\C}\g,dim_{\C}\g) $,  $\i$ et
$\i'$ sont Lagrangiennes. Un s-triple  est un triple de Manin
pour lequel la forme est sp\'eciale. \ste  Si $\i$ est sous $\p$ et
$\i'$ est sous
$\p'$, on dit que le triple de Manin est sous $(\p,\p')$

\end{defi}
\begin{rem}   D'apr\`es le Lemme 2 (v), si $\g$ est simple, la notion
de s-triple  et de triple de Manin coincident.
\end{rem}
On note $G$ le groupe  connexe , simplement
connexe, d'alg\`ebre de Lie $\g$. Si $\s$ est une sous-alg\`ebre de
$\g$, on note $S$ le sous-groupe analytique de $G$, d'alg\`ebre de
Lie $\s$. Comme $\g$ est complexe, les sous-groupes
paraboliques de $\g$ sont connexes (cf. [Bor], Th\'eor\'eme
11.16). Donc, si $\p$ est une sous-alg\`ebre parabolique de $\g$,
$P$ est le sous-groupe parabolique de $G$, d'alg\`ebre de Lie
$\p$.  \ste On remarque  que $G$ agit sur l'ensemble des
triples de Manin , en posant, pour tout triple de Manin
$(B,\i,\i')$ et tout
$g\in G$ : $$g (B,\i,\i'):=(B,Ad\>\>g(\i),Ad\>\>g(\i'))$$
  Notre  but est construire, par r\'ecurrence sur la dimension de
$\g^{der}$, tous les triples de Manin modulo cette action de $\g$.
\begin{prop}
Tout triple de Manin est conjugu\'e, sous l'action de $G$, \`a un
triple de Manin
 $(B,\i,\i')$  sous $(\p,\p')$ , avec $\b_{0}\subset \p$ et
$\b'_{0}\subset
\p'$ (un tel triple de Manin
 sera dit standard ).\ste
De plus $\p$ et $\p'$ sont uniques.
\end{prop}
\dem
On rappelle (cf. [Bor], Corollaire 14.13) que :$$
\>\>L\>\>'\>\>intersection
\>\>de\>\> deux \>\>sous-alg\grave{e}bres\>\> de \>\>Borel\>\> de\>\>
\g\>\>,
\>\>{\underline
\b}, \>\>{\underline \b'},$$\beq \>\> contient
\>\> une \>\>sous-alg\grave{e}bre \>\>de\>\> Cartan \>\>de\>\> \g
\eeq\ste
 Soit $(B, {\underline \i}, {\underline \i'}) $
 un
triple de Manin sous
$({\underline \p}, {\underline \p'})$. Soit ${\underline \b}$
(resp. ${\underline \b'}$) une sous-alg\`ebre de Borel de $\g$,
contenue dans
${\underline \p}$ (resp. ${\underline \p'})$. \ste
On a $\g={\underline \i}+{\underline \i'}\subset {\underline
\p}+ {\underline \p'}$. Donc ${\underline
\p}+ {\underline \p'}$ est \'egal \`a $\g$ et ${\underline
P}\>{\underline P'}$ est ouvert dans $G$. Mais ${\underline
P}\>{\underline P'}$ est r\'eunion de $({\underline B}, {\underline
B'})$-doubles classes, qui sont en nombre fini (Bruhat). L'une de
ces doubles classes contenues dans ${\underline
P}\>{\underline
P'}$ doit donc \^{e}tre ouverte. Soit $p\in {\underline
P}$ et $p' \in {\underline
P'}$, tels que ${\underline B}pp'{\underline B'}$ soit un ouvert  de
$G$. On pose
$B_{1}=p^{-1}{\underline B}p$, $B'_{1}=p'{\underline B'}p'^{-1}$.
Alors le sous-groupe de Borel de $G$, $B_{1}$ (resp.
$B'_{1}$), est contenu dans $P$ (resp. $P'$) et $B_{1}B'_{1}$ est
ouvert dans $G$. Donc, on a $\b_{1}+\b'_{1}=\g$ et l'intersection
de
$\b_{1}$ et $ \b'_{1}$ contient une  sous-alg\`ebre de Cartan de
$\g$, $\j_{1}$ (cf. (2.1)). Pour des raisons de dimension,
cette intersection est r\'eduite \`a $\j_{1}$. Alors $\b_{1}$ et $
\b'_{1}$ sont oppos\'ees relativement \`a $\j_{1}$ . D'apr\`es [Bor],
Proposition 11.19, il existe $g'\in G$ tel que
$Ad\>g'(\b_{1})=\b_{0}$, $Ad\>g'(\j_{1})=\j_{0}$. Alors
 $Ad\>g'(\b'_{1})$ est \'egal \`a $\b'_{0}$. Alors, notant
 $\i=Ad\> g'(  {\underline \i})$, $\i'=Ad\> g'(  {\underline \i'})$,
on voit que $(B,\i,\i') $ v\'erifie les propri\'et\'es voulues.\ste
L'unicit\'e de $\p$ r\'esulte du fait que deux sous-alg\`ebres
paraboliques de $\g$, conjugu\'ees par un \'el\'ement de $G$ et contenant
une m\^{e}me sous-alg\`ebre de Borel, sont \'egales (cf. [Bor],
Corollaire 11.17).\qed
\ste
{\em On fixe d\'esormais $\p$ (resp. $\p'$) une sous-alg\`ebre
parabolique de
$\g$, contenant $\b_{0}$ (resp. $\b'_{0})$. On note $\p=\l\oplus \n$
(resp. $\p'=\l'\oplus \n'$) la d\'ecomposition de Langlands de $\p$
(resp. $\p'$) telle que $\l$ (resp. $\l'$) contienne $\j_{0}$ (cf.
Lemme 3).  On note $\m=\l^{der}$, $\a$ le centre de $\l$. Si $\i$
est une sous-alg\`ebre de Lie r\'eelle de $\g$, Lagrangienne  pour une
forme de Manin, on notera $\h=\i\cap \m$, $\i_{\a}=\i\cap\a$,
${\tilde
\h}=\h\oplus\i_{a}$. On introduit des notations similaires
pour $\p'$.}\ste  Comme $\b_{0}\subset \p$
(resp. $\b'_{0}\subset \p' $), $\n$ (resp. $\n'$) est contenu dans
le radical nilpotent de $\b_{0}$ (resp. $\b'_{0}$). Ces derniers
sont d'intersection r\'eduite \`a z\'ero, donc :
\begin{equation}
\n\cap\n'=\{ 0\}
\end{equation}
D\'ecomposant $\p\cap\p'$ en sous-espaces poids sous $\j_{0}$, on
voit que :
\begin{equation}
\p\cap\p'=(\l\cap\l')\oplus (\n\cap\l')\oplus (\n'\cap\l).
\end{equation}

\begin{prop}
(i) Si un \'el\'ement de $G$ conjugue deux triples de Manin  sous $(\p,
\p')$, c'est un \'el\'ement de $P\cap P'$.\ste
(ii) Le groupe $L\cap L'$ est \'egal au sous-groupe analytique de
$G$, d'alg\`ebre de Lie $\l\cap \l'$. Notons $N_{L'}$ (resp. $N'_{L}$)
, le sous-groupe analytique de $G$, d'alg\`ebre de Lie $\n\cap\l'$
(resp. $\n'\cap\l$). Alors on a :
$$P\cap P'=( L\cap L')N_{L'}N'_{L}$$
De plus $N_{L'}$ et $N'_{L}$ commutent entre eux.
\end{prop}
\dem  Si $(B,\i,\i')$ et $(B, {\underline \i}, {\underline \i'})
$ sont deux triples de Manin sous $(\p, \p')$, conjugu\'es par un
\'el\'ement,
$g$, de
$G$, celui-ci conjugue le radical nilpotent de $\i$ avec celui de
${\underline  \i}$, donc normalise $\n$, puisque les deux triples
de Manin sont sous $(\p,\p')$. Mais un \'el\'ement du normalisateur,
$Q$, dans
$G$ de
$\n$, normalise le normalisateur dans $\g$ de $\n$, c'est \`a dire
$\p$, comme on l'a vu  plus haut ( cf.  (1. 17)).
Comme $P$ est connexe, les \'el\'ements de $Q$ normalisent $P$. Donc
$Q$ est inclus dans $P$ et $g\in P$. De m\^{e}me, on a  $g\in P'$. D'o\`u
(i)\ste
Montrons (ii). Il est clair que $P\cap P'$ est un sous-groupe de Lie
de $G$, d'alg\`ebre de Lie $\p\cap\p'$.
On a : $$[\n\cap\l',\n'\cap\l]\subset [\n,\l]\cap
[\n',\l']\subset\n\cap \n'$$
Donc $N_{L'}$ et $N'_{L}$ commutent entre eux, d'apr\`es (2.1). Alors
$ ( L\cap L')^{0}N_{L'}N'_{L}$ est un sous-goupe ouvert et connexe
de $P\cap P'$, donc on a :
\begin{equation}
(P\cap P')^{0}=( L\cap L')^{0}N_{L'}N'_{L}
\end {equation}
Soit $g\in P\cap P'$. Alors $Ad\>g (\j_{0})$ est une sous-alg\`ebre de
Cartan de $\g$, contenue dans $\p\cap\p'$, c'est donc une
sous-alg\`ebre de Cartan de
$\p\cap\p'$ (cf [Bou], Ch. VIII, Paragraphe 2.1, Exemple 3), donc
conjugu\'e,
par un \'el\'ement $g'$ de $(P\cap P')^{0}$, \`a $\j_{0}$, puisqu'il
s'agit d'alg\`ebres de Lie complexes. Donc
$ Ad\>g' g(\j_{0})=\j_{0}$ et $g'g$ est un \'el\'ement de $P\cap P'$. En
utilisant la d\'ecomposition de Bruhat de $G$ et $P$, pour $B_{0}$,
on voit que $g'g$ centralise le centre $\a$ de $\l$. De m\^{e}me
on voit que $g'g$ centralise le centre $\a'$ de $\l'$. Donc $g'g$
est un \'el\'ement du centralisateur, $L''$, de $\a+\a'$ dans $G$. Mais
$P'':=L''(N'_{L} N_{L'})$, est une d\'ecomposition de Langlands du
sous-groupe parabolique de $G$, d'alg\`ebre de Lie :
$$\p''=(\p'\cap\l)\oplus \n= (\l\cap\l')\oplus (\l\cap\n')\oplus
\n$$ Or $P''$ est connexe, puisque $G$ est complexe. Donc $L''$ est
connexe. Par ailleurs, il  contient $L\cap L'$ et a m\^{e}me alg\`ebre
de Lie que $L\cap L'$. Donc on a :
\begin{equation}
 L''=L\cap L'= (L\cap L')^{0}
\end{equation}
On conclut alors que $g'g\in (L\cap L')^{0}$. Donc $g$ est un
\'el\'ement de $ (P\cap P')^{0}$. Ce qui pr\'ec\`ede montre que :

$$P\cap P'=(P\cap P')^{0}$$
On ach\`eve de prouver $(ii)$, gr\^{a}ce \`a (2.4) et (2.5).\qed
Le Lemme suivant est une cons\'equence facile de r\'esultats de
Gantmacher (cf. [G]).
\begin{lem}
Si $\sigma$ et $\sigma'$ sont deux automorphismes involutifs et
antilin\'eaires d'une alg\`e-\ste bre de Lie semi-simple complexe, $\m$,
celle-ci contient au moins un \'el\'ement non nul et  invariant par ces
deux involutions.
\end{lem}
\dem
Avec nos hypoth\`eses
$\sigma\sigma'$ est un automorphisme
$\C$-lin\'eaire de
$\m$,
dont l'espace des points fixes, $\m^{
\sigma\sigma '}$, est un espace vectoriel complexe,  non r\'eduit \`a
z\'ero d'apr\`es  [G], Th\'eor\`eme  28. Mais $ \m^{
\sigma\sigma '}$ , est \'egal \`a $\{ X\in \m\vert
\sigma(X)=\sigma' (X)\}$ donc aussi \'egal \`a $\m^{
\sigma '\sigma }$. Si $X \in \m^{
\sigma\sigma '}$, on a donc $\sigma '(\sigma(X))= X$, soit encore
$\sigma '(\sigma(X))=\sigma (\sigma(X))$. Donc $\sigma(X)$ est
\'el\'ement de $\m^{
\sigma\sigma '}$. Par suite $\sigma$, restreint \`a $\m^{
\sigma\sigma '}$ est une involution antilin\'eaire de $\m^{
\sigma\sigma '}$. L'ensemble de ses points fixes est  une forme
r\'eelle de $\m^{
\sigma\sigma '}$, donc il est non r\'eduit \`a z\'ero. Mais cet ensemble
est \'egal \`a $\m^{\sigma}\cap\m^{\sigma '}$. \qed

\begin{prop}
Si $\sigma$ et $\sigma'$ sont deux af-involutions d'une alg\`ebre
de Lie semi-simple complexe, $\m$, elle contient au moins un \'el\'ement
non nul et  invariant par ces deux involutions.
\end{prop}
\dem
On note
$\m_{j}$,
$j= 1,
\dots,r$, les id\'eaux simples de $\m$. On d\'efinit une involution
$\theta$ de
$\{1,
\dots,r\}$ caract\'eris\'ee par  : $\sigma (\m_{j})=\m_{\theta (j)},
\>\> j =1,
\dots,r$. Nous allons d'abord \'etudier le cas suivant :
\begin{equation}
Il\>\> existe \>\> j \>\>tel\>\> que \>\>\theta(j)=\theta' (j)=j
\end{equation}
Dans ce cas, la restriction de
$\sigma$ et $\sigma'$ \`a $\m_{j}$, sont deux automorphismes
involutifs et
antilin\'eaires de $\m$, d'apr\`es le Corollaire du Lemme 6, qui
ont des points fixes non nuls en commun, d'apr\`es le Lemme
pr\'ec\'edent. La Proposition en r\'esulte, dans ce cas.
\ste Supposons maintenant :
\begin{equation}
Il\>\> existe \>\> j \>\>tel\>\> que \>\>\theta(j)=\theta' (j)
\not=j
\end{equation}
On note $j':=\theta(j)$. Il est clair que :
$(\m_{j}\times \m_{j'})^{\sigma}=\{(X,\sigma(X))\vert X\in \m_{j}\}$
et de m\^{e}me pour $(\m_{j}\times \m_{j'})^{\sigma'}$. Il existe un
\'el\'ement non nul, $X$ de $\m_{j}$ tel que
$({\sigma'}^{-1}{\sigma})(X)=X$, car ${\sigma'}^{-1}{\sigma}$ est
un automorphisme $\R$-lin\'eaire de $\m_{j}$ (cf. [G], Th\'eor\`eme 28).
Alors
$(X,\sigma(X))$ est un  \'el\'ement non nul de $(\m_{j}\times
\m_{j'})^{\sigma}\cap(\m_{j}\times
\m_{j'})^{\sigma'}$ , ce qui prouve la Proposition dans ce cas.
Il nous reste \`a \'etudier le cas suivant :
\begin{equation}
Pour \>\> tout  \>\> j,  \>\>\theta(j)\not=\theta' (j)
\end{equation}
On construit, pour tout $j$, par r\'ecurrence sur $n$, une suite
$j_{1}=j, j_{2},\dots,\ste j_{n},\dots$ , telle que :
\beq
pour\>\> tout\>\> n,\>\> j_{n+1}\not= j_{n}
\eeq
\beq
pour \>\>tout\>\> n,\>\> j_{n+1}= \theta( j_{n}) \>\>ou
\>\>\theta '(
\j_{n})
\eeq
Plus pr\'ecis\'ement, on pose
\beq j_{2}=\theta(1) \>\>si\>\> \theta(1)\not=1,\>\>
j_{2}=\theta'(1) sinon\eeq et, pour
$n\geq 2$, on pose :
$$ j_{n+1}= \theta( j_{n})\>\> si\>\> j_{n}= \theta'(j_{n-1})
\>\>et\>\>\theta (j_{n})\not=j_{n}
$$
$$j_{n+1}= \theta'( j_{n})\>\> si\>\> j_{n}= \theta'(j_{n-1})
\>\>et\>\>\theta (j_{n})=j_{n}$$
\beq
j_{n+1}= \theta'( j_{n})\>\> si\>\> j_{n}= \theta(j_{n-1})
et\>\>\theta' (j_{n})\not=j_{n}
\eeq
$$j_{n+1}= \theta( j_{n})\>\> si\>\> j_{n}= \theta(j_{n-1})
\>\>et\>\>\theta '(j_{n})=j_{n}$$
Ces relations d\'efinissent la suite $(j_{n})$, car,  \`a cause de
(2.8),  on a n\'ecessaire-\ste ment
$\theta(j_{n})\not=\theta'(j_{n})$ et
$\theta(j_{n-1})\not=\theta'(j_{n-1})$. Par ailleurs les relations
(2.9) et (2.10) sont v\'erifi\'ees, la premi\`ere r\'esultant d'une
r\'ecurrence imm\'ediate. On obtient \'egalement les relations
suivantes :
$$
Pour\>\> n\geq 2,  si \>\>\theta (j_{n})\not =j_{n}\>\> et
\>\>si\>\>
\theta '(j_{n})\not =j_{n}
\>\>on \>\>a\>\> :$$
\beq( j_{n-1},j_{n}, j_{n+1})\>\>est\>\> \acute{e} gal\>\>
\grave{a}
\>\>(\theta(j_{n}), j_{n},
\theta'(j_{n}))\>\>ou\>\>\grave{a}\>\>
(\theta'(j_{n}), j_{n},
\theta(j_{n}))
\eeq

\beq
Pour \>\> n\geq 2\>\> et \>\> si \>\>\theta (j_{n})=j_{n}\>\> ou
\>\>si\>\>
\theta '(j_{n})=j_{n}
\>\>on \>\> a \>\> : j_{n-1}=j_{n+1}
\eeq
Apr\`es renum\'erotation des $\g_{j}$, on peut supposer que le d\'ebut
de la suite $(j_{n})$, s'\'ecrit $j_{1}=1,j_{2}=2,\dots,
j_{p}=p,j_{p+1}=k<p$
\ste On fait d'abord la convention  suivante :
$$S\>'\> il\>\> existe \>\>j \>\>tel \>\>que \>\>\theta(j)=j
\>\>ou\>\>
\theta'(j)=j,\>\>on \>\>suppose\>\> qu\>'\>on
$$
$$ \>\>l\>'\>a \>\>choisi \>\>comme\>\>
premier\>\> \acute{e}l\acute{e} ment,\>\> et,\>\>
quitte\>\>\grave{a}
\>\> \acute{e} changer
\>\>le\>\> role $$ \beq \>\>de\>\>
\theta\>\> et\>\>
\theta',\>\> qu'il\>\> est\>\> fix \acute{e}\>\> par\>\>
\theta. \>\> On \>\> a \>alors \>\>\theta(1)= 1,\>\>
\theta'(1)=2\eeq
Traitons le cas o\`u $p=2$. Alors $j_{1}=j_{3}$, et
(2.8), (2.13) montrent que $\theta(2)$ ou $\theta'(2)$ est \'egal \`a 2.
Alors on doit avoir
$\theta(1)=1$, d'apr\`es (2.15),  puis $\theta'(1)=2$ d'apr\`es
(2.11). Comme
$\theta'(1)=2$, on a n\'ecessairement
$\theta(2)=2$.  Dans ce cas, un \'el\'ement
$(X_{1},X_{2})\in \m_{1}\oplus\m_{2}$ est invariant par
$\sigma$ et $\sigma'$ si et seulement si on a :
$$X_{1} =\sigma (X_{1} ) ,\>\>       X_{2}= \sigma
(X_{2}),\>\>
X_{2}=\sigma' ( X_{1})$$
ce qui \'equivaut au syst\`eme :
$$X_{1} =\sigma (X_{1} ),\>\> X_{1} =(\sigma'^{-1}\sigma\sigma
')(X_{1} ),\>\>X_{2}=\sigma' ( X_{1})$$
Mais la restriction de $\sigma $  \`a $\m_{1}$
(resp. $\m_{2}$) est un  automorphisme involutif antilin\'eaire,
puisque
$\theta(1)=1$ et $\theta(2)=2$ (cf. le Corollaire du Lemme 6). De
plus, la restriction de
$\sigma' $  \`a $\m_{1}$ est soit
$\C$-lin\'eaire, soit antilin\'eaire, d'apr\`es le Lemme 6. Alors, la
restriction de $\sigma'^{-1}\sigma\sigma
' $  \`a $\m_{1}$ est un  automorphisme involutif antilin\'eaire.
Alors, dans le cas $p=2$, la Proposition r\'esulte du Lemme 12. \ste
On suppose maintenant :
\beq
p>2
\eeq
On remarque d'abord que :
\beq
Si\>\>j=2,\dots, p-1,\>\> on \>\> a \>\> \theta(j)\not=j \>\> et
\>\>\theta'
(j)\not=j
\eeq
En effet, si on avait par exemple $\theta(j)=j$,  (2.14) conduirait
\`a $j-1=j+1$ une contradiction qui prouve (2.17).\ste
Montrons maintenant que : \beq
k=1 \>\> ou \>\> p-1
\eeq
Supposons  $k\not= 1$. Alors, on a $1<k\leq p-1$. Alors d'apr\`es
(2.13) et (2.14), on a l'\'egalit\'e d'ensembles  :

\beq \{\theta(k),\theta'(k)\}=\{k-1, k+1\},
\eeq ce qui  implique :
\beq
\theta(k)\not =k, \>\> \theta'(k)\not =k
\eeq
Comme $j_{p+1}=k$, on d\'eduit de (2.20) et (2.13) que la s\'equence
$(j_{p}, j_{p+1},j_{p+2})$ est \'egale soit \`a
$(\theta(k),k,\theta'(k))$,  soit \`a  $(\theta'(k),k,\theta(k))$,
c'est \`a dire , gr\^{a}ce \`a (2.19)), soit \`a $(k-1, k, k+1)$,
soit \`a
$(k+1, k, k-1)$. Mais $j_{p}=p$,
 est diff\'erent de $k-1$. Donc $p=k+1$, i.e. $k=p-1$. Ceci
ach\`eve de prouver (2.18).\ste
Traitons d'abord le cas :
\beq
k=1
\eeq
Comme $p>2$, on a $1\not= p-1$. Donc $j_{p-1}=p-1$, est diff\'erent
de $j_{p+1}=1$. Alors, (2.14), (2.13) impliquent l'\'egalit\'e
d'ensembles :
\beq \{\theta(p),\theta'(p)\}=\{p-1, 1\}
\eeq
Supposons, d'abord que $\theta'(1)=2$, ce qui implique,
d'apr\`es (2.11), que $\theta(1)=1$. Comme $p>2$, ni $\theta(p)$, ni
$\theta'(p)$ ne peut \^{e}tre \'egal \`a 1. Une contradiction avec
l'\'equation pr\'ec\'edente qui montre que l'on doit avoir, d'apr\`es
(2.11) :
\beq \theta(1)=2\eeq
Comme $p>2$, la seule possibilit\'e laiss\'ee par
(2.22) est :
\beq\theta (p-1)=p \>\>et\>\> \theta ' (p)=1\eeq
On d\'eduit de (2.23) et (2.17), joints \`a (2.13), que, pour
$j=1,\dots,p-1$, on a :
\beq \theta(j)=j+1, \>\>si \>\>j\>\> est\>\> impair \>\>(resp.\>\>
\theta'(j)=j+1\>\> si\>\> j\>\> est\>\> pair)
\eeq ce qui, joint \`a (2.24), implique que $p$ est pair. Notons
$p=2q$. \ste
On d\'eduit de (2.25) et (2.23) qu'un \'el\'ement $(X_{1}, \dots,X_{p})$
de
$\m_{1}\oplus\dots\oplus
\m_{p}$, est invariant \`a la fois par $\sigma$ et $\sigma'$ si et
seulement si le syst\`eme suivant est v\'erifi\'e:
$$\sigma(X_{1})=X_{2},\>\>
 \sigma'(X_{2})=X_{3}$$
$$\dots,\>\> \dots $$
$$\sigma(X_{2j-1})=X_{2j},\>\>
 \sigma'(X_{2j})=X_{2j+1}$$
$$\dots, \>\>\dots $$
$$\sigma(X_{2q-1})=X_{2q},\>\>
 \sigma'(X_{2q})=X_{1}$$
Notons $\tau $ la restriction de $(\sigma'\sigma)^{q}$ \`a
$\m_{1}$, qui est un automorphisme $\R$-lin\'eaire de $\m_{1}$.  Ce
syst\`eme poss\`ede une solution non nulle si et seulement si
l'\'equa-\ste tion :
$$X_{1}=\tau (X_{1}), \>\> X_{1}\in \m_{1}$$
poss\`ede une solution non nulle. C'est le cas, d'apr\`es [G],
Th\'eor\`eme 28. Ceci ach\`eve de prouver la Proposition dans le cas
$k=1$.\ste On suppose maintenant :
\beq k=p-1>1\eeq
Comme $(j_{p-1,}j_{p},j_{p+1})=(p-1,p,p-1)$, on d\'eduit
de (2.13), (2.14) et (2.8),  que l'on a  soit :
\beq \theta(p)=p,\>\>\theta'(p)=p-1
\eeq
soit :
\beq \theta(p)=p-1,\>\>\theta'(p)=p
\eeq
Alors, d'apr\`es notre convention  (2.15), on a $\theta(1)=1$.
Supposons (2.27) v\'erifi\'e. Comme ci-dessus, ceci joint \`a (2.23) et
(2.13), montre que $p$ est pair et que , pour $j=1,\dots,p-1$, on a
:
\beq \theta'(j)=j+1, \>\>si \>\>j\>\> est\>\> impair \>\>(resp.\>\>
\theta(j)=j+1\>\> si\>\> j\>\> est\>\> pair)
\eeq
On note $p=2q$.
On d\'eduit de (2.28) et (2.29) qu'un \'el\'ement $(X_{1}, \dots,X_{p})$
de
$\m_{1}\oplus\dots\oplus
\m_{p}$, est invariant \`a la fois par $\sigma$ et $\sigma'$ si et
seulement si le syst\`eme suivant est v\'erifi\'e :
$$\sigma(X_{1})=X_{1},\>\>
 \sigma'(X_{1})=X_{2}$$
$$\dots,\>\> \dots $$
$$\sigma(X_{2j})=X_{2j+1},\>\>
 \sigma'(X_{2j+1})=X_{2j+2}$$
$$\dots, \>\>\dots $$
$$\sigma(X_{2q-2})=X_{2q-1},\>\> \sigma'(X_{2q-1})=X_{2q}$$
 $$\sigma(X_{2q})=X_{2q}$$
Notant $\tau $ la restriction de $(\sigma'\sigma)^{q}$ \`a
$\m_{1}$, qui est un automorphisme $\R$-lin\'eai-\ste re de $\m_{1}$,
ce syst\`eme poss\`ede une solution non nulle si et seulement si
le sys-\ste t\`eme :
\beq X_{1}=\sigma (X_{1}),\>\>   X_{1}=(\tau^{-1}\sigma\tau)
(X_{1}),
\>\> X_{1}\in
\m_{1}\eeq  poss\`ede une solution non nulle. La restriction de
$\sigma$ \`a $\m_{1}$ et $\m_{p}$ est antilin\'eaire. Par ailleurs
$\tau$ est soit $\C$-lin\'eaire, soit antilin\'eaire, d'apr\`es le
Lemme 6.  Donc la restriction \`a $\m_{1}$ de  $\tau^{-1}\sigma\tau$
est antilin\'eaire. Il r\'esulte alors du Lemme 12, que (2.30) \`a une
solution non nulle. Ce qui ach\`eve la preuve de la Proposition dans
le cas \'etudi\'e. Le cas o\`u (2.28) est satisfait se traite de
mani\`ere
similaire, mais alors $p$ est impair.\ste
Ceci ach\`eve notre discussion et la preuve de la Proposition.\qed
\begin{theo}
Si $\g$ n'est pas commutative et si  $(B,\i,\i')$ est un triple de
Manin de $\g$, sous $(\p,\p')$,
$\l\cap\l'$ est diff\'erent de $\g$.
\end{theo}
\dem
Raisonnons par l'absurde et supposons qu'il existe un triple de
Manin , $(B,\i,\i')$, sous $(\g,\g)$, et que $\g$ ne soit pas
commutative. Alors $\h:=\i\cap{\g}^{der}$ (resp.
$\h':=\i'\cap\g^{der}$) est  l'espace des points fixes
d'une af-involutions $\sigma$ (resp. $\sigma '$) de
$\g^{der}$, d'apr\`es le Th\'eor\`eme 1. En appliquant la Proposition
pr\'ec\'edente, on aboutit \`a une contradiction avec l'hypoth\'ese
$\i\cap\i'=\{0\}$, ce qui ach\`eve de prouver le Th\'eor\`eme .\qed
\ste  Soit
$V$ un sous-espace
$\j_{0}$ invariant de
$\g$. On suppose  qu'il  est lasomme de sous-espaces poids de $\g$
pour
$\j_{0}$, ce qui s'\'ecrit aussi :
$$V=\sum_{\{ \lambda \in \j _{0}^{*}\vert
V^{\lambda}\not=\{0\}\}}\g^{\lambda}$$ Alors $V$ admet un
unique suppl\'ementaire $\j_{0}$-invariant, $V^{\perp}$, qui est \'egal
 \`a
la somme des sous-espaces poids de $\g$ qui ont une intersection
nulle avec $V$, soit encore :
$$V^{\perp} =\sum_{\{ \lambda \in \j _{0}^{*}\vert \g^{\lambda
}\cap V=\{0\}\}}\g^{\lambda}$$
On note $p_{V}$ (resp. $p^{V}$, la projection de $\g$ sur $V$
(resp. $V^{\perp}$) parall\`element \`a $V^{\perp}$ (resp.  $V$).
Tout sous-espace $\j_{0}$-invariant est stable sous $p^{V}$ et
$p_{V}$. \ste Si de plus $V$ est $\l$-invariant, $V^{\perp}$ est
aussi $\l$-invariant. En effet, comme $\l$ est r\'eductive dans $\g$,
$V$ admet un suppl\'ementaire $\l$-invariant qui n'est autre que
$V^{\perp}$. On voit aussi que dans ce cas, $V^{\perp}$ ne d\'epend
pas du choix de la sous-alg\`ebre de Cartan $\j_{0}$ de $\g$,
contenue dans $\l$.
On a le m\^{e}me fait pour
$\l'$ et
$\l\cap
\l'$.

\begin{theo}
Soit $B$ une forme de Manin r\'eelle (resp. complexe) sur $\g$ et $\i$,
$\i'$ des sous-alg\`ebres de Lie  Lagrangiennes de $\g$ pour $B$,
avec
$\i$ sous $\p$ et $\i'$ sous $\p'$. On a,  gr\^{a}ce au Th\'eor\`eme 1,
$\i=\h\oplus \i_{\a}\oplus \n$, o\`u $\h=\i\cap\m$, $\i_{\a}=\i
\cap\a$. On note ${\tilde \h}=\i\cap\l$. On fait de m\^{e}me pour $\i'$.
\ste Les conditions (i) et (ii) suivantes sont \'equivalentes :
\ste (i) $(B,\i,\i')$ est un triple de Manin
\ste (ii) Notant $\i_{1}=p^{\n'}({\tilde \h }\cap\p')$,
$\i'_{1}=p^{\n}({\tilde \h' }\cap\p)$, on a :\ste
a) $\i_{1}$ et $\i'_{1}$ sont contenues dans $\l\cap\l'$, et
$(B_{1}, \i_{1}, \i'_{1}) $  est un triple de Manin dans $\l\cap
\l'$, o\`u
$B_{1}$ d\'esigne la restriction de $B$ \`a $\l\cap \l'$.\ste
b) $\n\cap \h'$ et $\n'\cap\h$ sont r\'eduits \`a z\'ero.\ste
Si l'une de ces conditions est v\'erifi\'ee, on appellera $(B_{1},
\i_{1}, \i'_{1}) $ l'ant\'ec\'edent du triple de Manin  $(B,\i,\i')$.
\end{theo}
\dem
Montrons que  (i) implique (ii). Supposons que $(B,\i,\i')$ soit un
triple de Manin dans $\g$. Pour des raisons de dimension, ceci
\'equivaut \`a
$\i\cap\i'=\{ 0\}$. Ceci implique imm\'ediatement la propri\'et\'e b) de
(ii).\ste
Montrons ensuite que $\i_{1}$ est une sous-alg\`ebre de Lie
r\'eelle (resp. complexe) de
$\l\cap\l'$, et isotrope pour $B_{1}$. \ste
Etudiant les sous-espaces poids sous $\j_{0}$, on voit que :
\begin{equation}\l\cap \p' =(\l\cap\l')\oplus (\l\cap\n')
\end{equation}
Comme $\l\cap\l'$ est $\j_{0}$-invariant et que $p^{\n'}(\l\cap\n')$
est r\'eduit \`a z\'ero, on a : $$p^{\n'}(\l\cap\p')\subset \l\cap\l'$$
Il en r\'esulte que $\i_{1}$ est bien contenu dans $\l\cap\l'$. Par
ailleurs, la restriction de $p^{\n'}$ \`a $\p'$ est la projection sur
$\l'$, parall\`element \`a $\n'$. C'est donc un morphisme d'alg\`ebres de
Lie, ce qui implique que $\i_{1}$ est une sous-alg\`ebre de Lie
r\'elle de $\l\cap\l'$.\ste
Soit $X$, $X_{1}\in \i_{1}$. Ce sont des \'el\'ements de $\l\cap\l'$, et
il existe $N'$ et $N'_{1}\in \n'$ tels que $Y$ et $Y_{1}$ soient
\'el\'ements de $\i$, o\`u :
$$Y:=X+N',\>\> Y_{1}:=X_{1}+N_{1}'$$
Par ailleurs $\n'$ et $\p'$ sont orthogonaux pour $B$ (cf. la fin
de la d\'emonstration du Th\'eor\`eme 1). Un calcul imm\'ediat montre
alors
que $B(Y,Y_{1})$ est \'egal \`a $B_{1}(X, X_{1})$. Comme $Y, Y_{1}\in
\i$, $B(Y, Y_{1})$ est nul. Finalement, $\i_{1}$ est isotrope pour
$B_{1}$. On montre de m\^{e}me des propri\'et\'es similaires pour
$\i_{1}^{'}$. \ste
Montrons $\i_{1}+\i_{1}' =\l\cap\l'$.
Soit $X\in \l\cap\l'$. Alors $X=I+I'$, avec $I\in \i$, $I'\in \i'$.
Ecrivons $I=H+N,\>\>I'=H'+N'$ o\`u $ H\in {\tilde \h},\>\> H'\in
{\tilde \h',\>\> N\in \n,\>\> N'\in \n'} $. On a donc :
\begin{equation}
X=H+N+H'+N'
\end{equation}
ce qui implique  :
$H=X-H'-N'-N$. On voit ainsi que $H$ est \'el\'ement de
$(\p'+\n)\cap\l$. D\'ecomposant sous l'action de $\j_{0}$, on voit
que : \begin{equation}(\p'+\n)\cap\l=\p'\cap\l\end{equation}
Finalement
$H$ est \'el\'ement de
${\tilde \h}\cap \p'$. de m\^{e}me, on voit que $H'$ est \'el\'ement de
${\tilde \h}'\cap \p$. Par  ailleurs, $\n$ et $\n'$ sont des
sous-espaces $\j_{0}$-invariants et en somme directe avec
$\l\cap\l'$. Donc, appliquant $p_{\l\cap\l'}$ \`a (2.32), on a :
$$X=p_{\l\cap\l'}(H)+p_{\l\cap\l'}(H')$$
De (2.31), on d\'eduit que la restriction de $p_{\l\cap\l'}$ \`a
$\l\cap\p'$ est \'egale \`a la restriction de $p^{\n'}$ \`a $\l\cap\p'$.
Donc, on a : $$p_{\l\cap\l'}(H)= p^{\n'}(H)\in \i_{1}$$
on obtient de m\^{e}me :
$$p_{\l\cap\l'}(H)\in \i'_{1}$$ et l'on conclut que :
$$ X\in
\i_{1}+\i'_{1}
$$
Ceci ach\`eve de prouver que : \beq\l\cap\l'=\i_{1}+\i'_{1}\eeq
Par ailleurs :$$ \l\cap\l'\>\> est
\>\>le\>\> centralisateur\>\> d\>'\>un\>\> \acute{e} l
\acute{e}ment\>\> semi-simple
\>\>de\>\>
\g,\>\> dont$$
$$\>\> l\>'\>image \>par \>\>la\>\> repr \acute{e} sentation\>\>
adjointe\>\> n\>'\>a \>\>que\>\> des\>\> valeurs
\>\>propres\>\>$$\beq
r\acute{e}
 elles\>\>\>\>\>\>\>\>\>\>\>\>\>\>\>\>\>\>\>\>\>\>\>\>\>\>\>
\>\>\>\>\>\>\>\>\>\>\>\>\>\>\>\>\>\>
\>\>\>\>\>\>\>\>\>\>\>\>\>\>\>\>
\>\>\>\>\>\>\>\>\>\>\>\>\>\>\>\>\>\>\>\>
\>\>\>\>\>\>\>\>\>\>\>\>\>\>\>\>\>\>\>\eeq
En effet
$(\l
\cap\l')\oplus ((\n'\cap\l)\oplus \n)$ est une d\'ecomposition de
Langlands d'une sous-alg\`ebre parabolique de $\g$.\ste Alors  la
restriction
$B_{1}$ de
$B$ \`a
$\l\cap\l'$ est une forme de Manin  (cf. Corollaire du Lemme 2
), et
$\i_{1}$,
$\i'_{1}$, qui sont isotropes pour $B_{1}$, sont de dimensions
r\'eelles inf\'erieures ou \'egales \`a la dimension complexe de
$\l\cap\l'$.
La somme dans (2.34) est n\'ecesssairement directe, ce qui ach\`eve de
prouver que (i) implique (ii). \ste
Montrons que (ii) implique (i). Supposons satisfaites  les
conditions a) et b) de (i). Montrons que $\i\cap\i'$ est r\'eduit \`a
z\'ero. Soit $X$ un \'el\'ement de $\i\cap\i'$. Alors :
\begin{equation}X=H+N=H'+N', \>\>o\grave{u} \>\> H\in {\tilde
\h},\>\>  H'\in
{\tilde \h',\>\> N\in \n,\>\> N'\in \n'}
\end{equation}

On a alors :
$$H=H'+N'-N\in\l\cap(\p'+\n)$$ (2.33) implique que $H\in \l\cap\p'$.
De m\^{e}me, on montre que $H'\in \l'\cap\p$. Appliquant $p_{\l\cap\l'}$
\`a (2.36), on voit que :
$$p_{\l\cap\l'}(X)=p_{\l\cap\l'}(H)=p_{\l\cap\l'}(H')$$
et, gr\^{a}ce \`a la premi\`ere partie de la d\'emonstration, cela conduit
 \`a :
$$p_{\l\cap\l'}(X)=p^{\n'}(H)=p^{\n}(H')\in \i_{1}\cap\i'_{1}$$
Donc on a : $$p^{\n'}(H)=p^{\n}(H')= 0$$
Mais $p^{\n'}$ est injective sur ${\tilde \h}\cap \p'$, car
$\n'\cap{\tilde\h}$ est r\'eduit \`a z\'ero. En effet
$\n'\cap{\tilde\h}$ est contenu dans $\n'\cap\l$. On voit que cette
derni\`ere intersection est \'egal \`a $\n'\cap\m$. Donc
$\n'\cap{\tilde\h}$ est \'egal \`a $\n'\cap \h$, qui est r\'eduit
\`a z\'ero,
d'apr\`es b). Donc $H$ est nul et il en va de m\^{e}me de $H'$. Alors $X$
est un \'el\'ement de $\n\cap\n'$, qui est r\'eduit \`a z\'ero, d'apr\`es nos
hypoth\`eses sur $\p$, $\p'$. Donc $X$ est nul et $\i\cap\i'$ est
r\'eduit \`a z\'ero. Alors  la somme $\i+\i'$ est directe, et l'on a
$\g=\i\oplus \i'$ pour des raisons de dimension. Ceci ach\`eve de
prouver le Th\'eor\`eme.\qed

\begin{prop}
Si $(B,\i,\i')$ est un triple de Manin sous $(\p,\p')$, d'ant\'ec\'edent
$(B_{1},\i_{1},\i'_{1})$, et si $g=nn'x\in P\cap P'$, o\`u $x\in L
\cap
L' $, $n\in N_{L'}$, $n'\in N'_{L}$, l'ant\'ec\'edent de
$(B,Ad\>g(\i),Ad\>g(\i'))$ est \'egal \`a $(B_{1},Ad\>
x(\i_{1}),Ad\>
x(\i'_{1}))$.
\end{prop}
\dem
Ecrivons ${\underline \i}=Ad\>g(\i)$ et ${\underline {\tilde
\h}}={\underline \i}\cap\l$, etc.  On note
$(B_{1},{\underline \i}_{1},{\underline \i}'_{1})$, l'ant\'ec\'edent de
$(B,Ad\>g(\i),\ste Ad\>g(\i'))$. On a, gr\^{a}ce \`a la Proposition 2
:
$$g=n'nx=n'x (x^{-1}n x)$$
Donc : $$Ad\>g(\i)= Ad\> n'x(\i)$$
puisque $x^{-1}n x\in N \subset I$. Comme $n'x\in (L\cap
L')N'_{L}\subset L$, cela implique : $${\underline {\tilde
\h}}= Ad\>n'x ({\tilde
\h}),$$ o\`u ${\tilde
\h}=\i\cap \l$. Mais $n'x$ est aussi \'el\'ement de $P'$. Alors, on a :
$${\underline {\tilde
\h}}\cap \p'= Ad\>n'x ({\tilde
\h}\cap \p')$$
D'o\`u l'on d\'eduit :
$${\underline \i}_{1}=p^{\n'}(Ad\>n'x ({\tilde
\h}\cap \p'))$$
Mais il est clair que la restriction de $p^{\n'}$ \`a $\p'$, n'est
autre que la projection sur $\l'$, parall\'element \`a $\n'$. Cette
restriction entrelace l'action adjointe de $P'$ sur $\p'$ avec
l'action naturelle de $P'$ sur $\l'$, identifi\'e au quotient de
$\p'$ par $\n'$ ($N'$ agit trivialement). Il en r\'esulte :
$${\underline \i}_{1}=Ad\> x(p^{\n'} ({\tilde
\h}\cap \p'))=Ad\> x(\i_{1})$$
comme d\'esir\'e. On traite de mani\`ere similaire  ${\underline
\i'}_{1}$.
\qed

\begin{prop}
 Tout triple de Manin sous $(\p,\p')$ est conjugu\'e, par un
\'el\'e-\ste ment de
$P\cap P'$ \`a un triple de Manin , $(B,\i,\i')$, sous $(\p,\p')$,
d'ant\'ec\'edent  $(B_{1},\i_{1},\i'_{1})$ pour lequel
 il existe une \sac fondamentale ${\tilde \f}$ (resp. ${\tilde
\f'}$), de
$\i$ (resp. $\i'$), contenue dans $ \i_{1}$ (resp.
$\i'_{1}$). On dit que le triple $(B,\i,\i')$ est li\'e \`a son
ant\'ec\'edent, avec lien $({\tilde \f},{\tilde \f'})$.
\end{prop}
\dem
D\'emontrons (i).
Soit $(B,{\underline \i},{\underline \i}')$ un triple de Manin pour
$\g$, sous
$(\p,\p')$. On note ${\underline \h}={\underline \i}\cap \m$,
etc.. On note  ${\underline \sigma}$ (resp. ${\underline
\sigma}'$), l'af-involution  de
$\m$ (resp. $\m'$) ayant ${\underline \h}$ (resp. ${\underline \h}'$)
pour espace de points fixes. On d\'efinit de m\^eme ${\underline
\i}_{\a}$ . Comme ${\underline
\i}+{\underline
\i'}=
\g$, on a :
$${\underline \i}+\p'=\g$$ Appliquant $p_{\l}$ \`a cette \'egalit\'e, on
en d\'eduit  :
$$({\underline \i}\cap \l)+(\p'\cap\l)=\l$$
On applique encore la projection de $\l$ sur $\m$, parall\`element \`a
$\a$ pour obtenir :$${\underline \h}+(\p'\cap \m)=\m$$
En cons\'equence, ${\underline H }(P'\cap M)^{0}$ est ouvert dans $M$.
Or
$(P'\cap M)^{0}$ est le sous groupe parabolique de $M$, d'alg\`ebre
de Lie
$\p'\cap\m$. Par ailleurs $\sigma$ \'etant une af-involution, $\m$
est le produit d'id\'eaux $\m_{j}$, invariants par $\sigma$ et sur
lesquels induit :
\ste soit une conjugaison par rapport \`a une forme r\'eelle,\ste  soit
''l' \'echange des facteurs ''  de deux id\'eaux
isomorphes dont $\m_{j}$
est la somme.\ste Il r\'esulte alors de [M2], [M1], que
${\underline \h}\cap\p'$ contient une \sac fondamentale ${\underline
\f}$ de
${\underline \h}$ et une sous-alg\`ebre de Borel, $\b$, de $\m$,
contenant
${\underline
\f}$, contenue dans $\p'\cap\m$, et telle que  :
\begin{equation}
{\underline \sigma}({\underline \b})+{\underline \b}=\m
\end{equation}
D'apr\`es le Lemme 11, :
$${\tilde {\underline \f}}:={\underline \f}+{\underline
\i}_{\a}$$ est une \sac fondamentale de ${\underline \i}$
et  le centralisateur dans $\g$,
${\underline
\j}$, de ${\tilde {\underline \f}}$  est une
\sac  de $\g$, contenue dans
$\l$. De la d\'efinition des af-involutions, il r\'esulte que toute \sac
de
${\underline \h}$ contient des \'el\'ements r\'eguliers de $\m$. Il
r\'esulte alors de [Bor], Proposition 11.15, que ${\underline
\j}\cap\m$ est contenu dans $\b$. Donc
${\underline \j}=({\underline \j}\cap\m)\oplus
\a$ est contenu dans $\p'\cap\l$. C'est une \sac de $\p'\cap\l$,
donc elle est conjugu\'ee \`a $\j_{0}$, par un \'el\'ement du sous-groupe
analytique de $G$, d'alg\`ebre de Lie $\p'\cap\l$. Mais, d'apr\`es
(2.31), on a  : $\p'\cap\l=(\l\cap\l')\oplus (\n'\cap\l)$ et ce
sous-groupe analytique est \'egal \`a $(L\cap L')N'_{L}$, puisque
$L\cap L'$ est connexe, d'apr\`es la Proposition 2.\ste
Donc, il existe $n'\in N'_{L}$, $x\in L\cap L'$, tels que
$$Ad\>xn'({\underline
\j})=\j_{0}$$
soit encore :
\begin{equation} Ad\> n'({\underline \j})=Ad\> x^{-1}(\j_{0})\subset
\l\cap
\l'
\end{equation}
On trouve de m\^{e}me ${\tilde {\underline \f'}}$, $ {\underline
\b'}$, $ {\underline \j'}$ et $x'\in L\cap L'$, $n\in N_{L'}$,
v\'erifiant des propri\'et\'es similaires. On pose :
$$u=nn', \>\> \i=Ad\> u({\underline \i}), \>\> \i'=Ad\> u({\underline
\i'}) $$
Comme $n$ et $n'$ commutent et que $\n$ est un id\'eal de $\i$, on a :
$$Ad\> u({\underline \i})=Ad\> n'({\underline \i})$$
et de m\^{e}me :
$$Ad\> u({\underline \i'})=Ad\> n({\underline \i'})$$
On pose alors :\beq {\tilde \f}= Ad\> n'({\tilde {\underline \f}}),
\>\> {\tilde \f'}= Ad\> n({\tilde {\underline \f}}'), \>\> \b= Ad\>
n'({\underline \b}) , \>\> \b'= Ad\>
n({\underline \b'})
\eeq
On voit alors que $(B,\i,\i')$ est un triple de Manin, conjugu\'e par
$u$ \`a
$(B,{\underline \i},{\underline \i}')$ et   sous
$(\p,\p')$. \ste
On va voir que ${\tilde \f}$ a les propri\'et\'es voulues.
D'abord, comme ${\tilde {\underline \f}}$ est une \sac fondamentale
de
${\underline \i}$,  par conjugaison, on en
d\'eduit que ${\tilde \f}$ est une \sac fondamentale de $\i$. Le
centralisateur,  $\j$, de ${\tilde \f}$ v\'erifie  \beq
\j= Ad\> n'({\underline \j})
\eeq
donc est contenu dans $\l\cap\l'$, d'apr\`es (2.38). Alors, d'apr\`es
le Lemme 11, on a bien
${\tilde \f}=\f\oplus ({\tilde \f}\cap\a)$, o\`u
$\f={\tilde \f}\cap \h$, et ${\tilde \f}\cap\a$ est \'egal \`a
$\i_{\a}$.\ste On a vu que
${\tilde \f}$ est contenu dans $\i\cap\l\cap\l'$, donc dans
${\tilde
\h}\cap\p'$. De plus
$p^{\n'}$ est l'identit\'e sur $\l'$. Donc $ {\tilde \f}$ est contenu
dans $\i_{1}$. Par ailleurs, comme ${\tilde \f}$ est une \sac
de $\i$, contenue dans ${\tilde \h}\cap\p'$, c'est  une \sac de
${\tilde \h}\cap\p'$ (cf. [Bou], Ch. VII, Paragraphe 2.1, Exemple
3), et par projection , c'est une \sac de $\i_{1}$ (cf. l.c.,
Corollaire 2 de la Proposition 4). \ste Il reste \`a voir que cette
\sac de $\i_{1}$ est fondamentale.\ste
On suppose que $\i_{1}$ est sous $\p_{1}$. D'apr\`es le Lemme 11
(iv), il existe une unique d\'ecompo-\ste sition de Langlands
$\p_{1}=
\l_{1}\oplus\n_{1}$, telle que
$ \l_{1}$ contienne ${\tilde \f}$. On note $
\m_{1}={\l_{1}}^{der}$. Il suffit de voir que
${\tilde
\f}\cap
\m_{1}$ est une \sac fondamentale de $
\h_{1}:=\i_{1}\cap \m_{1}$. Pour cela, il suffit de voir qu'aucune
racine de
$ {\tilde \f}\cap
\m_{1}$ dans $
\m_{1}$ n'est r\'eelle. D'apr\`es le Lemme 11 (iv), ${\tilde \f}=
({\tilde \f}\cap
\m_{1})\oplus ({\tilde \f}\cap
\a_{1})$, o\`u $
\a_{1}$ est le centre de $
\l_{1}$. Alors, une racine $\alpha $ de $
{\tilde \f}\cap
\m_{1}$ dans $
\m_{1}$, prolong\'ee par z\'ero sur ${\tilde \f}\cap
\a_{1}$ est une racine de ${\tilde \f}$ dans $\m$. Mais alors,
comme $\f$ est une \sac fondamentale de $\h$, $\alpha $ n'est pas
r\'eelle sur $\f$. Ceci  prouve que ${\tilde \f}$ est une \sac
fondamentale de
$\i_{1}$. On montre de m\^eme que ${\tilde \f'}$ est une \sac
fondamentale de $\i'$ et $\i'_{1}$. Ceci ach\`eve la preuve de la
Proposition.\qed

\begin{theo}
Tout triple de Manin r\'eel (resp. complexe) sous $(\p,\p')$ est
conjugu\'e, par un \'el\'ement de $P\cap P'$, \`a un triple de Manin
r\'eel
 (resp. complexe) sous $(\p,\p')$, $(B,\i,\i')$,
 dont tous les ant\'ec\'edents successifs,  $(B,\i_{1},\i'_{1}),
(B,\i_{2},\i'_{2}),$ $ \dots$, sont des triples de Manin
 standard dans $\g_{1}=\l\cap \l', \g_{2}, \dots$,
relativement \`a l'intersection de $\b_{0}, \b'_{0}$, avec
$\g_{1}, \g_{2}, \dots$, et tel que l'intersection $\f_{0}$ (
resp. $\f'_{0}$) de $\j_{0}$ avec $\i$ (resp. $\i'$) soit une
\sac fondamentale de $\i$ (resp. $\i'$), contenue dans $\i_{1},
\i_{2}, \dots $ (resp.  $\i'_{1},
\i'_{2}, \dots $) . \ste Un triple satisfaisant ces propri\'et\'es
sera appel\'e triple fortement standard.\ste Le plus petit entier,
 $k$, tel que $\g_{k}=\j_{0}$, est appel\'e la hauteur du triple
fortement standard.
\end{theo}
{\em D\'emonstration :} On proc\`ede par r\'ecurrence sur la
dimension de $\g^{der}$. Si celle-ci est nulle, le Th\'eor\`eme est
clair. Supposons l'assertion d\'emontr\'ee pour les alg\`ebres
r\'eductives
dont l'id\'eal d\'eriv\'e est de dimension strictement inf\'erieure
\`a celle
de $\g^{der}$. D'apr\`es la Proposition 5, le triple donn\'e est
conjugu\'e, par un \'el\'ement de $P\cap P'$, \`a un triple de Manin
${\cal T}'$, sous $(\p,\p')$, li\'e \`a son ant\'ec\'edent
${\cal T}'_{1}$.
D'apr\`es l'hypoth\`ese de r\'ecurrence, ce dernier est conjugu\'e par
un
\'el\'ement, $g_{1}$,  de $L\cap L'$, a un triple de Manin
fortement standard : ${\underline{\cal T}}_{1}:=g_{1} {\cal
T}'_{1}$. Par transport de structure, ${\underline{\cal
T}}=g_{1} {\cal T}'$ est li\'e \`a son ant\'ec\'edent
${\underline{\cal T}}_{1}$. On note ${\underline{\cal T}}=(B,
{\underline \i},{\underline \i '}) $,
${\underline{\cal T}}_{1}=(B,
{\underline \i}_{1},{\underline \i '}_{1}) $, et $(\f, \f')$ un
lien entre ces triples. Comme ${\underline{\cal T}}_{1}$ est
fortement standard, $\f_{1}:= \j_{0}\cap {\underline \i}_{1}$
(resp. $\f'_{1}:= \j_{0}\cap {\underline \i}'_{1}$) est une
\sac fondamentale de ${\underline \i}_{1}$ (resp.
${\underline \i '}_{1}$). Alors $\f$ et
$\f_{1}$ (resp. $\f'$ et
$\f'_{1}$ ) sont des \sac fondamentales de ${\underline \i}_{1}$
 (resp. ${\underline \i '}_{1}$), donc conjugu\'ees par un
\'el\'ement $i_{1}$ de ${\underline I}_{1}$ (resp.  $i'_{1}$ de
${\underline I '}_{1}$), i.e. :
$$\f=i_{1} (\f_{1}), \f'=i'_{1}(\f'_{1})$$
On dispose d'une suite exacte de groupes :
 $$0\rightarrow N'\rightarrow P' \rightarrow L'\rightarrow 0 $$
 o\`u la fl\`eche de $P'$ dans $L'$ est le morphisme dont la
diff\'erentielle est la restriction de $p^{\n'}$ \`a $\p$. La
d\'efinition  de ${\underline \i_{1}}$ montre que la restriction de ce
morphisme \`a  $({\underline {\tilde H}}\cap P')^{0}$ est un morphisme
surjectif sur ${\underline I}_{1}$, de noyau contenu dans $N'_{L}$.
Comme
${\underline{\tilde H}}$ est contenu dans ${\underline I}$, on en
d\'eduit qu'il existe
$n'\in N'_{L}$ et $i\in {\underline l}$ tels que  : $$i_{1}=n'i$$
De m\^eme on trouve $n\in N_{L'}$, $i'\in {\underline I'}$ tels que
: $$i'_{1}=ni'$$
Montrons que le triple de Manin ${\cal T}=(B,\i,\i')$,  d\'efini
par  ${\cal T}=n^{-1}n'^{-1} {\underline {\cal T}}$, convient.
D'abord c'est un triple sous $(\p,\p')$, puisque $n,\> n'\in
P\cap P'$, conjugu\'e du triple initial. Par ailleurs $n$ et $n'$
commutent (cf. Proposition 2) et $N$ est contenu dans ${\underline
I}$. Donc, on a : $$
\i=n'^{-1} ({\underline \i})=n'^{-1}i^{-1}({\underline \i})=
i_{1}^{-1}({\underline \i})$$
Comme $\f=i_{1}(\f_{1})$ est une \sac fondamentale de $
{\underline
\i}$, $i_{1}^{-1} (\f)= \f_{1}$ est une \sac fondamentale de
$\i$, par transport de structure. De plus $\f_{1}$ est contenue
dans ${\underline \i}_{1}$. De m\^eme $ \f'_{1}$ est une \sac
fondamentale de $\i'$, contenue dans ${\underline \i}'_{1}$.
Par ailleurs,  comme ${\cal {\underline T}}_{1}$ est fortement
standard, $\j_{0}$ est la somme directe de $\f_{1}$ et
$\f'_{1}$. Comme $\i$ et $\i'$ ont une intersection r\'eduite \`a
z\'ero, il en r\'esulte que $\f_{1}$ (resp. $\f'_{1}$) est \'egal \`a
l'intersection de $\j_{0}$ avec $\i$ (resp. $\i'$). Par
ailleurs, d'apr\`es la Proposition 4, l'ant\'ec\'edent de
$(B,\i,\i')$ est \'egal \`a ${\cal {\underline T}}_{1}$. Comme ce
dernier est fortement standard, ce qui pr\'ec\`ede suffit \`a prouver
que $(B,\i,\i')$ l'est aussi. \qed
\begin{rem}
Le Th\'eor\`eme 4 r\'eduit la
classification des triples de Manin \`a celle des triples fortement
standard.
\end{rem}
\begin{prop}
Si deux triples de Manin fortement standard sont conjugu\'es par
un \'el\'ement de $G$, ils sont de m\^eme hauteur. Ceci permet de
d\'efinir la hauteur d'un triple de Manin  comme la hauteur d'un
triple fortement standard auquel il est conjugu\'e.
\end{prop}
\dem
On proc\`ede par r\'ecurrence sur la dimension de $\g^{der}$. Si
celle-ci est nulle, la Proposition est vraie. Sinon, d'apr\`es la
Proposition 2 et la Proposition 4,  si deux triples de Manin
fortement standard sont conjugu\'es par un \'el\'ement, leurs
ant\'ec\'edents sont conjugu\'es par un \'el\'ement de $L\cap L'$. La
Proposition r\'esulte alors imm\'ediatement de l'application de
l'hypoth\`ese de r\'ecurrence. \qed
Etablissons quelques propri\'et\'es des triples fortement standard. On
utilisera les deux  Lemmes  suivant.
\begin{lem}
Soit $B$ une forme de Manin $\C$-bilin\'eaire. Soit $\i$ une
sous-alg\`ebre Lagrangienne complexe sous $\p$. Soit
$\alpha$ un poids non nul de $\j_{0}$ dans $\g$, tel que
$\g^{\alpha}$ soit contenu dans $\i$. Alors $\g^{\alpha}$ est
contenu dans $\n$.
\end{lem}
{\em D\'emonstration :} On emploie les notations du Th\'eor\`eme 1.
On note $\f=\j_{0}\cap \h$. On sait que si $\beta$ est une
racine de $\j_{0}$ dans $\m$, $\sigma
(\m^{\beta})=\m^{\beta'}$, avec $\beta'\not= \beta$ (cf. Lemme
5, pour les f-involutions). De plus si $\beta_{\vert
\f}=\gamma_{\vert \f}$, pour un autre poids de $\j_{0}$ dans
$\m$, on a $\beta$ \'egal \`a $\gamma$ ou $\gamma'$. On note :
$$\h_{\beta}:=\{X+\sigma
(X)\vert X \in \m^{\beta}\}$$
Soit $R_{*}$ un sous-ensemble de l'ensemble  $R$, des poids non
nuls  de
$\j_{0}$ dans $\m$, tel que
$R_{*}$ et
$\{\beta'\vert \beta \in R_{*}\}$ forme une partition de $R$.
Alors, on a : \beq \h\oplus
\i_{\a}=(\f+\i_{\a}) \oplus (\oplus _{\beta \in R_{*}}\h_{\beta})
\eeq
 qui est une d\'ecomposition en somme directe de
repr\'esentations de $\f$ qui sont deux \`a deux sans
sous-quotients simples isomorphes, toutes \'etant
irr\'educ-\ste tibles, sauf peut-\^etre la premi\`ere. Si $\alpha$ n'est
pas un poids de $\j_{0}$ dans $\n$, comme il est non nul, c'est
un poids de $\j_{0}$ dans $\m$. En \'etudiant l'action de $\f$,
on est conduit \`a :
$$ \g^{\alpha}\subset \h_{\alpha}$$
Comme $\alpha\not= \alpha'$, c'est impossible. Une
contradiction qui ach\`eve de prouver le Lemme. \qed
\begin{lem}
Soit $(B,\i,\i')$ un triple de Manin complexe fortement standard, et
soit
$(B_{1},\i_{1},\i'_{1})$ son ant\'ec\'edent, que l'on suppose sous
$(\p_{1},\p'_{1})$. On note $\p_{1}=\l_{1}\oplus \n_{1}$ la
d\'ecomposition de Langlands de $\p_{1}$ telle que $\l_{1}$ contienne
$\j_{0}$ et on note
$\m_{1}$ l'id\'eal d\'eriv\'e de $\l_{1}$.\ste Alors $\m_{1}$ est \'egal
 \`a
l'id\'eal d\'eriv\'e de $(\m\cap \m')\cap \sigma (\m\cap \m')$ et
l'involution, $\sigma_{1}$, de $\m_{1}$, dont $\h_{1}:=\i_{1}\cap
\m_{1}$ est l'ensemble des points fixes, est \'egale \`a la restriction
de $\sigma$ \`a $\m_{1}$.
\end{lem}
\dem
On r\'eutilise les notations du Lemme pr\'ec\'edent. Etudiant la
d\'ecomposition en repr\'esen-\ste tations irr\'eductibles sous $\f$,  de
${\tilde \h}\cap \p'$, et utilisant (2.41), on voit que :
$${\tilde \h}\cap \p'=(\f+\i_{\a}) \oplus (\oplus _{\beta \in
R_{*}, \h_{\beta}\subset \p'}\h_{\beta})$$
Mais $\h_{\beta}\subset \p'$ si et seulement si $\g^{\beta} $ et
$\g^{\beta  '} $ sont contenus dans $\p'$. Si $\beta$ satisfait
cette condition on a :  $$p^{\n'}(\h_{\beta})=\g^{\beta}\>\> si\>\>
\g^{\beta}\subset \m'\>\> et\>\> \g^{\beta}\subset \n'$$
$$p^{\n'}(\h_{\beta})=\h_{\beta}\>\> si\>\>
\g^{\beta},  \g^{\beta '}\subset \m'$$
Notons, pour $V$ sous-espace vectoriel complexe de $\g$, invariant
sous $\j_{0}$, \ste $\Delta (V, \j_{0})$, l'ensemble  des poids non nuls
de
$\j_{0}$ dans $V$. Alors on a :

\beq
\i_{1}=\u_{1}\oplus \v_{1}\oplus \i_{\a}
\eeq
o\`u :
\beq
\u_{1}=\f\oplus_{\beta\in R_{*}\cap \Delta (\l\cap \l',\j_{0})
,\>\beta '\in R_{*}\cap \Delta (\l\cap \l',\j_{0})}\h_{\beta}
\eeq
et \beq
\v_{1}=\oplus_{\beta\in \Delta (\l\cap \l',\j_{0}), \>
\beta ' \in \Delta (\l\cap \n',\j_{0})}\g^{\beta}
\eeq
On remarque que :
\beq
\u_{1}= ((\m\cap \l')\cap \sigma (\m\cap\l'))^{\sigma}
\eeq
Donc $\u_{1} $ est une sous-alg\`ebre de Lie de $\m\cap\l'$,
r\'eductive, dont le centre est contenu dans $\f$.
Notons $\w_{1}=(\m\cap \m')\cap \sigma (\m\cap\m')+\j_{0}$.
On voit que $\q_{1}:=\w_{1}\oplus \v_{1}$ est une
une sous-alg\`ebre parabolique, de $\l\cap\l'$,
dont le radical nilpotent est $\v_{1}$, et la d\'ecomposition
ci-dessus est une d\'ecomposition de Langlands de $\q_{1}$ avec
$\j_{0}$ contenu dans $\w_{1}$.
\ste Par ailleurs
$\i_{1}$ est contenue dans $\q_{1}$, d'apr\`es (2.42), (2.43) et
(2.44), Tenant compte du fait que, pour $\beta\in R$,
$[\f,\m^{\beta}]=\m^{\beta}$, d'apr\`es la d\'efinition de $\f$ et des
f-involutions, on d\'eduit de (1.3) que le radical  nilpotent de
$\i_{1}$ contient $\v_{1}$. Par ailleurs, comme le centre de
$\u_{1}$ est contenu dans $\j_{0}$, le radical et, a fortiori, le
radical nilpotent  de
$\i_{1}$ est contenu dans $\j_{0}\oplus \v_{1}$. Mais ce radical
nilpotent ne rencontre pas $\j_{0}$ (cf. Th\'eor\`eme 1), donc il est
\'egal \`a
$\v_{1}$. Alors $\p_{1}=\q_{1}$, $\l_{1}=\w_{1}$. Donc $\m_{1}$ a
la forme annonc\'ee. Comme $\h_{1}=\u_{1}\cap \m_{1}$, on d\'eduit de
(2.45) que $\sigma_{1}$ a la forme annonc\'ee. \qed

\section{Classification des triples de Manin complexes}

\setcounter{equation}{0} Dans toute cette partie triple de Manin
voudra dire triple de Manin complexe. On rappelle qu'on a fix\'e
$\j_{0}$ et
$\b_{0}$, $\p$, $\p'$. \ste
   On d\'efinit : $$\C^{+}:= \{\lambda \in \C^{*}\vert
Re\> \lambda <0, \>\> ou\>\>  Re\> \lambda=0\>\> et \>\> Im\>
\lambda >0\}, \>\> \C^{-}= \C^{*}\setminus \C^{+}$$  Si $B $ est
une forme de Manin complexe sur $\g$, on note $\g_{+}$ (resp.
$\g_{-}$) la somme de ses id\'eaux simples, $\g_{i}$, pour
lesquels la restriction de $B$ \`a $\g_{i}$ est \'egal \`a
$K^{\g_{i}}_{\lambda_{i}}$, avec $\lambda _{i}\in \C^{+}$ (resp.
$\C^{-}$).\begin{lem}
Aucune sous-alg\`ebre de Lie semi-simple complexe de $\g_{+}$
(resp. $\g_{-}$) n'est isotrope pour $B$.)
\end{lem}
On fait la d\'emontration pour $\g_{+}$, celle pour $\g_{-}$ \'etant
identique. Soit $\s$  une sous alg\`ebre de Lie semi-simple
complexe de $\g_{+}$. On note $\k_{\s}$ une forme r\'eelle
compacte de $\s$.  Alors $\k_{\s}$ est contenu dans une forme
r\'eele compacte de $\g_{+}$. En effet, le sous-groupe analytique
de $G_{+}$ d'alg\`ebre de Lie $\k_{\s}$ est compact, comme groupe
de Lie connexe d'alg\`ebre de Lie semi-simple compacte. Il est
donc contenu dans un sous-groupe compact maximal $K$ de
$G_{+}$, et son alg\`ebre de Lie, $\k$, convient. On note $\g_{i}$,
$i=1, \dots , p$, les id\'eaux simples de $\g_{+}$. Alors
$\k=\oplus _{i=1, \dots , p}\k_{i}$, o\`u $\k_{i}=\k \cap
\g_{i}$. Le Lemme r\'esultera de la preuve de :
\beq
B(X,X)\not= 0, \>\> X\in \k \setminus \{0\}
\eeq
Pour cela on remarque que :
\beq \sum_{i=1,\dots , p}\lambda_{i}x_{i}\not= 0, \>\> si\>\>,
pour \>\> i=1, \dots, p, \>\> \lambda_{i}\in \C^{+}\>\> et \>\>
x_{i}\geq 0, \>\> non \>\> tous \>\> nuls
\eeq
Maintenant, si $X\in \k \setminus \{0\}$ et $X=\sum_{i=1,\dots,
p}X_{i}$, o\`u $X_{i}\in \k_{i}$, $K_{\g_{i}}(X_{i}, X_{i})$ est
n\'egatif o\`u nul,  et non nul pour au moins un indice $i$. Alors
(3.1), et donc le Lemme,  r\'esulte de (3.2) et de la d\'efinition de
$\g_{+}$.
\begin{lem}
Soit $\i$ une sous-alg\`ebre Lagrangienne sous $\p$. On note
$\h=\i\cap \m$, $\m_{+}=\m\cap \g_{+}$, $\m_{-}=\m\cap \g_{-}$.
\ste
On a $\m=\m_{+}\oplus \m_{-}$. Par ailleurs la f-involution,
$\sigma$,  dont $\h$ est l'espace des points fixes, induit un
morphisme bijectif, $\tau$, entre $\m_{+}$ et $ \m_{-}$ tel que
:  $$B(\tau (X),\tau (X))=-B(X, X), \>\> X\in \m_{+}$$
\end{lem}
{\em D\'emonstration :} L'involution $\sigma$ permute, sans
point fixe, d'apr\`es le Lemme 5, les id\'eaux simples de $\m$.
Ceux ci sont contenus dans des id\'eaux simples de $\g$, donc
contenus soit dans $\m_{+}$, soit dans $\m_{-}$. Donc
$\m=\m_{+}\oplus \m_{-}$. Si $\sigma$ envoyait un id\'eal simple de
$\m_{+}$ dans un autre id\'eal de $\m_{+}$, l'alg\`ebre de Lie
$\m_{+}$,  donc aussi
$\g_{+}$, contiendrait une sous-alg\`ebre semi-simple complexe
isotrope. C'est impossible, d'apr\`es le Lemme pr\'ec\'edent. Donc
$\sigma$ envoie tout id\'eal simple de $\m_{+}$ dans $\m_{-}$. Le
Lemme en r\'esulte imm\'ediatement. \qed
{\bf Notations} \ste
{\em On notera $\j_{+}=\j_{0}\cap \g_{+}$,
$\a_{+}=\a\cap
\j_{+}$. on d\'efinit de m\^eme  $\j_{-}$ et $\a_{-}$. La restriction de
$\j_{0}$ \`a $\j_{+}$ identifie les racines de $\j_{0}$ dans
$\g_{+}$ \`a celles de $\j_{+}$. On note ${\tilde R}_{+}$
l'ensemble de celles-ci , $\Sigma_{+}$, l'ensemble
des racines simples de l'ensemble de  racines positives,
${\tilde R}_{+}^{+}$, de ${\tilde R}_{+}$, form\'e des \'el\'ements de
${\tilde R}_{+}$ qui sont des poids de $\j_{+}$ dans $\b_{0}\cap
\g_{+}$. On d\'efinit de m\^eme ${\tilde R}_{-}$, relativement \`a
$\g_{-}$
et
$\j_{-}$.\ste
 On d\'efinit aussi
$\Sigma_{-}$, l'ensemble
des racines simples de l'ensemble de  racines positives,
${\tilde R}_{-}^{+}$, de ${\tilde R}_{-}$, form\'e des \'el\'ements de
${\tilde R}_{-}$ qui sont des poids de $\j_{-}$ dans $\b'_{0}\cap
\g_{-}$.
\ste On d\'efinit, pour
$\i$ comme dans le Lemme pr\'ec\'edent,
$R_{+}$, l'ensemble des racines de $\j_{+}$ dans $\m_{+}$,
$\Gamma_{+}=\Sigma_{+}\cap R_{+}$. \ste Puis on
d\'efinit comme ci-dessus $R_{-}$ et $\Gamma_{-}$.\ste La
restriction de $\tau$ \`a $\a^{+}:=\m_{+}\cap\j_{0}$ d\'efinit
une bijection entre $\a^{+}$ et  $\a^{-}:=\m_{-}\cap\j_{0}$,
dont l'inverse de la transpos\'ee induit une bijection, not\'ee
$A$, ente $R_{+}$ et $R_{-}$. On notera, pour $\alpha \in
R_{+}$, $A\alpha$, au lieu de $A(\alpha)$. }
\ste \ste
Soit $\alpha\in {\tilde R}$. On note $H_{\alpha}$ l'\'el\'ement de
$\j_{0}$, tel que $\alpha (H_{ \alpha})=2$ et qui est orthogonal au
noyau de $\alpha$ pour la forme de Killing de $\g$. Soit
$\alpha\in R_{+}$ et
$\beta=A\alpha$. Alors, on voit facilement que
$H_{\alpha}$ (resp. $H_{\beta}$) est  \'el\'ement de $\a^{+}$ (resp.
 $\a^{-}$).   Comme $\tau$ transporte la forme
de Killing de $\m_{+}$ sur celle de $\m_{-}$, et que celles-ci sont
proportionnelles \`a la restriction de la forme de Killing de $\g$, on
a :
$$\tau (H_{\alpha})=H_{\beta}$$
Le Lemme 15 implique donc :
\beq
B(H_{A\alpha}, H_{A\beta})=-B(H_{\alpha},H_{ \beta}), \>\> \alpha,
\beta
\in
\Gamma_{+}
\eeq
Soit $\i'$ une autre sous-alg\`ebre Lagrangienne de $\g$, pour
laquelle on introduit des objets similaires, not\'es avec des
$'$.\ste
On notera $C$, ou parfois $''A^{-1}A'\>''$,  l'application
d\'efinie sur la partie, \'eventu-\ste ellement vide :
\beq dom\> C:=\{\alpha\in R'_{+}\vert A'\alpha \in R_{-}\}\eeq
par :
\beq C\alpha=A^{-1}A'\alpha,\>\> \alpha \in dom \> C \eeq
l'image de $C$ \'etant \'egale \`a :
\beq Im \>C=\{\alpha \in R_{+}\vert A\alpha \in R'_{-}\} \eeq
\begin{lem} Soit $(B,\i,\i')$,  un triple fortement standard.
Avec les notations pr\'ec\'edentes,  pour tout $\alpha \in
dom\> C$, il existe $n\in \N^{*}$
tel que :$$\alpha, \dots, C^{n-1}\alpha \in dom\> C\>\> et
\>\> C^{n}\alpha\notin dom\> C $$
\end{lem}
{\em D\'emonstration :} Soit $\alpha\in dom\> C$. Supposons que
pour tout $n\in \N^{*}$, $C^{n}\alpha$ soit d\'efini et \'el\'ement
de $dom \> C$. Comme $R'_{+}$ est un ensemble fini, il existe
$n_{1}$, $n'_{1}\in \N^{*}$, distincts, tels que $C^{n_{1}}\alpha =
C^{n'-{1}}\alpha$. D'o\`u l'on d\'eduit l'existence de $n\in
\N^{*}$ tel que $C^{n}\alpha=\alpha$. On note :
$$\alpha_{1}=\alpha, \alpha_{2}=C\alpha, \dots ,
\alpha_{n}=C^{n-1}\alpha$$
Montrons que :
\beq \{\alpha_{1}, \dots , \alpha_{n}\}\subset R_{+}\cap R'_{+},
\>\> A(\{\alpha_{1}, \dots , \alpha_{n}\})=
 A'(\{\alpha_{1}, \dots , \alpha_{n}\})
\eeq
 En effet, pour tout $i$, $\alpha_{i}$ est \'el\'ement de $dom\>
C\cap Im \> C$, ce qui implique la premi\`ere inclusion. Par
ailleurs, pour $i=1, \dots, n-1$, comme $''A^{-1}A'\>
''\alpha_{i}=\alpha_{i+1}$, on a :
$$A'\alpha_{i}=A\alpha _{i+1}, \>\> i=1, \dots , n-1$$
Maintenant, comme $''A^{-1}A'\>
''\alpha_{n}=\alpha_{1}$, on a aussi
$$A'_{\alpha_{n}}=A_{\alpha_{1}}$$ Ceci ach\`eve de prouver (3.7). \ste
On note $R_{+}^{0}$ (resp. $R_{-}^{0}$), l'intersection de
$R_{+}$ (resp. $R_{-}$) avec le sous-espace vectoriel r\'eel,
$V_{0}$, de $\j_{0}^{*}$, engendr\'e par les $\alpha_{i}$ (resp.
$A\alpha_{i}$), $i=1,\dots, n$. L'identification de $\j_{0}$ \`a
$\j_{0}^{*}$, \`a l'aide de la forme de Killing de $\g$, fait
appara\^{i}tre $(V_{0})_{\C}$ comme (le dual d') un sous-espace
vectoriel complexe $\a^{+}_{0}$  de $\a^{+}$. On d\'efinit de
m\^eme $\a^{-}_{0}$. On note $\m_{+}^{0}$ (resp. $\m_{-}^{0}$) la
sous-alg\`ebre de Lie de $\g$ engendr\'ee par les espaces
radiciels $ \m_{+}^{\alpha}$, $\alpha \in R^{0}_{+}$ (resp.
  $ \m_{-}^{\alpha}$, $\alpha \in R^{0}_{-}$). On voit
imm\'ediatement que $\m_{+}^{0}$ est semi-simple et que :
$$ \m_{+}^{0}= (\oplus_{\alpha \in R^{+}_{0}}  \m_{+}^{\alpha})
\oplus \a^{+}_{0}$$
 car $R_{+}^{0}$ est un syst\`eme de racines dans $\a^{+}_{0}$
([Bou], Ch. 7, Par. 1, Proposition 4). Il en va  de m\^eme pour
$\m_{-}^{0}$. Alors
$\tau
$ et
$\tau'$ induisent un isomorphisme entre $\m_{+}^{0}$ et
$\m_{-}^{0}$. Donc $\tau'^{-1}\tau $ induit un automorphisme de
$\m_{0}^{+}$, qui a donc un point fixe non nul $X$. Alors, on
a :
$$\tau(X)=\tau'(X)\>\> et \>\> X+ \tau(X)\in
\h\cap\h'\subset\i\cap \i'$$
Une contradiction qui ach\`eve de prouver la propri\'et\'e voulue
pour $C$.
\qed
\begin{prop}
 Soit $(B,\i,\i')$,  un triple fortement standard. \ste
Alors, pour tout $\alpha\in R_{+}$ (resp. $R'_{+}$), $\alpha$
et $A\alpha $ (resp. $A'\alpha $), sont de m\^eme signes
(relativement aux  ensembles de racines positives ${\tilde
R}^{+}_{+}$,  ${\tilde R}^{+}_{-}$, d\'efinis plus haut)
\end{prop}
{\em D\'ebut de la d\'emonstration :} On raisonne par r\'ecurrence sur
la dimension de $\g^{der}$. Si celle-ci est nulle, le r\'esultat
est clair. On suppose maintenant que celle-ci n'est pas nulle,
et que la Proposition est vraie pour les alg\`ebres r\'eductives
dont l'id\'eal d\'eriv\'e est de dimension strictement inf\'erieure \`a
celle de $\g^{der}$. Nous allons commencer par \'etablir
plusieurs Lemmes.

\begin{lem}
Avec les notations pr\'ec\'edentes, on a :\ste
(i) Si $\alpha\in R_{+}$ et $\alpha
\notin Im\> C$, $\alpha $ et $A\alpha$ sont de m\^emes signes.
\ste  (ii)  Si $\alpha\in R'_{+}$, et $\alpha \notin Dom\> C$,
$\alpha $ et $A'\alpha$ sont de m\^eme signes.
\end{lem}
{\em D\'emonstration : } Montrons (i). Raisonnons par l'absurde, et
supposons que $\alpha$ et $\beta: =A\alpha$ soient de signes
oppos\'es.  L'hypoth\`ese sur $\alpha$ \'equivaut \`a : \beq \alpha\in
R_{+}, \>\> A\alpha
\notin R'_{-}
\eeq Quitte \`a changer $\alpha$ en $-\alpha$, on peut supposer
$\alpha$ positive. Soit
$X$ un \'el\'ement non nul de $\m^{-\alpha}_{+}\subset \b'_{0}$. Alors
$\tau(X)\in \m_{-}^{-\beta}\subset \b'_{0}$, d'apr\`es notre
hypoth\`ese sur $\alpha$, et la d\'efinition des ensembles de
racines positives,  ${\tilde
R}^{+}_{+}$,  ${\tilde R}^{+}_{-}$. Enfin $Y:=X+ \tau(X)$ est
un \'el\'ement non nul  de $\h$.\ste Comme $\m_{-}^{-\beta}\subset
\b'_{0}$ et que $-\beta=-A\alpha\notin  R'_{-}$, on a
$\m_{-}^{-\beta}\subset \n'$. \ste Supposons
d'abord $\alpha\notin R'_{+}$.
Comme $\alpha$ est positive, $\m^{-\alpha}_{+}$ est contenu dans $
\b'_{0}$ et
$\alpha\notin R'_{+}$,
 implique, comme ci-dessus, que  $\m^{-\alpha}_{+}$ est contenu dans
$\n'$. Alors $X+\tau(X)$ est un \'el\'ement non nul de de
$\h\cap\n'\subset \i\cap \i'$. Une contradiction qui montre
qu'on doit avoir $\alpha\in R'_{+}$. \ste Alors $X\in \m'$,
$\tau(X)\in \n'$,  $Y\in \h\cap \p'$ et $p^{\n'}(Y)=X$. Donc
$\m^{-\alpha}_{+}=\g^{-\alpha}$ est contenu dans $\i_{1}$, o\`u
$(B_{1}, \i_{1}, \i'_{1})$ est le triple ant\'ec\'edent de $(B,
\i,\i')$. Appliquant  le Lemme 13 \`a $\i_{1}$, on voit
qu'alors $\g^{-\alpha}$ est contenu dans $\n_{1}$. Mais comme
le triple $(B,
\i,\i')$ est fortement standard les poids de $\j_{0}$ dans
$\n_{1}$ doivent \^etre des poids de $\j_{0}$ dans $\b_{0}$, ce
qui n'est pas le cas de $-\alpha$. Ceci ach\`eve de prouver (i).
Pour (ii), l'hypoth\`ese se traduit par une condition analogue \`a
(3.8), en \'echangeant le r\^ole de $\i$ et $\i'$. On d\'eduit donc
(ii) de (i). \qed
\begin{lem}
Si $\xi\in R_{+}\cap Im\> C$, il existe $\alpha \in dom\> C$,
$\alpha\notin Im \> C$, et $n\in \N^{*}$ tels que :
$$\alpha, C\alpha, \dots C^{n-1}\alpha \in dom \> C, \>\>et
\>\> C^{n}\notin dom \> C$$
et v\'erifiant :
$$\xi= C^{i}\alpha, \>\> pour \>\> un \>\> i\in \{1,\dots ,
n\}$$
\end{lem}
\dem Comme $C$ est une bijection de $dom \>C$ sur $Im\> C$, car
$A$ et $A'$ sont injectives, on note $C^{-1}$ la bijection
r\'eciproque. Echangeant le r\^ole de $\i$ et $\i'$ dans le Lemme
17 , on voit qu'il existe $n'\in \N^{*}$ tel que : $$\xi,
C^{-1}\xi, \dots, (C^{-1})^{n'-1}\xi\in Im \> C, \>\>
(C^{-1})^{n'}\xi \notin Im\> C$$
On pose $\alpha= (C^{-1})^{n'}\xi\in dom \> C$. On a aussi
$\alpha \notin Im\> C$. On choisit, gr\^ace au Lemme 17, $n\in
\N^{*}$ tel que :
$$\alpha, C\alpha, \dots C^{n-1}\alpha \in dom \> C, \>\>et
\>\> C^{n}\alpha\notin dom \> C$$ Alors $\alpha$ et $n$ ont
clairement les propri\'et\'es voulues. \qed
\begin{lem}
Soit $\alpha \in R_{+}\cap R_{-}$ tel que $A\alpha\in R_{-}\cap
R'_{-}$ (resp. $A'\alpha\in R_{-}\cap
R'_{-}$). Alors $\alpha$ et $A\alpha$ (resp. $A'\alpha$) sont
de m\^eme signe.
\end{lem}
\dem Prouvons l'assertion sur $A\alpha$ et soit $\alpha$ comme
dans l'\'enonc\'e. Soit $X\in \g^{\alpha}$. Les hypoth\`eses
impliquent que $\g^{\alpha}$ et $\g^{A\alpha}$ sont  contenus
dans $\m\cap\m'$. Donc, on a : $$\g^{\alpha}\subset (\m\cap\m')\cap
\sigma (\m\cap \m')$$
et d'apr\`es le Lemme  14, on a $\g^{\alpha}\subset \m_{1}$. De plus,
d'apr\`es ce m\^eme Lemme, l'involution $\sigma_{1}$ est la
restriction de
$\sigma$ \`a $\m_{1}$. On va appliquer l'hypoth\`ese de
r\'ecurrence du d\'ebut de la d\'emonstration de la Proposition 7.
Pour cela on remarque $(\l\cap \l')_{+}=(\l\cap\l')\cap
\g_{+}$ et de m\^eme pour $(\l\cap \l')_{-}$. Donc les racines de
$\j_{0}$ dans $\l\cap \l'$ qui sont positives dans
$\l\cap\l'$ sont positives dans $\g$. L'application de l'hypoth\`ese
de  r\'ecurrence montre alors que $\alpha$ et $A\alpha$ sont de m\^eme
signe. Ceci ach\`eve de prouver (i). Alors l'assertion sur
$A'\alpha$ r\'esulte de celle sur $A\alpha$, par \'echange du r\^ole de
$\i$ et $\i'$\qed
 \begin{lem} \ste
(i) Soit $\alpha\in R_{+}\cap
Im\> C$, $\alpha\notin R'_{+}$. Alors $\alpha$ et $A\alpha$ sont
de m\^eme signe.\ste (ii)  Soit $\alpha\in R'_{+}\cap dom\> C$,
$\alpha\notin R_{+}$. Alors $\alpha$ et $A'\alpha$ sont de m\^eme
signe.\end{lem}
\dem D\'emontrons (i). Soit $\alpha $ comme dans l'\'enonc\'e.
Quitte \`a changer $\alpha$ en $-\alpha$, on peut supposer que
$\alpha$ est n\'egative. Raisonnons par l'absurde et supposons
$A\alpha $ positive. Soit $X\in \g^{\alpha}$. Nos hypoth\`eses
montrent que $\g^{\alpha}$ est contenu dans $\n'$. Par
ailleurs, \'ecrivant $\alpha= ''A^{-1}A'\>''\beta$, o\`u $\beta \in
dom\> C$, on a $A\alpha= A'\beta \in R'_{-}$. Donc
$\g^{A\alpha}$ est contenu dans $\m'_{-}$. Alors :
$$X+\sigma(X)\in \h\cap \p', \>\>
p^{\n'}(X+\sigma(X))=\sigma(X)$$
Il en r\'esulte que $\g^{A\alpha}$ est contenu dans $\i_{1}$,
donc dans $\n_{1}$, d'apr\`es le Lemme 13. Alors $\g^{A\alpha}$
doit \^etre contenu dans $\b_{0}$, puisque le triple $(B_{1},
\i_{1}, \i'_{1})$ est  standard. Mais, comme $A\alpha\in R_{-}$
et est positive, ce n'est pas le cas. Une contradiction qui
ach\`eve de prouver (i). (ii) se d\'eduit de (i) par l'\'echange de
$\i$ et $\i'$.\qed
{\em Fin de la d\'emonstration de la Proposition 7: }\ste
Soit $\alpha \in R_{+}$ et montrons que $A\alpha $ est de m\^eme
signe que $\alpha$. Distinguons 3 cas. \ste
1) Si $\alpha \notin Im \> C$, cela r\'esulte du Lemme 18.\ste
2) Si $\alpha \in Im\> C$ et $\alpha \notin R'_{+}$, cela
r\'esulte du Lemme pr\'ec\'edent. \ste
3) Si $\alpha \in Im\> C$ et $\alpha \in R'_{+}$, on \'ecrit
$\alpha= ''A^{-1}A'\>''\beta$, o\`u $\beta \in
dom\> C$. Alors on a $A\alpha= A'\beta \in R_{-}\cap R'_{-}$.
On
conclut gr\^ace au Lemme 20. \ste  On vient donc de montrer
que pour tout $\alpha\in R_{+}$, $\alpha$ et $A\alpha$
sont de m\^eme signe.\ste On d\'emontre un \'enonc\'e similaire pour
$A'$ en \'echangeant le r\^ole de $\i$ et $\i'$. Ceci ach\`eve la
d\'emonstration de la Proposition.\qed \ste
On rappelle que $H_{\alpha}$, $\alpha \in {\tilde R}$ a \'et\'e d\'efini
avant (3.3). C'est la coracine correspondant \`a $\alpha$.  On
rappelle qu'un  syst\`eme de g\'en\'erateurs de Weyl de
$\g^{der}$ est une famille ${\cal W}=(H_{\alpha},X_{\alpha},
Y_{\alpha})_{\alpha\in
\Sigma}$, o\`u
$\Sigma=\Sigma_{+}\cup \Sigma_{-}$,  telle que,  pour tout $\alpha,
\beta \in \Sigma$, on ait :
\beq [X_{\alpha}, Y_{\beta}]=\delta_{\alpha\beta} H_{\beta} \eeq
\beq [H_{\alpha}, X_{\beta}] = N_{\alpha\beta }X_{\beta} \eeq
\beq [H_{\alpha}, Y_{\beta}]=-N_{\alpha\beta }Y_{\beta}\eeq
o\`u: \beq N_{\alpha\beta}=\beta(H_{\alpha})=2K_{\g}(H_{\alpha},
H_{\beta})/K_{\g}(H_{\alpha},
H_{\alpha})\eeq
On a alors :
$$N_{\alpha\beta}\in -\N, \>\> si\>\>\alpha \not= \beta$$
et : \beq ad X_{\alpha}^{1-N_{\alpha\beta}}X_{\beta}=ad
Y_{\alpha}^{1-N_{\alpha\beta}}Y_{\beta} =0,\>\> si\>\> \alpha \not=
\beta\eeq
 Un tel syst\`eme existe et on obtient les autres par
conjugaison par les \'el\'ements de $J_{0}$.
\begin{defi}
On appelle donn\'ee de Belavin-Drinfeld g\'en\'eralis\'ee, relativement
 \`a
$B$, la donn\'ee de
\ste $(A,A', \i_{\a},\i_{\a'})$, o\`u :\ste
1) $A$ est une bijection d'une partie $\Gamma_{+}$ de $\Sigma_{+}$
sur une partie $\Gamma_{-}$ de $\Sigma_{-}$, telle que :
\beq
B(H_{A\alpha}, H_{A\beta})=-B(H_{\alpha},H_{ \beta}), \>\> \alpha,
\beta
\in
\Gamma_{+}
\eeq
2) $A'$ est une bijection d'une partie $\Gamma'_{+}$ de $\Sigma_{+}$
sur une partie $\Gamma'_{-}$ de $\Sigma_{-}$, telle que :
 \beq
B(H_{A'\alpha}, H_{A'\beta})=-B(H_{\alpha},H_{ \beta}), \>\> \alpha,
\beta
\in
\Gamma'_{+}
\eeq
3) On d\'efinit  $C=''A^{-1}A' \>''$ comme dans (3.4), (3.5).
Alors $C$ satisfait la ''condition de sortie '' :\ste
Pour tout $\alpha \in
dom\> C$, il existe $n\in \N^{*}$
tel que $\alpha, \dots, C^{n-1}\alpha \in dom\> C\>\> et
\>\> C^{n}\alpha\notin dom\> C $.
\ste
4)  $\i_{\a}$ (resp. $\i_{\a'}$)  est un sous-espace vectoriel
complexe de
$\j_{0}$,  contenu et Lagrangien dans l' orthogonal, $\a$
(resp. $\a'$) pour la forme de Killing de
$\g$ (ou pour $B$), \`a l'espace engendr\'e par les
$H_{\alpha}$, $\alpha \in \Gamma:=\Gamma_{+}\cup\Gamma_{-}$ (resp.
$\Gamma':=\Gamma'_{+}\cup\Gamma'_{-}$).\ste
 5) Notons
$\f$ le sous-espace de
$\j_{0}$ engendr\'e par la famille $H_{\alpha}+H_{A\alpha}$, $\alpha\in
\Gamma_{+}$. On d\'efinit de m\^eme $\f'$. Alors : \beq (\f\oplus\i_{\a}
)\cap (\f'\oplus \i_{\a'})=\{0\}\eeq
On notera alors $R_{+}$ le sous-syst\`eme de racines de ${\tilde R}$
form\'e des \'el\'ements de ${\tilde R}$ qui sont combinaison lin\'eaire
d'\'el\'ements de $\Gamma_{+}$. On d\'efinit de m\^eme $R_{-}$, $R'_{+}$,
$R'_{-}$. On notera encore $A$ (resp. $A'$) le prolongement par
$\R$-lin\'earit\'e de $A$ (resp. $A'$), qui d\'efinit une bijection de
$R_{+}$ sur $R_{-}$ (resp. $R'_{+}$ sur $R'_{-}$).

\end{defi}
\begin{lem}
Si $A$ v\'erifie la condition 1) ci dessus, il existe un unique
isomorphisme, $\tau $, de la sous-alg\`ebre $\m_{+}$ de $\g$, engendr\'ee
par les $X_{\alpha},H_{\alpha}, Y_{\alpha}$, $\alpha \in
\Gamma_{+}$, sur la sous-alg\`ebre $\m_{-}$ de $\g$, engendr\'ee
par les $X_{\alpha},H_{\alpha}, Y_{\alpha}$, $\alpha \in
\Gamma_{-}$, tel que :
$$\tau (H_{\alpha})=H_{A\alpha}, \tau (X_{\alpha})=X_{A\alpha},
\tau (Y_{\alpha})=Y_{A\alpha}, \alpha \in
\Gamma^{+}$$
\end{lem}
En effet $\m_{+}$ et $\m_{-}$sont semi-simples  et les
familles donn\'ees sont des syst\`emes   de g\'en\'era- \ste teurs
de Weyl
de ces alg\`ebres. Alors, d'apr\`es [Bou], Chapitre VIII, Paragraphe 4.3,
Th\'eor\`eme 1, il suffit, pour montrer le Lemme,  de voir que les
relations, du type (3.9) \`a (3.13),  satisfaites par ces g\'en\'erateurs se
correspondent, c'est \`a dire qu'il faut montrer :
\beq \alpha(H_{\beta})=A\alpha(H_{A\beta}), \>\> \alpha,
\beta
\in
\Gamma_{+}
\eeq
Montrons d'abord que, la deuxi\`eme d\'efinition de $N_{\alpha\beta}$
(deuxi\`eme \'egalit\'e de (3.12)), on peut
remplacer $K_{\g}$  par n'importe qu'elle autre forme $\g$-invariante
non d\'eg\'en\'er\'ee. En effet si $\g^{\alpha}$ n'est pas dans le m\^eme
id\'eal simple que $\g^{\beta}$, $B(H_{\alpha},H_{
\beta})=K_{\g}(H_{\alpha},H_{ \beta})=0$. Sinon $B$ et $K_{\g}$ sont
proportionnelles sur l'id\'eal simple contenant $\g^{\alpha}$ et
$\g^{\beta}$. D'o\`u notre assertion. Alors (3.17) r\'esulte
imm\'ediatement
de la condition 1).\qed
\begin{prop}
Soit ${\cal BD} =(A,A', \i_{\a}, \i_{\a'})$ une  donn\'ee de
Belavin-Drinfeld g\'en\'eralis\'ee, relative \`a $B$.  On note $\p$ la
sous-alg\`ebre parabolique de $\g$, contenant
$\b_{0}$ et $\m$. Sa d\'ecomposition de Langlands $\p=\l\oplus
\n$, o\`u $\l$ contient $\j_{0}$, v\'erifie $\l=\m\oplus \a$. On
utilise les notations du Lemme pr\'ec\'edent. On note $\i:=\h\oplus
\i_{\a}\oplus
\n$, o\`u
$\h:=\{X+\tau (X)\vert X\in
\m_{+}\}$. On d\'efinit de m\^eme $\i'$. \ste
Alors $(B,\i,\i')$ est un triple de Manin fortement standard.
 On dira que ce triple de Manin est associ\'e \`a la donn\'ee de
Belavin-Drinfeld g\'en\'eralis\'ee, ${\cal BD}$, et au syst\`eme de
g\'en\'erateurs de Weyl, ${\cal W}$. On le notera
${\cal T}_{{\cal BD, W}}$.\ste
On note ${\cal W}_{1}=(H_{\alpha},X_{\alpha}, Y_{\alpha})_{\alpha\in
\Gamma\cap \Gamma'}$. C'est un syst\`eme de g\'en\'erateurs de Weyl de
$(\l\cap\l')^{der}$.
Notons $\Gamma_{1+}=\Gamma_{+}\cap \Gamma'_{+}\cap
A^{-1}(\Gamma_{-}\cap \Gamma'_{-})$. On note $A_{1}$ la restriction
de
$A$ \`a  $\Gamma_{1+}$ et  $\Gamma_{1-}$ son image. On d\'efinit de
m\^eme $A'_{1}$. On note $\a_{1}$ l'intersection des noyaux des
\'el\'ements de $\Gamma_{1+}\cup \Gamma_{1-}$ et on note
$\i_{\a_{1}}=\i_{\a}\oplus \t_{1}$, o\`u $\t_{1}$ est l'intersection
de $\f$ avec $\a_{1}$. On d\'efinit de m\^eme
$\i_{\a'_{1}}$. Alors
$(A_{1}, A'_{1},
\i_{\a_{1}},
\i_{\a'_{1}})$ est une donn\'ee de Belavin-Drinfeld g\'en\'eralis\'ee,
${\cal BD}_{1}$,  pour
$\l\cap\l'$ et l'ant\'ec\'edent du triple ${\cal
T}_{{\cal BD, W}}$ est \'egal \`a ${\cal
T}_{{\cal BD}_{1},{\cal  W}_{1}}$.
\end{prop}
\dem
On proc\`ede par r\'ecurrence sur la dimension de $\g^{der}$. Si
celle-ci est nulle,  le r\'esultat est clair. On suppose le r\'esultat
est vrai pour les alg\`ebres r\'eductives dont l'id\'eal d\'eriv\'e est
de
dimension strictement inf\'erieure \`a celle de $\g^{der}$.\ste
 Montrons  que
$\h$ est isotrope pour
$B$. Etudions la forme bilin\'eaire, $B'$, sur $\m_{+}$, d\'efinie par :
$$B'(X,Y)=-B(\tau(X),\tau(Y)), \>\> X,Y\in \m_{+}$$ C'est
clairement une forme bilin\'eaire invariante sur $\m_{+}$, qui
coincide sur $\a^{+}=\j_{0}\cap \m_{+}$ avec la restriction de
$B$ \`a $\m_{+}$, d'apr\`es (3.14) et le Lemme 22. Comme $\a^{+}$ est
une
\sac de
$\m_{+}$, le Lemme 1 permet de voir que $B'$ coincide sur $\m_{+}$
avec cette restriction. Comme $\g_{+}$ et $\g_{-}$ sont orthogonaux
pour
$B$, il en r\'esulte que $\h$ est isotrope. Par ailleurs, la
d\'efinition de  $\h$ montre que c'est l'espace des points
fixes d'une f-involution, $\sigma$, de $\m$, dont la
restriction \`a $\m_{+}$ est \'egale \`a $\tau$. En particulier
$\sigma$ permute $\m_{+}$ et $\m_{-}$. L'application du Th\'eor\`eme
1 montre que $\i$ est une sous-alg\`ebre de Lie de $\g$,
Lagrangienne pour $B$.\ste On montre de m\^eme que $\i'$ est
Lagrangienne pour $B$.\ste Pour montrer que $(B,\i,\i')$ est un
triple de Manin de $\g$,   on se propose d'appli-\ste quer le
Th\'eor\`eme
3.\ste  Montrons que $\h\cap \n'=\{0\}$. Soit $X\in \h\cap \n'$.
Alors, comme $\n'$ est orthogonal \`a $\j_{0}$, pour la forme de
Killing de $\g$, il en est ainsi de $X$.\
 Notons $R_{+}$
l'ensemble des poids non nuls  de $\j_{0}$ dans
$\m_{+}$. Ceux-ci sont en particulier nuls sur $\j_{-}$.  Alors
on a : \beq X=\sum _{\alpha\in R_{+}}(X(\alpha )+\tau (X(\alpha ))),
\>\>o\grave{u}\>\> X(\alpha ) \in \m_{+}^{\alpha}\eeq
On note que $\tau(X(\alpha ))$ est de poids $\beta$ sous
$\j_{0}$, o\`u $\beta$ est un poids de $\j_{0}$ dans $\m_{-}$,
donc nul sur $\j_{+}$.
La d\'ecomposition (3.18) apparait alors comme une d\'ecomposition
de $X$ en vecteurs poids sous $\j_{0}$, pour des poids deux
\`a deux distincts. Comme $X\in \n'\subset  \b'_{0}$, et que
$\m_{+}\in \g_{+}$, il faut que,  pour tout $\alpha \in R_{+}$,  avec
$\alpha$ positive, on ait $X(\alpha )=0$. D'autre part, si
$\alpha \in R_{+}$,  avec
$\alpha$ n\'egative, $\tau(X(\alpha))$ est un \'el\'ement de
$\m_{-}^{\beta}$, avec $\beta $ n\'egative. En effet $\beta$ est
l'image de $\alpha$ par le prolongement $\R$-lin\'eaire de $A$ \`a
l'espace vectoriel r\'eel engendr\'e par $\Gamma_{+}$. Alors
$\tau(X(\alpha))$ appartient \`a $\b_{0}$. Comme $X\in \n'$, on
doit avoir $\tau(X(\alpha))=0$, donc aussi $X(\alpha)=0$.
Finalement $X=0$, comme d\'esir\'e. Donc $\h\cap \n'=\{0\}$. On
montre de m\^eme $\h'\cap \n=\{0\}$.\ste  On note $${\tilde
\h}=\h\oplus \i_{\a},\>\> \f=\j_{0}\cap \h,\>\> {\tilde \f}=\f\oplus
\i_{\a}$$ Si
$\alpha\in R_{+}$, on note : $$\h_{\alpha}=\{X+\tau(X)\vert X\in
\m_{+}^{\alpha}\}$$  Alors, d'apr\`es (2.41) o\`u l'on prend
$R_{*}=R_{+}$,  on a :
\beq {\tilde \h}={\tilde \f}\oplus (\oplus_{\alpha_{\in
R_{+}}}\h_{\alpha})\eeq
Cette  d\'ecomposition
apparait  comme une d\'ecomposition de ${\tilde \h}$ en
repr\'esenta-\ste tions de $\f$ sans sous-quotients simples
\'equivalents, les $\h_{\alpha}$ \'etant de plus irr\'educ-\ste tibles.
Comme $\p'$contient ${\tilde \f}$, on en  d\'eduit :
\beq {\tilde \h}\cap \p'={\tilde \f}\oplus (\oplus_{\alpha\in
R_{+}, \h_{\alpha}\subset \p'}\h_{\alpha})\eeq
On note encore  $A$ le prolongement $\R$-lin\'eaire de $A$ au
sous-espace vectoriel r\'eel de $\j_{0}^{*}$ engendr\'e par
$\Gamma_{+}$. On a alors, pour tout $\alpha \in R_{+}$, $\sigma
(\g^{\alpha})=\g^{A\alpha}$ . Comme $\p'$ est somme
de sous-espaces poids sous $\j_{0}$, on a $\h_{\alpha}\subset \p'$ si
et seulement si $\alpha $ et $A\alpha$ sont des poids de
$\j_{0}$ dans $\p'$ c'est \`a dire :
\beq \alpha \in R_{+}\cap (R'_{+}\cup {\tilde R}^{-}) \>\> et
\>\>  A\alpha \in R_{-}\cap (R'_{-}\cup {\tilde R}^{+})\eeq
Calculons $p^{\n'}(\h_{\alpha}) $ pour $\alpha$ satisfaisant
(3.21). On distingue quatre cas.
\ste 1) $$ \alpha \in R_{+}\cap R'_{+},  \>\>A\alpha \in R_{-}\cap
R'_{-}$$ \ste Dans ce cas , si $X\in \g^{\alpha}$, on a $X\in
\m'_{+}$, $\tau(X)\in \m'_{-}$. Alors $p^{\n'}(X+\tau (X))=
X+\tau (X)$. Donc :$$ p^{\n'}(\h_{\alpha})=\h_{\alpha}$$
2) $$ \alpha \in R_{+}\cap R'_{+},  \>\>A\alpha\in {\tilde
R}^{+}\setminus R'_{-}$$  On trouve comme ci-dessus que  si
$X\in \g^{\alpha}$, on a $X\in
\m'_{+}$, $\tau(X)\in \n'$. Alors $p^{\n'}(X+\tau (X))=X$. D'o\`u l'on
d\'eduit que : $$ p^{\n'}(\h_{\alpha})=\m_{+}^{\alpha}$$
 3) $$ \alpha \in {\tilde R}^{-}\setminus R'_{+},\>\> A\alpha \in
R_{-}\cap R'_{-}$$ Si $X\in \g^{\alpha}$, on a $X\in
\n'$, $\tau(X)\in \m'_{-}$. Alors $p^{\n'}(X+\tau (X))=
\tau (X)$. Donc : $$ p^{\n'}(\h_{\alpha})=\m_{-}^{A\alpha}$$
4) $$ \alpha \in {\tilde R}^{-}\setminus R'_{+},\>\>A\alpha\in
{\tilde R}^{-}\setminus R'_{-}$$ Alors on aurait
$\h_{\alpha}\subset
\n'$ et   cette possibilit\'e est exclue, d'apr\`es ce qu'on a vu
plus haut.
\ste \ste Notons :
\beq \n_{1}:= (\oplus _{\alpha \in R^{+}_{+}\cap R'_{+},
A\alpha\notin R'_{-}}\m_{+}^{\alpha}) \oplus (\oplus _{\alpha
\in R^{-}_{+}\setminus R'_{+}, A\alpha \in
R'_{-}}\m_{-}^{A\alpha})
\eeq
et \beq \m_{1}= ((\m\cap\m' \cap \sigma (\m\cap \m'))^{der}
, \>\> \l_{1}=\m_{1}+\j_{0}\eeq
L'analyse des poids sous $\j_{0}$ montre que
$$\l_{1}= \j_{0}\oplus (\oplus_{ \alpha_{\in R_{1}}}\g^{\alpha})$$
o\`u :
$$R_{1}=(R_{+}\cap R'_{+}\cap A^{-1}(R_{-}\cap R'_{-}))\cup
(R_{-}\cap R'_{-}\cap A(R_{+}\cap R'_{+}))$$
En particulier $\l_{1}$ et $\n_{1}$ ont une intersection r\'eduite
\`a z\'ero. \ste
On note $\sigma_{1}$ la restriction de $\sigma $ \`a
$\m_{1}$,  qui est clairement une f-involution. En effet $\m_{1}$
est la somme directe de son intersection avec $\g_{+}$ et $\g_{-}$,
et
$\sigma_{1}$ permute ces deux id\'eaux. On note
$\h_{1}$ l'ensemble des points fixes de
$\sigma_{1}$. On d\'eduit de la d\'efinition de $\m_{1}$ ( cf. (3.23)),
et
de $\sigma_{1}$ que :
\beq \h_{1}=(\oplus _{\alpha\in R_{+}\cap R'_{+},A\alpha \in
R_{-}\cap R'_{-}}\h_{\alpha})\oplus (\f\cap \m_{1})
\eeq Notons :
$$\i_{1}=p^{\n'}({\tilde
\h}\cap
\p')$$  Gr\^ace \`a (3.20) et ce qui pr\'ec\'ede on voit que $\i_{1}$
est la
somme de son intersection, ${\tilde \f}$, avec $\j_{0}$, avec son
intersection avec l'orthogonal de $\j_{0}$ pour la forme de Killing
de $\g$. Tenant compte du fait que $A$ pr\'eserve le signe des racines,
on voit, gr\^ace \`a la discussion ci-dessus et \`a (3.20), que
l'intersection de
$\i_{1}$ avec l'orthogonal de $\j_{0}$ pour la forme de Killing
de $\g$, coincide avec celle de $\h_{1}\oplus \n_{1}$. On en
d\'eduit facilement l'\'egalit\'e :
\beq \i_{1}=(\h_{1}+{\tilde \f})\oplus \n_{1}
\eeq
 Montrons que
$\p_{1}:=\l_{1}\oplus
\n_{1}$ est une sous-alg\`ebre parabolique de $\l\cap \l'$. \ste
D'apr\`es la d\'efinition de $\p_{1}$, et de celle de $\m_{1}$,
$\n_{1}$ (cf (3.22), (3.23)), $\p_{1}$ est la somme de son
intersection
$\p_{1+}$ avec
$\g_{+}$ avec celle avec $\g_{-}$, $\p_{1-}$. Il suffit donc
d'\'etudier s\'epar\'ement ces intersections. On ne traite que $\p_{1+}$,
$\p_{1-}$ se traitant de la m\^eme mani\`ere.
Soit :
$$E=\{\alpha \in R^{+}_{+}\cap R'_{+}\vert
A\alpha\notin R'_{-}\}, \>\>F=R_{+}\cap R'_{+}\cap A^{-1}(R_{-}\cap
R'_{-})$$
Alors $E$, $F$, $E\cup F$  sont des parties
closes  du syst\`eme de racines de $\j_{0}$ dans $\l\cap \l'\cap\g_{+}$,
$R_{+}\cap R'_{+}$, $E\cup F$ en \'etant une partie parabolique,
contenant $R_{+}^{+}\cap R'^{+}_{+}$ et dont $E$ est un id\'eal (cf
[War], 1.1.2.9, 1.1.2.13, pour la terminologie). Cela r\'esulte du
fait que $R^{+}_{+}$, $R'_{+}$ et $R'_{-}$ sont
des parties closes, d'apr\`es leur d\'efinition (cf. fin de la
D\'efinition 5). Le seul point non imm\'ediat est le fait que si
$\alpha\in E$, $\beta
\in F $ et
$\alpha +\beta \in {\tilde R}$, alors $\alpha+\beta$ appartient \`a
$E$. D'apr\`es la d\'efinition de $R'_{-}$, on voit que nos hypoth\`eses
impliquent que $A(\alpha+\beta)\notin R'_{-}$. Il reste \`a voir que
$\alpha+\beta$ est positive. L'\'etude de ses  composantes dans la
base $\Sigma$, montre que l'une de celles-ci, correspondant \`a une
racine $\gamma\in \Gamma_{+}\cap A^{-1}\Gamma'_{-}$, est strictement
positive, car pour l'une au moins de ces racines , la composante de
$\alpha $ est strictement positive tandis que celle de $\beta$ est
nulle. Cela implique que $\alpha +\beta$ est positive et
finalement dans $E$.\ste
Joint \`a la d\'efinition de $\p_{1+}$, et celles de $\m_{1}$ et
$\n_{1}$ , cela implique que $\p_{1+}$ est une
sous-alg\`ebre parabolique de $\l\cap\l'\cap \g_{+}$. \ste On conclut
que
$\p_{1}$ est   une
sous-alg\`ebre parabolique de $\l\cap\l'$, qui contient
$\b_{0}\cap\l\cap\l'$. La d\'emonstration montre aussi que $\n_{1}$
en est le radical nilpotent.
 \ste Notons $
\f_{1}=\f\cap\m_{1}$ et  $\t_{1}$ l'intersection de  $\f$ avec le
centre $\a_{1}$ de $\l_{1}$.
Montrons que  :
\beq \f=\f_{1}\oplus \t_{1} \eeq
En effet tout \'el\'ement de $\f$ est de la forme $X+\tau (X)$, o\`u $X$
est un \'el\'ement de $\j_{0}\cap\m_{+}$. On note que $\m$ et $\m_{1}$
sont la somme directe de leurs intersection avec $\g_{+}$ et
$\g_{-}$. Il en va de m\^eme de $\j_{0}$, $\a$, $\a_{1}$. Alors
$X$ est la  somme d'un \'el\'ement de $\j_{0}\cap \m_{+}$ avec un
\'el\'ement de $\a_{1}\cap\g_{+}$. La d\'ecomposition ci-dessus en
r\'esulte aussitot.\ste
On pose :
\beq \i_{\a_{1}}= \t_{1}\oplus \i_{\a}\subset \a_{1}
\eeq
de sorte que (3.25) se r\'e\'ecrit
$$\i_{1}=\h_{1}\oplus \i_{\a_{1}}\oplus \n_{1}$$
On a :
$$ dim_{\C}\j_{0}=dim_{\C}(\j_{0}\cap\m_{+}) +
dim_{\C}(\j_{0}\cap\m_{-})+ dim_{\C}\a=dim_{\R}\f+ dim_{\C}\a$$
Comme $\i_{\a}$ est Lagrangienne dans $\a$, tenant compte de ce qui
pr\'ec\`ede, on a :
\beq  dim_{\C}\j_{0}=dim_{\R}\f_{1}+dim_{\R}\t_{1}+dim_{\R}
\i_{\a}\eeq d'o\`u l'on d\'eduit : \beq
dim_{\C}\j_{0}=dim_{\C}(\j_{0}\cap \m_{1})+dim_{\C}\a_{1}\eeq
Montrons :
 \beq dim_{\C}(\j_{0}\cap \m_{1})=dim_{\R}\f_{1} \eeq
En effet, $\f_{1}=\f\cap \m_{1}$, v\'erifie :
$$\f_{1}=\{X+\tau(X)\vert X\in \j_{0}\cap \m_{1+}\}$$
et  $$\j_{0}\cap\m_{1}=(\j_{0}\cap\m_{1+})\oplus
(\j_{0}\cap\m_{1-})$$ les deux facteurs \'etant \'echang\'es par $\tau$.
(3.30) en r\'esulte. On  d\'eduit de (3.28) \`a (3.30) que :
$$dim_{\R}\i_{\a_{1}}=dim_{\C}\a_{1}$$
Comme $\t_{1}\subset \h$ et $\i_{\a}$ sont isotropes pour $B$,
et orthogonales pour $B$, puisque le centre d'une alg\`ebre r\'eductive
est orthogonal \`a son id\'eal d\'eriv\'e pour toute forme invariante, ce
qui pr\'ec\`ede montre que :\beq \i_{\a_{1}}\>\> est
\>\>Lagrangienne\>\> dans \a_{1}\eeq  Cela implique que $\i_{1}$ est
Lagrangienne pour $B$.\ste  On introduit de la m\^eme mani\`ere
$\i'_{1}$. On prend $A_{1}$ \'egal \`a la restriction
de $A$ \`a $\Gamma_{1+}=\Gamma_{+}\cap \Gamma'_{+}\cap
A^{-1}(\Gamma_{-}\cap \Gamma'_{-})$. On d\'efinit de m\^eme $A'_{1}$.
Montrons que
$(A_{1},A'_{1},
\i_{a_{1}},\i_{\a'_{1}})$ est une donn\'ee de
Belavin-Drinfeld g\'en\'eralis\'ee. Les conditions 1), 2), 3)
r\'esultent imm\'ediatement des conditions satisfaites par
$(A,A',\i_{\a},\i_{\a'})$. La condition 4) r\'esulte de
(3.31). Enfin 5) r\'esulte de la condition 5) pour ${\cal
BD}$, joint \`a  (3.26) et (3.27). L'application de l'hypoth\`ese de
r\'ecurrence montre que $(B_{1}, \i_{1}, \i'_{1})$ un triple de Manin
associ\'e \`a
${\cal BD}_{1}$ et ${\cal W}_{1}$. Le Th\'eor\`eme 3 permet de conclure
que
$(B,\i,\i')$ est un triple de Manin, d'ant\'ec\'edent  $(B_{1}, \i_{1},
\i'_{1})$. Ceci ach\`eve la preuve de la Proposition.\qed
\begin{theo}
\ste (i) Tout triple de la forme ${\cal T}_{{\cal BD}, {\cal W}}$, o\`u
${\cal BD}$ est une donn\'ee de Belavin-Drinfeld g\'en\'era-\ste
lis\'ees, et
${\cal W}$ un ensemble de g\'en\'erateurs de Weyl de $\g$, est un
triple de Manin fortement standard.\ste (ii) R\'eciproquement
tout triple fortement standard est de cette forme.
\ste (iii)  Soit ${\cal T}_{{\cal BD}, {\cal W}}$ (resp
${\cal T}_{{\underline {\cal BD}}, {\cal {\underline W}}}$), o\`u
${\cal BD}$, $ {\underline {\cal BD}}$ sont des  donn\'ees de
Belavin-Drinfeld g\'en\'eralis\'ees, et  ${\cal W}$, $ {\cal
{\underline W}}$ des  ensembles de g\'en\'erateurs de Weyl de $\g^{der}$.
Ces triples de Manin sont conjugu\'es sous $G$, si et seulement si
${\cal BD}={\underline {\cal BD}}$. Alors ils sont conjugu\'es par
l'\'el\'ement de
$J_{0}$ qui conjugue
${\cal W}$ et  ${\cal {\underline W}}$.
\end{theo}
\dem
Le point (i) a \'et\'e vu \`a la Proposition pr\'ec\'edente. \ste
Prouvons (ii). Soit $(B,\i,\i')$ un triple fortement standard
sous $(\p,\p')$. On utilise les notations qui suivent le Lemme 16.
D'apr\`es la Proposition 7, l'application $A$ est une bijection de
$R_{+}$ sur $R_{-}$ qui pr\'eserve les signes des racines. Donc elle
induit une bijection de $\Gamma_{+}$ sur $\Gamma_{-}$, not\'ee
encore $A$. On a des propri\'et\'es analogues pour $A'$.  On note
$\i_{\a}=\i\cap\a$, $\i'_{\a'}=\i'\cap\a'$. On voit facilement que
$\f$ est l'espace vectoriel engendr\'e par les
$H_{\alpha}+H_{A\alpha}$, $\alpha\in \Gamma_{+}$ et que $\f\oplus
\i_{\a}=\i\cap\j_{0}$. On a des propri\'et\'es similaires pour $\i'$.
Alors (3.16) r\'esulte du fait que $\i$ et $\i'$ ont une intersection
r\'eduite \`a z\'ero. D'apr\`es le Th\'eor\`eme 1, le Lemme 17 et (3.3),
$(A, A',\i_{\a},
\i'_{\a'})$ est une donn\'ee de Belavin-Drinfeld g\'en\'eralis\'ee,
not\'ee
${\cal BD}$. Il reste \`a trouver un syst\`eme de g\'en\'erateurs de Weyl de
$\g$, v\'erifiant les relations du Lemme 22, avec $\tau$ comme dans
le Lemme 16. Si
$\alpha\in
\Gamma'_{+}$ est comme dans le Lemme 18, i.e. $\alpha  \in dom\> C$,
$\alpha\notin Im\> C$, on choisit
$X_{\alpha}\in \g^{\alpha}$, $Y_{\alpha}\in \g^{-\alpha}$,
tels que $[X_{\alpha}, Y_{\alpha}]=H_{\alpha}$. \ste
Puis notant
$\alpha_{i}=C^{i}\alpha$, $i=0, \dots , n$, on d\'efinit par
r\'ecurrence sur $i$, pour $Z=X$ ou $Z=Y$ :\beq
Z_{\alpha_{i+1}}:= \tau^{-1}
\tau '(Z_{\alpha_{i}})
\eeq l'expression du membre de droite
\'etant bien d\'efinie, pour $i=0,\dots,n-1$, car alors
$\alpha_{i}\in dom\> C$.\ste
Puis on pose, pour $Z=X$ ou $Z=Y$ :
\beq Z_{A\alpha_{i}}=\tau (Z_{\alpha_{i}}),\>\>
si\>\>\alpha_{i}\in\Gamma_{+}, \>\>Z_{A'\alpha_{i}}=\tau '
(Z_{\alpha_{i}}),\>\> si\>\>\alpha_{i}\in\Gamma'_{+}\eeq
Cette d\'efinition est coh\'erente, car si
$\beta =A{\alpha_{i}}=A'{\alpha_{i'}}$, on a $\alpha_{i'}\in
dom\> C$ et $C {\alpha_{i'}}={\alpha_{i}}$. Mais alors on a $i=
i'+1$. Il r\'esulte de (3.32) que les deux d\'efinitions de $Z_{\beta}$
de (3.33) coincident. \ste
Maintenant, soit un \'el\'ement $\beta$ de $\Gamma_{+}\cup \Gamma_{-}$
qui n'est pas de la forme $\alpha_{i}$, pour un $\alpha$ comme
ci-dessus. En particulier on a $\beta\notin dom\> C$, $
\beta\notin Im\> C$. On
choisit $X_{\beta }\in \g^{\beta}$, $Y_{\beta}\in \g^{-\beta}$,
tels que $[X_{\beta}, Y_{\beta}]=H_{\beta}$, puis on pose, pour
$Z=X$ ou $Z=Y$  : \beq Z_{A\beta}= \tau ( Z_{\beta}) \>\> si\>\>
\beta \in \Gamma_{+}, Z_{A'\beta}= \tau' ( Z_{\beta}) \>\> si\>\>
\beta \in \Gamma'_{+}
\eeq Montrons que cette d\'efinition est coh\'erente. Supposons que
$A\beta$ soit d\'efini et  \'egal \`a l' un des $A'\alpha_{i}$ ci
dessus. On aurait alors $\beta \in Im \>C$, ce qui est
impossible. Comme $ \beta\notin dom\> C$, on voit de m\^eme que
$A'\beta$, s'il est d\'efini ne peut-\^etre \'egal \`a l'un des
$A\alpha_{i}$.
Enfin si $A\beta=A'\beta'$, pour deux \'el\'ements
$\beta$, $\beta'$ comme ci-dessus,on aurait $\beta \in Im\>
C$, ce qui n'est pas. Les autres \'egalit\'es \`a envisager pour voir
la coh\'erence de (3.33) et (3.34) \'etant exclues, d'apr\`es la
bijectivit\'e de $A$ et $A'$, cette coh\'erence est donc prouv\'ee.\ste
Enfin si $\alpha\in \Sigma$, $\alpha\notin \Gamma_{+}\cup
\Gamma_{-}$, on choisit $X_{\alpha}\in \g^{\alpha}$, $Y_{\alpha}\in
\g^{-\alpha}$, tels que $[X_{\alpha}, Y_{\alpha}]=H_{\alpha}$.
Il est alors facile de voir que la famille
$(H_{\alpha},X_{\alpha} , Y_{\alpha})_{\alpha\in \Sigma}$ est
un syst\`eme de g\'en\'erateurs de Weyl, ${\cal W}$. De plus la
d\'efinition de $\tau$, $\tau'$ et (3.33) (3.34) montrent que le triple
$(B,\i, \i')$ est \'egal \`a ${\cal T}_{{\cal BD}, {\cal W}}$. Ceci
ach\`eve la preuve de (ii).
\ste Prouvons (iii). Supposons les deux
triples conjugu\'es par un \'el\'ement de $G$.  On note le premier triple
$(B,\i,\i')$ qu'on suppose sous $(\p,\p')$, on note ${\cal
BD}=(A,A',\i_{\a},
\i_{\a'})$ et on introduit des notations similaires pour le deuxi\`eme
triple, en soulignant. Montrons que  :
\beq \i\cap\j_{0}={\underline \i}\cap\j_{0},\>\>
\i'\cap\j_{0}={\underline \i'}\cap\j_{0}\eeq
On proc\`ede par r\'ecurrence sur la dimension de $\g^{der}$, le
r\'esultat \'etant clair si celle-ci est nulle. Par ailleurs, comme
deux sous-alg\`ebres paraboliques standard conjugu\'ees sous $G$ sont
\'egales, on a   \beq (\p,\p')=({\underline \p},{\underline \p'}) \eeq
et les deux triples sont conjugu\'es par un \'el\'ement de $P\cap P'$. De
la Proposition 4, on d\'eduit que les ant\'ec\'edents des deux triples,
qui sont fortement standard, sont conjugu\'es par un \'el\'ement de
$L\cap
L'$. L'application de l'hypoth\`ese de r\'ecurrence conduit au
 r\'esultat
voulu, car si
$(B_{1},\i_{1},\i'_{1})$ est l'ant\'ec\'edent de $(B,\i,\i')$, la
d\'efinition des triples fortement standard montre que
$\i\cap\j_{0}=\i_{1}\cap\j_{0}$, etc. L'\'egalit\'e (3.36) montre que
$\m={\underline \m}$, donc $\Gamma= {\underline \Gamma}$. De (3.35),
on d\'eduit l'\'egalit\'e de
$\f:=\m\cap\i\cap \j_{0}$ avec ${\underline \f}:=\m\cap
{\underline \i}\cap
\j_{0}$. Comme la d\'efinition de  ${\cal T}_{{\cal BD}, {\cal W}}$
montre que $\f$ est engendr\'e par $(H_{\alpha}+ H_{A\alpha})_{\alpha
\in \Gamma_{+}}$ et de m\^eme pour ${\underline \f}$,  l'\'egalit\'e de
$A$ et ${\underline A}$ en r\'esulte imm\'ediatement. Il en va de m\^eme
de l'\'egalit\'e de
$A'$ et ${\underline A'}$.
Comme $\i_{\a}= \i\cap\a$ et de m\^eme pour ${\underline
\i}_{{\underline \a}}$, on d\'eduit l'\'egalit\'e de ces espaces de
 (3.35)
, car
$\a={\underline \a}$ et $\i_{\a}=\i\cap\j_{0}\cap\a$ et de m\^eme pour
${\underline
\i}_{{\underline \a}}$. On proc\`ede de m\^eme pour $\i'_{\a'}$,
${\underline
\i'}_{{\underline \a'}}$.  Donc
${\cal BD}$ est \'egal \`a ${\underline {\cal BD}}$, comme d\'esir\'e.
Maintenant, il est clair qu'un \'el\'ement de $J_{0}$ qui conjugue
${\cal W}$ et
${\underline {\cal W}}$, conjugue les deux triples.\qed
\begin{rem}
Soit $\g_{1}$ une alg\`ebre de Lie simple complexe, $\j_{1}$, une \sac
de $\g_{1}$, $\g=\g_{1}\times \g_{1}$, $\j_{0}= \j_{1}\times \j_{1}$
et $B$ la forme $\C$-bilin\'eaire sur $\g$, $\g$-invariante, \'egale \`a
$K_{\g_{1}}$ sur le premier facteur et \`a $-K_{\g_{1}}$ sur le
deuxi\`eme
facteur. La classification de Belavin et Drinfeld de certaines
$R$-matrices (cf. [BD], Th\'eor\`eme 6.1) se r\'eduit, d'apr\`es l.c.
\'equations 6.1 \`a 6.5 , et [S], Propositions 1 et 2, \`a la
classification
des triples de Manin pour
$\g$,
$(B,\i,\i')$, o\`u
$\i$ est \'egal \`a la diagonale $diag(\g_{1})$ de $\g_{1}\times \g_{1}$,
modulo la conjugaison par la diagnale de $G_{1}\times G_{1}$ .\ste
Ceci
se fait simplement \`a l'aide du Th\'eor\`eme pr\'ec\'edent.\end{rem}
\ste On
remarque que
$\g_{+}=\g_{1}\times \{0\}$, $\g_{-}=\{0\} \times \g_{1} $. On fixe
pour cela une sous-alg\`ebre de Borel $\b_{1}$ de $\g_{1}$. On note
$\b'_{1}$ la sous-alg\`ebre de Borel oppos\'ee, relativement \`a $\j_{1}$.
On pose $\b_{0}=\b_{1}\times \b'_{1}$. Soit ${\cal W}_{1}$ un syst\`eme
de g\'en\'erateurs de Weyl de $\g_{1}$, relativement \`a l'ensemble,
$\Sigma_{1}$, des racines simples de l'ensemble des racines de
$\j_{1}$ dans
$\b_{1}$. On note ${\cal W}$ le syst\`eme de g\'en\'erateurs de Weyl de
$\g$ \'egal \`a $({\cal W}_{1}\times \{0\})\cup (\{0\} \cup {\cal
W}_{1}$. D'apr\`es le Th\'eor\`eme pr\'ec\'edent, il existe une
unique donn\'ee
de Belavin-Drinfeld g\'en\'eralis\'ee, ${\cal BD}=(A, A',\i_{a},
\i_{\a'})$
telle que
$(B,\i,\i')$ soit conjugu\'e \`a ${\cal T}_{ {\cal BD}, {\cal W}}$.
Il est
alors facile de voir que $\i$ est sous $\g$. Alors $\a$ et $\i_{\a}$
sont r\'eduits \`a z\'ero et  $\Gamma_{+}=\Sigma_{1}\times \{0\}$. Par
ailleurs, la conjugaison des triples se traduit par le fait que
l'isomorphisme
$\tau
$ du Lemme 22, de $\g_{1}$ sur $\g_{1}$, est un automorphisme
int\'erieur de $\g_{1}$. Par ailleurs, il pr\'eserve $\j_{1}$ et induit
une permutation, $A$, de l'ensemble $\Sigma_{1}$.  Cette permutation
doit donc  \^etre triviale, i.e.
$A$ est l'identit\'e (cf. [Bou], Chapitre VIII, Paragraphe 5.2). Alors
$\tau
$ est l'identit\'e et la premi\`ere sous-alg\`ebre isotrope de $
{\cal T}_{
{\cal BD}, {\cal W}}$ est la diagonale. L'\'el\'ement de $G$ qui
conjugue
les deux triples stabilise donc la diagonale. C'est donc un \'el\'ement
de la diagonale.
\ste Il  r\'esulte de la discussion pr\'ec\'edente et du
Th\'eor\`eme 5, que
l'ensemble
${\cal T}_{ {\cal BD}, {\cal W}}$, o\`u ${\cal BD}$ d\'ecrit l'ensemble
des donn\'ees de Belavin-Drinfeld g\'en\'eralis\'ees telles que
$A$ est l'identit\'e de
$\Sigma_{1}$, $\i_{\a}=\{0\}$, classifie les triples de Manin pour
$\g$, $(B,\i,\i')$, o\`u $\i$ est \'egal \`a la diagonale $diag(\g_{1})$
 de
$\g_{1}\times \g_{1}$, modulo la conjugaison par la diagnale de
$G_{1}\times G_{1}$. Ceci redonne le r\'esultat de Belavin-Drinfeld
[BD], Th\'eor\`eme 6.1.

\section{Triples de Manin r\'eels pour une alg\`ebre semi-simple
complexe}
\setcounter{equation}{0}
\subsection{}
Soit $\g_{1}$ une alg\`ebre de Lie semi-simple complexe,
soit $\g_{1\R}$ une forme r\'eelle d\'eploy\'ee de $\g_{1}$ et
$\b_{1}$ une sous-alg\`ebre de Borel de $\g_{1}$, complexifi\'ee
d'une sous-alg\`ebre de Borel, $\b_{1\R}$, de  $\g_{1\R}$. Soit
$\j_{1\R}$ une sous-alg\`ebre de Cartan d\'eploy\'ee de $\g_{1\R}$,
contenue dans $\b_{1\R} $ et $\j_{1}$ sa complexifi\'ee. On note
$X\mapsto {\overline X}$ la conjugaison de $\g_{1}$ par rapport
\`a sa forme r\'eelle $\g_{1\R}$. On note $\eta$ l'application de
 $\g_{1}$ dans $\g:=\g_{1}\times \g_{1}$, d\'efinie par :
$$\eta(X)=(X,{\overline X}), \>\> X\in \g_{1}$$
Alors $\eta(\g_{1})$ est une forme r\'eelle de $\g$ et la
conjugaison, $j$, par rapport \`a cette forme r\'eelle v\'erifie :
$$j(X,Y)=({\overline Y}, {\overline X}),\>\> X,Y\in \g_{1}$$
 On note $\b_{0}:=\b_{1}\times \b_{1}$
et $\j_{0}:=\j_{1}\times \j_{1}$ \ste
Si $V$ est un sous-espace vectoriel r\'eel de $\g_{1}$, on
noptera $V_{\C}=\eta(V)+i\eta(V)\subset \g$. On utilisera les
Notations  qui suivent le Lemme 16 pour $\g$.
\ste Si $B$ est une forme $\R$-bilin\'eaire
invariante  sur $\g_{1}$, on note $B_{\C}$, l'unique forme
$\C$-bilin\'eaire
invariante  sur $\g$, telle que :
\beq B_{\C}(\eta(X), \eta (X') )=B(X,X'),
\>\> X,X'\in \g_{1} \eeq
On voit ais\'ement que, si $\s$, est un id\'eal simple de $\g_{1}$,
et si la restriction de  $B$ \`a $\s$ est \'egal \`a $Im \>\lambda
K_{\s}$, pour $\lambda\in \C$, on a :
\beq B_{\C}((X,Y), (X',Y'))=\lambda K_{\s}(X,X')-{\overline
\lambda} K_{\s}(Y,Y'),\>\>X,X',Y,Y'\in \s \eeq
La d\'emonstration de la Proposition suivante est imm\'ediate.
\begin{prop}
On fixe une une forme de Manin r\'eelle sur $\g_{1}$.
L'application qui \`a un un triple de Manin   r\'eel dans $\g_{1}$,
$(B,\i,\i')$, associe $(B_{\C},\i_{\C},\i'_{\C})$, est une
bijection entre l'ensemble des triples de Manin r\'eels de
$\g_{1}$, associ\'es \`a $B$, et l'ensemble des triples de Manin
complexes de $\g$, associ\'es \`a $B_{\C}$, et pour lesquels chacune
des sous-alg\`ebres Lagrangiennes est stable par $j$.\ste
Cette bijection transforme triples fortement standard,
relativement \`a $\b_{1}$, $\j_{1}$, en triples fortement
standard relativement \`a $\b_{0}$, $\j_{0}$. En outre elle transforme
les triples conjugu\'es par $G_{1}$ en triples conjugu\'es par
$\eta (G_{1})$
\end{prop}
{\bf Hypoth\`ese }\ste
{\em On suppose que, pour tout d\'eal simple de $\g_{1}$,
la restriction de  $B$ \`a $\s$ est \'egal \`a $Im \>\lambda
K_{\s}$, pour un  $\lambda$ r\'eel.}\ste
On note $\g_{1+}$ la somme des id\'eaux simples $\s$ de
$\g_{1}$, tels que  la restriction de  $B$ \`a $\s$ est \'egal \`a $Im
\>\lambda K_{\s}$, pour $\lambda\in \R^{+}$. On d\'efinit de
m\^eme $\g_{1-}$. Alors, au vu de  (4.2), on a, pour la forme
$B_{\C}$ : \beq \g_{+}= \g_{1+}\times \g_{1-},\>\>
\g_{-}= \g_{1-}\times \g_{1+}\eeq
On note ${\tilde R}_{1*}$  l'ensemble des
racines de $\j_{1}$ dans $\g_{1*}$ , o\`u $*$ vaut + ou -. Alors, avec
les Notations qui suivent le Lemme 16,  on a :
 $$ {\tilde R}_{+}= ({\tilde R}_{1+}\times \{0\})\cup (\{0\}\times
  {\tilde R}_{1-})$$
On note $\Sigma_{1+}$ (resp. $\Sigma_{1-}$) les racines simples
de $\j_{1+}:=\j_{1}\cap \g_{1+}$ dans $\b_{1}\cap \g_{1+}$
(resp. $\j_{1-}:=\j_{1}\cap \g_{1-}$ dans $\b'_{1}\cap
\g_{1-}$, o\`u $\b'_{1}$ est la sous-alg\`ebre de Borel de
$\g_{1}$ oppos\'ee \`a $\b_{1}$, relativement \`a $\j_{1}$)), qu'on
identifie \`a des racines de $\j_{1}$ dans $\g_{1}$. Alors on a :
\beq \Sigma_{+}=(\Sigma_{1+}\times \{0\})\cup (\{0\}\times
\Sigma_{1-}),\>\>
\Sigma_{-}=((-\Sigma_{1-})\times \{0\})\cup (\{0\}\times
(-\Sigma_{1+}))\eeq
Soit ${\cal W}_{1}$ un syst\`eme de g\'en\'erateurs de Weyl de
$\g_{1}$, relativement \`a $\Sigma_{1}$, dont tous les \'el\'ements
sont dans $\g_{1\R}$. Soit ${\cal W}$ le syst\`eme de g\'en\'erateurs
de Weyl de
$\g$, relativement \`a
$\Sigma:=\Sigma_{+}\cup\Sigma_{-}$, et d\'efini comme suit, en tenant
compte de (4.4) . \ste  Si $ \alpha\in
\Sigma_{1+}$ on pose :
\beq H_{(\alpha,0)}=(H_{\alpha}, 0),\>\>
X_{(\alpha,0)}=(X_{\alpha}, 0),\>\> Y_{(\alpha,0)}=(Y_{\alpha},
0)\eeq
\beq  H_{(0,-\alpha)}=(0,-H_{\alpha}),\>\>
X_{(0,-\alpha)}=(0,Y_{\alpha}),\>\>
Y_{(0,-\alpha)}=(0,X_{\alpha})
\eeq
Si $ \alpha\in
\Sigma_{1-}$, on pose  :
\beq  H_{(0,\alpha)}=(0,H_{\alpha}),\>\>
X_{(0,\alpha)}=(0,X_{\alpha}),\>\>Y_{(0,\alpha)}=(0,Y_{\alpha})
\eeq
\beq H_{(-\alpha,0)}=(-H_{\alpha}, 0),\>\>
X_{(-\alpha,0)}=(Y_{\alpha}, 0),\>\> Y_{(-\alpha,0)}=(X_{\alpha},
0)\eeq
On notera  $\alpha\mapsto \alpha ^{f}$ l'\'echange des
facteurs dans  $\j_{1}^{*}\times\j_{1}^{*}=\j_{0}^{*}$. On fait
de m\^eme dans $J_{0}$. \begin{prop}\ste On fixe  une forme de
Manin r\'eelle, $B$, sur $\g_{1}$. Soit ${\cal
BD}=(A,A',\i_{\a},\i_{\a'})$ une donn\'ee de Belavin-Drinfeld
g\'en\'eralis\'ee pour $\g$ et $B_{\C}$.\ste  (i) Pour   $t\in J_{0}$,
le
triple de Manin $(B_{\C},\i,\i')={\cal T}_{{\cal BD}, t{\cal W}}
$ est le complexifi\'e d'un triple r\'eel de $\g_{1}$, relativement
\`a $B$, si et seulement si :\ste  1) l'application $\alpha\mapsto
-\alpha^{f}$ induit  une bijection de $\Gamma_{+}$ sur
$\Gamma_{-}$ et l'on a, en prolongeant $A$ par $\R$-lin\'earit\'e :
$$(A\alpha)^{f}= A^{-1}(\alpha^{f}),
\>\>\alpha\in \Gamma_{+}$$
2) A' v\'erifie des conditions similaires.\ste
3) $\i_{\a} $ et $\i_{\a'}$ sont stables par $j$
4)$$ Posant \>\> u:={\overline t}(t^{f})^{-1},
\>\>on \>\> a\>\>u^{\alpha}=u^{A\alpha },
\alpha\in \Gamma_{+}
$$
5) L'\'el\'ement  $t$ de $J_{0}$ v\'erifie des conditions similaires
relativement \`a
$A'$.
\
\ste (ii) On fixe  ${\cal BD} $ et $t$ v\'erifiant les conditions
ci-dessus. Soit $t_{1}\in J_{1}$ et $t'=(t_{1},{\overline
t_{1}})t$. Alors ${\cal BD} $ et $t'$ v\'erifient les conditions
ci-dessus, et les triples complexes ${\cal T}_{{\cal BD}, t{\cal
W}}$, ${\cal T}_{{\cal BD}, t'{\cal W}}$ sont les complexifi\'es de
triples r\'eels de $\g_{1}$, conjugu\'es par $t_{1}$.
\end{prop}
\dem
Si l'automorphisme $\R$-lin\'eaire de $\g$, $j$, laisse $\i$
invariant, il laisse invariant le radical nilpotent $\n$ de
$\i$, donc aussi $\p$, qui est le normalisateur de $\n$. Par
ailleurs $\j_{0}$ est contenu dans $\p$ et est invariant par
$j$. Donc $j(\l)=\l$, d'o\`u l'on d\'eduit $j(\m)=\m$.\ste
Alors $\i=\h\oplus\i_{\a}\oplus\n$ est invariant par $j$ si et
seulement si :
\beq j(\h)=\h,\>\>  j(\i_{\a})=\i_{\a},\>\> j(\n)=\n \eeq
La seconde \'egalit\'e de l'\'equation pr\'ec\'edente  conduit
\`a 3). \ste
L'\'egalit\'e  (cf (3.19)) :
\beq \h=\f\oplus (\oplus_{\alpha\in R_{+}}\h_{\alpha})\eeq
et la stabilit\'e de $\j_{0}$, et de son orthogonal pour la forme
de Killing  par $j$, montre que la premi\`ere \'egalit\'e de
(4.9) implique :
$$j(\f)=\f$$
Mais $\f$ est engendr\'e par les $H_{\alpha}+H_{A\alpha}$,
$\alpha \in \Gamma_{+}$. De plus $j(H_{\alpha}+H_{A\alpha})$
est \'egal \`a $H_{\alpha^{f}}+H_{(A\alpha)^{f}}$. Ce dernier doit
\^etre une combinaison lin\'eaire de  $H_{\beta}+H_{A\beta}$,
 $\beta \in \Gamma_{+}$. Mais on a : $$H_{\beta}\in \g_{+}, \>\>et
\>\> H_{A\beta}\in \g_{-},\>\> H_{\alpha^{f}}\in \g_{-},\>\>
H_{(A\alpha)^{f}}\in \g_{+}$$
car $f$ \'echange $\g_{+}$ et $\g_{-}$ (voir (4.3)).
Comme $f$ envoie chaque \'el\'ement de $\Sigma$ sur l'oppos\'e d'un
\'el\'ement de $\Sigma$
(cf. (4.4)), la projection, sur
$\g_{+}$ et
$\g_{-}$, de l'\'ecriture  de  $H_{\alpha^{f}}+H_{(A\alpha)^{f}}$ dans
la base de $\f$, montre qu'il existe $\beta
\in \Gamma_{+}$ telle que :
\beq (A\alpha)^{f}=-\beta,\>\> (\alpha)^{f}=-A\beta\eeq
Ceci implique imm\'ediatement la condition 1).
Comme :
$$j(\g^{\alpha})=\g^{\alpha ^{f}},\>\>\alpha\in {\tilde
R}$$
au vu de (4.10) et (4.11),la stabilit\'e de $\h$ par $j$ implique
alors  :
$$j(\h_{\alpha})=\h_{-\beta }$$
o\`u $\alpha\in \Gamma_{+}$ et $ \beta$ est comme ci-dessus. Mais, la
d\'efinition de $\i$ (cf. Proposition 8 et Lemme 22) montre que
 $\h_{\alpha}$ a pour base :
$$U_{\alpha}:=t^{\alpha}X_{\alpha}+t^{A\alpha}X_{A_{\alpha}}$$
et $\h_{-\beta }$ a pour base:
$$V_{\beta }:=t^{-\beta }Y_{\beta}+t^{-A\beta }Y_{A_{\beta}}$$
Par ailleurs la d\'efinition de $j$ et celle de ${\cal W}$
montrent que : $$j(X_{\alpha})=Y_{-\alpha^{f}}, \>\>
j(X_{A\alpha})=Y_{-(A\alpha)^{f}}$$
Donc, on a :
$$j(U_{\alpha}):={\overline
t}^{\alpha}Y_{-\alpha^{f}}+{\overline
t}^{A\alpha}Y_{-(A_{\alpha})^{f}}$$
 L'\'ecriture de la proportionnalit\'e
de $ j(U_{\alpha})$ \`a $V_{\beta}$, conduit \`a :
$${\overline t}^{\alpha}t^{A\beta}={\overline
t}^{A\alpha}t^{\beta}, \alpha \in \Gamma_{+}$$
En tenant compte de (4.11), ceci implique la condition 4).
On proc\`ede
de m\^eme pour $\i'$.\ste  Ce qui pr\'ec\`ede montre que les
conditions 1)
\`a 5) sont n\'ecessaires pour que le triple donn\'e soit le
complexifi\'e
d'un triple r\'eel.
R\'eciproquement, montrons que si ces conditions sont satisfaites
le triple donn\'e est bien le complexifi\'e d'un
triple r\'eel. En  effet, la deuxi\`eme \'egalit\'e de (4.9) est alors
satisfaite. \ste
Comme les racines de $\j_{0}$ dans $\n$ sont celles de
$\j_{0}$ dans $\b_{1}\times \b_{1}$  qui ne sont pas
combinaison lin\'eaires d'\'el\'ements de $\Gamma_{+}\cup
\Gamma_{-}$, on d\'eduit la troisi\`eme \'egalit\'e de (4.9) de la
condition 1) et (4.11). Il reste \`a v\'erifier la premi\`ere
\'egalit\'e de
(4.9). D'abord, il r\'esulte  de 1) et de la discussion ci-dessus que :
\beq
\f\>\> est\>\> stable \>\>par\>\> j\eeq   Par ailleurs les  conditions
1) et 4),  et la discussion ci-dessus montre que : \beq
j(\h_{\alpha})=\h_{-(A\alpha )^{f}},\>\> \alpha \in
\Gamma_{+}\eeq On montre de m\^eme que :
\beq j(\h_{-\alpha})=\h_{(A\alpha )^{f}},\>\> \alpha \in
\Gamma_{+}\eeq
Mais $\h$, qui est isomorphe \`a sa projection dans $\g_{+}$, est
engendr\'ee par les $\h_{\alpha}$, $\h_{-\alpha}$, $\alpha \in
\Gamma_{+}$. Alors (4.12) \`a (4.14), joints \`a 1) montrent
que
$\h$ est stable par $\i$, ce qui ach\`eve de prouver que $\i$
est stable par $j$. On proc\`ede de m\^eme pour $\i'$. Ceci ach\`eve
de prouver (i). (ii) est une cons\'equence imm\'ediate de la
Proposition pr\'ec\'edente \qed

On se fixe une donn\'ee de Belavin-Drinfeld g\'en\'eralis\'ee, qui
v\'erifie les propri\'et\'es 1) \`a 4) de la Proposition p\'ec\'edente.
On
introduit $C=''A^{-1}A'\> ''$, comme dans la D\'efinition 5. On
note $\Gamma_{0}:=dom \> C\cup Im\> C$. Si $\alpha \in
\Gamma_{0}$, on note
${\cal C}(\alpha)$, l'ensemble des $\beta \in \Gamma_{0}$ tels qu' il
existe $n\in \N$ tel que  :
$$\alpha, C\alpha, \dots,C^{n-1}\alpha \in \Gamma_{0}, \>\>
et \>\>\beta=C^{n}\alpha $$
ou  :
$$\beta, C\beta , \dots, C^{n-1}\beta \in \Gamma_{0}, \>\>
et\>\> \alpha=C^{n}\beta $$
Si $\alpha\in \Sigma_{+}$ n'appartient pas \`a $\Gamma_{0}$, on
pose  ${\cal C}(\alpha)=\{\alpha\}$. Suivant Panov [P1], on
appelle les ${\cal C}(\alpha)$ des chaines.
Il est clair que les chaines sont disjointes ou confondues, et
forment une partition de $\Sigma_{+}$.  On appellera
$C$-\'equivalence la relation d'\'equivalence sur $\Sigma_{+}$, dont
les classes d'\'equivalence sont les chaines. On dira que, pour deux
\'el\'ements distincts $\alpha$ et $\beta$, de $\Sigma_{+}$, $\alpha$
est $C$-li\'e \`a $\beta$  si $\alpha  \in dom\> C$ et $\beta=C\alpha$.
On d\'efinit, pour $\alpha\in \Gamma_{0}$ : $$\check{{\cal C}}
(\alpha):= \{A^{-1}(-\beta^{f})\vert \beta \in {\cal
C}(\alpha)\cap \Gamma_{+}\}\cup  \{A'^{-1}(-\beta^{f})\vert
\beta \in {\cal C}(\alpha)\cap \Gamma'_{+}\}$$
Si $\alpha\notin \Gamma_{0}$, on pose
$\check{\cal  C}(\alpha)
:={\cal C}(\alpha)$. \begin{lem}\ste  (i)
Pour tout $\alpha\in \Sigma_{+}$, $\check{{\cal C}}(\alpha)$ est
de la forme ${\cal C}(\beta)$, et l'op\'eration  $\check{}$ est une
involution de l'ensemble des chaines. \ste (ii) Si $t$ v\'erifie les
conditions 4)
 et 5)
de la Proposition pr\'ec\'edente, on a, avec les notations de
celle-ci, pour tout $\alpha\in \Gamma_{0}$: $$  u^{\beta}=
u^{\beta'}, \beta,\>\> \beta'\in {\cal C}(\alpha)$$
$$ u^{\beta }={\overline {u^{\alpha}}}, \>\>
\beta \in \check{{\cal C}}(\alpha)$$ \end{lem}
\dem
Montrons (i). Il suffit d'\'etudier le cas $\alpha\in \Gamma_{0}$.
Soit $\alpha\in \Sigma_{+}$. On note ${\cal A}={\cal
C}(\alpha)\cap \Gamma_{+}$,  ${\cal A}'={\cal
C}(\alpha)\cap \Gamma'_{+}$. De la d\'efinition des chaines,
il r\'esulte que, pour toute paire d'\'el\'ements distincts $\beta$,
$\beta'$, de ${\cal A}$,  il existe $n\in \N^{*}$,
$\beta_{1},\dots, \beta_{n+1}\in {\cal A}$, o\`u, pour
$i=1,\dots ,n$,  $\beta_{i}$ est $C$-li\'e \`a $\beta_{i+1}$,
avec $\{ \beta,\beta '\}=\{ \beta_{1},\beta_{n+1}\}$.    \ste On
remarque que les conditions 1) et 2) impliquent, par un calcul
imm\'ediat, que : $$Si\>\> \beta, \beta' \in {\cal A} \>\> sont\>\>
C-li\acute{e}s,\>\> A^{-1}(-\beta ^{f})\>\> et\>\>
 A^{-1}(-\beta'^{f})\>\> sont\>\>
C-li\acute{e}s $$
On a le m\^eme \'enonc\'e pour ${\cal A}'$. \ste
De ce qui pr\'ec\`ede, il r\'esulte  que les \'el\'ements de
$A^{-1}{\cal A}$ (resp. $A'^{-1}{\cal A}'$) sont
$C$-\'equivalents entre eux.   Pour achever de prouver (i), il suffit
de traiter le cas o\`u ${\cal A}$ et ${\cal A'}$ sont non vides et
d'intersection vide. De la d\'efinition des chaines , et du
fait que $dom\>C$ (resp. $Im\> C$) est un sous-ensemble de
$\Gamma'_{+}$ (resp. $\Gamma_{+}$), cela n'est possible que si
${\cal C}(\alpha)$ n'a que deux \'el\'ements et ${\cal A}$ (resp.
${\cal A'}$)  un \'el\'ement $\beta$ (resp. $\beta'$), o\`u $\beta' $
et $\beta$ sont $C$-li\'es. En outre :
\beq {\beta }\in \Gamma_{+}\setminus \Gamma'_{+},\>\>
{\beta '}\in \Gamma'_{+}\setminus \Gamma_{+}\eeq
 Mais alors $A'\beta'$ et
$A\beta$ sont \'egaux et on note $\gamma$ leur valeur commune .
Utilisant les conditions 1) et 2) de la
Proposition 10, on en d\'eduit que :
$$A^{-1}(-\beta^{f})=A'^{-1}(-\beta'^{f})=-\gamma^{f},\>\>
\check{{\cal C}} (\alpha) =\{-\gamma^{f}\}$$
Montrons que la classe de $C$-\'equivalence de $-\gamma^{f}$ est
r\'eduite \`a un \'el\'ement. Pour cela il suffit de voir que
$-\gamma^{f}$ n'est $C$-li\'e \`a aucun \'el\'ement  et qu'aucun
\'el\'ement ne lui est $C$-li\'e. Raisonnons par l'absurde et
supposons par exemple que $\delta$ soit $C$-li\'e \`a
$-\gamma^{f}$. On a alors :
\beq \delta \in  \Gamma'_{+}\>\> et  \>\>
A(-\gamma^{f})=A'\delta\eeq  On d\'eduit des conditions 1), 2) de la
Proposition 10 :
\beq A^{-1}\gamma=A'^{-1}(-\delta^{f})\eeq
D'apr\`es (4.16) et la condition 1), $A'^{-1}(-\delta^{f}$ est un
 \'el\'ement de $\Gamma'_{+}$. Par ailleurs, d'apr\`es la d\'efinition
de $\gamma$, $A^{-1}\gamma$ est \'egal \`a $\alpha$. Joint \`a  (4.17),
cela montre que $\alpha$ est \'el\'ement de $\Gamma'_{+}$, ce qui
contredit (4.15). Donc aucun \'el\'ement n'est $C$-li\'e \`a
$-\gamma^{f}$.\ste
On montre de m\^eme que $-\gamma^{f}$ n'est li\'e \`a aucun
\'el\'ement. Cecci ach\`eve de prouver (i). \ste
Si $\beta$
est $C$-li\'e \`a $\beta'$, on a $A\beta'=A'\beta$. Par ailleurs,
d'apr\`es les conditions 4) et 5) de la proposition 10, on a :$$
u^{\beta}=u^{A'\beta},\>\>u^{\beta'}=u^{A\beta'}$$ La
premi\`ere \'egalit\'e de (ii) en r\'esulte.\ste
Pour la deuxi\`eme \'egalit\'e, on remarque, que si $\alpha\in
\Gamma_{+}$, on a, d'apr\`es la condition 4) de la Proposition
10 :
$$u^{\alpha}=u^{A\alpha}$$
Mais, il r\'esulte de la d\'efinition de $u$ que :
\beq u^{f}={\overline u}^{-1}\eeq
On en d\'eduit le r\'esultat voulu.\qed
 Le th\'eor\`eme suivant  a \'et\'e sugg\'er\'e par  un  r\'esultat de A. Panov
(cf. [P1], Th\'eor\`eme 6.13) et la preuve que nous en  donnons plus
loin.
 \begin{theo}
Soit $B$ une forme de Manin r\'eelle sur $\g_{1}$. Tout triple de
Manin r\'eel dans $\g_{1}$, relativement \`a $B$, est conjugu\'e par un
\'el\'ement de
$G_{1}$ \`a un triple fortement standard dont le complexifi\'e est de la
forme
${\cal T}_{{\cal BD}, t{\cal W}}$, o\`u $t$ est un \'el\'ement de $J_{0}$
et ${\cal BD}=(A,A',\i_{\a},\i_{\a'})$ une donn\'ee de
Belavin-Drinfeld g\'en\'eralis\'ee pour $\g$ et $B_{\C}$, qui v\'erifie,
outre les conditions 1) \`a 5) de la Proposition 10 : \ste
1)
$$u^{2}=1$$  2) Pour tout $\alpha\in
Gamma^{+}\setminus \Gamma_{0}$ : $u_{1}^{\alpha}=1$.
\ste
3) $t=(v_{1}, 1)$, o\`u $u_{1}\in J_{1}$.\ste

Alors $v_{1}^{2}=1$
\end{theo}

\begin{rem}
Pour une classe de conjugaison sous $G_{1}$ de triples r\'eels
de $\g_{1}$, relativement \`a $B$, la donn\'ee ${\cal BD}$ est uniquement
d\'etermin\'ee, d'apr\`es le Th\'eor\`eme 5  et il n'y a qu'un nombre
fini de
choix possibles pour
$t$, car les \'el\'ements de carr\'e 1 de $J_{1}$ sont en nombre fini.
\end{rem}
\dem
Tenant compte du fait que $u$ n'est pas chang\'e par la
mutiplication de $t$ par un \'el\'ement de la forme
$(t_{1},{\overline t}_{1})$, on voit qu 'il suffit  de trouver
$t$ satisfaisant toutes les conditions \`a l'exception de 3).
D'apr\`es la Proposition 10, il existe ${\underline t}$ et ${\cal
BD}$ satisfaisant les conditions 1) \`a 5) de celle-ci (en y
changeant  $t$ en ${\underline t}$).
On v\'erifie ais\'ement que  si $t'\in J_{0}$ satisfait:
\beq t'^{\alpha}=t'^{A\alpha}, \alpha \in \Gamma_{+}\eeq
\beq t'^{\alpha}=t'^{A'\alpha}, \alpha \in \Gamma'_{+}\eeq
 on a : $${\cal
T}_{{\cal BD}, {\underline t}{\cal W}}={\cal
T}_{{\cal BD}, {\underline t}t'{\cal W}}$$
Il reste \`a choisir $t'$ v\'erifiant  (4.19) et (4.20) de telle sorte
que $t={\underline t}t'$ v\'erifie les propri\'et\'es voulues.
\ste On choisit un sous-ensemble $\Theta $ de $\Sigma_{+}$, tel que
toute classe de $C$-\'equivalence soit de la forme ${\cal
C}(\alpha)$ ou $\check{{\cal
C}} (\alpha)$ pour un unique $\alpha\in \Theta$.\ste
On remarque que si $A\alpha= A'\beta$, $\alpha$ et $\beta$ sont
$C$-\'equivalents.
On caract\'erise alors $t'$ par $t^(\alpha)$, $\alpha \in \Sigma$.
 On choisit  pour tout
$\alpha\in \Sigma_{+}$, une racine carr\'ee  $z_{\alpha}$ de
${\overline {({\underline u}^{\alpha})}}^{-1}$, o\`u ${\underline u}={\overline
{\underline t}} ({\underline t}^{f})^{-1}$. On pose alors, pour
$\alpha\in
\Theta$  :
\beq t'^{\beta}=z_{\alpha}
,
\>\>si
\>\>
\beta
\in {\cal C}(\alpha), \>\>et \>\> si\>\> \alpha\notin
\Gamma_{0}\>\>\eeq
\beq t'^{\beta}=z_{\alpha} \>\>(resp.\>\> {\overline {z_{\alpha}}}),
\>\>si
\>\>
\beta
\in {\cal C}(\alpha) (resp.\>\>\beta  \in \check{{\cal C}}
(\alpha) )\>\> et
\>\>si  \>\> {\cal C}(\alpha)\not=\check{{\cal
C}} (\alpha) \eeq
\beq t'^{\beta}=\vert u_{\alpha}\vert ^{-1/2}\>\> si\>\> \alpha\in
\Gamma_{0}\>\>  est\>\> tel \>\> que \>\> {\cal C}(\alpha) =
\check{{\cal
C}} (\alpha), \>\>et \>\>\beta \in {\cal C}(\alpha) \eeq Ceci
d\'etermine $t'^{\beta}$ pour $\beta \in
\Sigma_{+}$, et l'on pose :
\beq t'^{{\beta}^{f}}={\overline {(t'^{\beta})^{-1}}}, \beta\in
\Sigma_{+} \eeq
Ainsi $t'$ est enti\`erement caract\'eris\'e par les relations (4.21)
\`a
(4.24). \ste Il faut voir que
$t'$ v\'erifie (4.19) et(4.20). Soit $\beta\in \Gamma_{+}$. On suppose
$\beta\in {\cal C}(\alpha)$, $\alpha\in \Theta$. D'apr\`es (4.24), on a
:
$$t'^{A\beta}={\overline t'^{-(A\beta)^{f}}}$$
Mais, d'apr\`es la condition 1) de la Proposition 10, et la
d\'efinition de $\check{{\cal
C}} (\alpha)$, on a :
$-(A\beta)^{f}=A^{-1}(-\beta^{f})\in \check{{\cal
C}} (\alpha)$, donc
, d'apr\`es (4.22) et (4.23), on a :
 $${\overline {t'^{-(A\beta)^{f}}}}=t'^{\beta}$$
On traite de m\^eme le cas o\`u $\beta \in \check{{\cal
C}} (\alpha)$, et
alors
$\alpha\in \Theta$. Ceci prouve que (4.19) est v\'erifi\'e. On prouve
(4.20) de la m\^eme mani\`ere.
\ste Calculons $v_{\beta}:={\overline
t'^{\beta}}(t'^{-1})^{{\beta}^{f}}$, $\beta\in \Sigma_{+}$
 A l'aide de
(4.22) \`a (4.24) et du Lemme 23, on voit que pour $\alpha\in \Theta$
:
$$v_{\beta}=u^{-\beta},\>\> si\>\> \beta
\in {\cal C}(\alpha)\cup \check{{\cal
C}} (\alpha) , \>\>et \>\> si\>\> \alpha\notin
\Gamma_{0}\>\> ou\>\> si\>\> {\cal C}(\alpha)\not=\check{{\cal
C}} (\alpha) $$
$$ v_{\beta}=\vert u\vert ^{-\beta},\>\>   si\>\> \alpha\in
\Gamma_{0}\>\>  est\>\> tel \>\> que \>\> {\cal C}(\alpha)=\check{{\cal
C}} (\alpha), \>\>et \>\>\beta \in {\cal C}(\alpha) $$
Par ailleurs il r\'esulte du Lemme 23 :
$$u^{\beta}\>\> est\>\> r\acute{e}el\>\>  si\>\> \alpha\in
\Gamma_{0}\>\>  est\>\> tel \>\> que \>\> {\cal C}(\alpha)=\check{{\cal
C}} (\alpha), \>\>et \>\> si\>\>\beta \in {\cal C}(\alpha) $$
Alors il est clair que $t:=t'{\underline t}$ v\'erifie la  condition
 2) du Th\'eor\`eme et l'\'el\'ement $u$ correspondant v\'erifie
:$$u^{\alpha}=1 \>\> ou -1, \alpha\in \Sigma_{+}$$
$$u^{f}={\overline u}^{-1}$$
Il en r\'esulte que $u^{2}=1$ comme d\'esir\'e.  \qed
\subsection{Une autre d\'emonstration d'un r\'esultat d'A. Panov}
\ste On suppose maintenant que $\g_{1}$ est une alg\`ebre de Lie
complexe simple et on pose
$B_{1}=Im\> K_{\g_{1}}$. On fixe une forme r\'eelle $\h_{1}$ de
$\g_{1}$ et
$\sigma_{1}$ la conjugaison de $\g_{1}$ par rapport \`a $\h_{1}$. \ste
On note $G_{1}$ le groupe adjoint de $\g_{1}$ (et non son
recouvrement universel) et, si $\e$ est une sous-alg\`ebre de Lie
r\'eelle
de $\g_{1}$, on note $E$ le sous-groupe analytique de $G_{1}$,
d'alg\`ebre de Lie $\e$. \ste On s'int\'eresse aux triples de Manin
r\'eels
de
$\g_{1}$,
$(B_{1},\i_{1},\i'_{1})$ tels que $\i_{1}$ soit \'egal \`a $\h_{1}$,
qu'on appelle ''$\h_{1}$-triple'' , \`a conjugaison pr\`es par les
\'el\'ement
de
$G_{1}$, ou, ce qui revient au m\^eme, par ceux de
$G_{1}^{\sigma_{1}}$, qui est l' ensemble des \'el\'ements de $G_{1}$
commutant \`a $\sigma_{1}$. Ces triples d\'ecrivent les structures de
big\`ebres de Lie sur $\h_{1}$, dont le double est isomorphe \`a
$\g_{1}$ (cf. [P1]).
\ste Soit
$\f_{1}$ une \sac fondamentale  de $\h_{1}$ et soit $\j_{1}$
la complexifi\'ee de $\f_{1}$ dans $\g_{1}$. On choisit
 une sous-alg\`ebre de Borel, $\b_{1}$, de $\g_{1}$, contenant $\j_{1}$
et telle que $\sigma_{1}(\b_{1})$ soit \'egal \`a la sous-alg\`ebre
 de Borel de $\g_{1}$, $\b'_{1}$, oppos\'ee \`a $\b_{1}$,
relativement \`a $\j_{1}$. On note $\Sigma_{1}$ l'ensemble  des racines
simples de
$\j_{1}$ dans $\b_{1}$.\ste
On note $\theta$ l'involution de $\Sigma_{1}$, caract\'eris\'ee par
:
\beq \sigma_{1}(\g_{1}^{-\alpha})=\g_{1}^{\theta(\alpha)} \eeq
En calulant de deux mani\`eres diff\'erentes $t(\sigma_{1}(X))$, pour
$X\in \g^{\alpha}$, on trouve  :
\beq  {\overline 
{(t^{\sigma_{1}})^{\alpha}}}={\overline {t^{-\theta(\alpha)}}}\eeq Il
est facile de voir  qu'on peut choisir un syst\`eme de g\'en\'erateurs
de Weyl de
$\g_{1}$, ${\cal W}_{1}$,  relativement \`a
$\Sigma_{1}$, tel que :
\beq   \sigma_{1}(H_{\alpha})=H_{-\theta(\alpha)},
\>\> \sigma_{1}(X_{\alpha})=\varepsilon_{\alpha}Y_{\theta(\alpha)},
\>\>
\sigma_{1}(Y_{\alpha})=\varepsilon_{\alpha}X_{\theta(\alpha)}\eeq o\`u :
\beq \varepsilon_{\alpha}=1\>\> ou\>\>-1, \>\> et \>\>
\varepsilon_{\alpha}=1\>\> si \>\> \theta(\alpha)\not= -\alpha\eeq
On note $\g :=\g_{1}\times
\g_{1}$,
$\j_{0}=\j_{1}\times
\j_{1}$,
$\b_{0}=\b_{1}\times \b'_{1}$. On note $\eta_{1}$
l'aplication  de
$\g_{1}$ dans $\g$ d\'efinie par :
$$\eta_{1}(X)=(X,\sigma_{1}(X)),\>\> X\in \g_{1}$$
dont l'image par $\eta_{1}$ est une forme r\'eelle de $\g$. On note
$j_{1}$ la conjugaison de $\g$ par rapport \`a cette forme r\'eelle, qui
v\'erifie :
$$j_{1}(X,Y)=(\sigma_{1}(Y),\sigma_{1}(X)),\>\> X,Y\in \g_{1}$$
On note ${\cal W}=({\cal W}_{1}\times \{0\})\cup( \{0\} \times {\cal
W}_{1}$). On note $B$ la forme de Manin sur $\g$, \'egale \`a
$K_{\g_{1}}$ sur le premier facteur et \`a $-K_{\g_{1}}$ sur le
deuxi\`eme facteur facteur. Alors $\g_{+}=\g_{1}\times \{0\}$,
$\g_{-}=  \{0\}\times \g_{1}$ et $\Sigma_{+}= \Sigma_{1}$,
$\Sigma_{-}=
 \Sigma_{1}$\ste
Si $V$ est un sous-espace vectoriel r\'eel de $\g_{1}$, on notera :
$$V_{\C}:=\eta_{1}(V)+i\eta_{1}(V)$$
On remarque que $\h_{1\C}$ est \'egal \`a la diagonale,
$diag(\g_{1})$, dans $\g_{1}\times \g_{1}$.\ste Si
${\cal T}_{1}=(B_{1},\i_{1},\i'_{1})$ est un triple de Manin r\'eel dans
$\g_{1}$, on appelle triple complexifi\'e de ${\cal T}_{1}$, le triple
de Manin complexe dans $\g$, ${\cal T}_{1\C}=
(B,\i_{1\C},\i'_{1\C})$.
Alors $\i_{1\C}$ et $\i'_{1\C}$ sont  stables par $j_{1}$. De
plus si on a  deux triples r\'eels dans $\g_{1}$, comme ci-dessus, ils
sont conjugu\'es par un \'el\'ement de $G_{1}$ si et seulement si leurs
complexifi\'es sont conjugu\'es par un \'el\'ement de
$\eta_{1}(G_{1})=\{(g_{1},g_{1}^{\sigma_{1}})\vert g_{1}\in G_{1}\}$.
On remarque que si deux sous-espaces vectoriels complexes de $\g$
sont conjugu\'es par un \'el\'ement de $\eta_{1}(G_{1})$, si l'un est
stable
par $j_{1}$, l'autre l'est aussi.\ste
De plus le complexifi\'e d'un triple fortement standard relativement \`a
$\b_{1}$,
$\j_{1}$, est fortement standard relativement \`a $\b_{0}$,
$\j_{0}$, car il est facile de voir que le complexifi\'e  d'un
ant\'ec\'edent est l'ant\'ec\'edent du complexifi\'e.\ste Comme tout
triple de
Manin r\'eel dans
$\g_{1}$ est conjugu\'e  par un \'el\'ement de
$G_{1}$, \`a un triple fortement standard, d'apr\`es le Th\'eor\`eme 4,
on a
imm\'ediatement :
\begin{lem}\ste
(i) Si ${\cal T}_{1}$ est un ''$\h_{1}$-triple", son complexifi\'e est
conjugu\'e par un \'el\'ement  de
$\eta_{1}(G_{1})$ \`a un triple fortement standard ${\cal
T}=(B,\i,\i')$ tel que : \ste
1) $\i$ est conjugu\'e par un \'el\'ement de $\eta_{1}(G_{1})$ \`a
$diag(\g_{1})$.\ste
2) $\i$ et $\i'$ sont stables par $j_{1}$.\ste
(ii) Si ${\cal
T}=(B,\i,\i')$ est un triple  de Manin fortement standard v\'erifiant
1), 2) , il existe un ''$\h_{1}$-triple'', unique \`a conjugaison sous
$G_{1}$ pr\`es (o\`u plut\^ot $ G_{1}^{\sigma_{1}}$ pr\`es) tel que son
complexifi\'e  soit conjugu\'e \`a ${\cal T}$ par un \'el\'ement de
$\eta_{1}(G_{1})$.  Deux triples fortement standard, associ\'es \`a
$B$, v\'erifant 1), 2),  conjugu\'es sous
$\eta_{1}(G_{1})$, conduisent \`a la m\^eme classe de conjugaison sous $
G_{1}^{\sigma_{1}}$ de ''$\h_{1}$-triple''. \end{lem}
\begin{lem}\ste
(i)Tout triple de Manin complexe, fortement standard dans $\g$ pour
$\b_{0}$, $\j_{0}$, associ\'e \`a
$B$, v\'erifiant les conditions 1), 2), du Lemme pr\'ec\'edent est
conjugu\'e par un \'el\'ement de $\eta_{1}(G_{1})$ \`a un triple ${\cal
T}_{{\cal BD}, (1,t)({\cal W})}$, v\'erifiant les m\^emes
conditions, o\`u
${\cal BD}$ est une donn\'ee de Belavin-Drinfeld
$(A,\i_{\a},A',\i'_{\a'})$, relativement \`a $B$ et
$t$ est un \'el\'ement de $J_{1}$.\ste
(ii) Le fait que  ${\cal
T}_{{\cal BD},(1,t)({\cal W})}$ v\'erifie les conditions 1),
2) du Lemme pr\'ec\'edent \'equivaut \`a:\ste 1) Il existe un \'el\'ement
$g_{1}\in
G_{1}$, tel que
$t=g_{1}^{\sigma_{1}}g_{1}^{-1}$. En particulier on a
$t^{\sigma_{1}}=t^{-1}$ \ste 2) $\Gamma_{+}=\Sigma_{1} $, et
$A$ est l'identit\'e.
\ste 3) $\i_{\a}=\{0\}$\ste
4) $$\theta(\Gamma'_{+})=\Gamma'_{-},\>\>et \>\> \theta
(A'\alpha)=A'^{-1}(\theta(\alpha)),\>\
\varepsilon_{A'\alpha}t^{A'\alpha}=
\varepsilon_{\alpha}t^{\alpha} ,
\>\> \alpha\in
\Gamma'_{+}$$ 5) L'espace $\i'_{\a'}$ est stable par $j_{1}$.\ste
\ste
\end{lem}
\dem
Soit ${\underline {\cal T}}={\cal
T}_{{\cal BD}, (t_{1},t_{2})({\cal W})}$ un triple de Manin
fortement standard, v\'erifiant les conditions 1), 2) du Lemme
pr\'ec\'edent, o\`u $t_{1},t_{2}\in J_{1}$. Alors  ${\underline {\cal
T}}$ est conjugu\'e par $(t_{1},t_{1}^{\sigma_{1}})\in \eta_{1}(G_{1})$
\`a
${\cal T}:={\cal
T}_{{\cal BD},(1,t)({\cal W})}$, o\`u $t=t_{1}^{-\sigma_{1}}t_{2}$,
qui v\'erifie les conditions 1) et 2) du Lemme pr\'ec\'edent.
Montrons que cela implique les propri\'et\'es 1) \`a 5) ci-dessus.\ste
Ecrivons ${\cal T}=(B,\i,\i')$. On obtient 1) en \'ecrivant que $\i$ est
conjugu\'e par un \'el\'ement de
$\eta_{1}(G_{1}) $ \`a $diag(\g_{1})$.\ste
Alors $\i$ est semi-simple, donc sous $\g$. Par suite, avec les
notations du Th\'eor\`e
me 1, on a $\m=\g$ et
$\Gamma_{+}=\Sigma_{1}$. Alors, utilisant les notations du
Lemme 22,  $\i=\{(X,\tau(X))\vert X\in \g_{1}\}$. Comme $\i$ est
conjugu\'e par un \'el\'ement de $\eta_{1}(G_{1})$ \`a $diag(\g_{1})$,
cela implique que l'automorphisme $\tau $ de $\g_{1}$ est
int\'erieur. Or, par d\'efinition, $\tau$ pr\'eserve $\j_{1}$, et
l'inverse du transpos\'e de sa restriction \`a $\j_{1}$ induit $A$
sur $\Sigma_{1}$. Comme $A$ pr\'eserve $\Sigma_{1}$, $A$ doit \^etre
l'identit\'e. On a donc prouv\'e 2). 3) r\'esulte du fait que $\m=\g$,
donc $\a=\{0\}$.\ste  On traduit maintenant le fait que $\i'$
est stable par $j_{1}$. Cela implique que $\p'$ est stable par
$j_{1}$.  Comme $\j_{0}$ est stable par $j_{1}$, joint \`a
(4.25), cela implique aussi que $\m'$ est stable par $j_{1}$.
Cela conduit imm\'ediatement \`a la premi\`ere \'egalit\'e de 4).
Les deux autres sont obtenues en traduisant la
stabilit\'e de $\h'$ par $j_{1}$. Plus pr\'ecis\'ement on pose :
$$U=(X_{\alpha}, t^{A'\alpha}X_{A'\alpha})\in \h',\>\> \alpha \in
\Sigma_{1}$$
On a, en tenant compte de (4.27) et de l'antilin\'earit\'e de
$\sigma_{1}$ :  $$j_{1}(U)=({\overline
{t^{A'\alpha}}}\varepsilon_{A'\alpha}Y_{\theta(A'\alpha)},
\varepsilon_{\alpha}Y_{\theta(\alpha)})$$
La stabilit\'e de $\h'$ par $j_{1}$, implique que
$\theta(A'\alpha)$ est \'el\'ement de $\Gamma'_{+}$ et
$j_{1}(U)$ doit \^etre un multiple scalaire de :
$$V= (Y_{\theta(A'\alpha)},t
^{A'\theta(A'\alpha)}Y_{A'\theta(A'\alpha)})$$
Les deux premi\`eres \'egalit\'es de 4) en r\'esulte imm\'ediatement
et l'on
obtient en outre :
$$\varepsilon_{A'\alpha}t^{A'\alpha}=
\varepsilon_{\alpha}t^{\theta(\alpha)} $$
En utilisant (4.26) et le fait que $t^{\sigma_{1}}=t^{-1}$, on aboutit \`a
la troisi\`eme \'egalit\'e de 4).\ste
La condition 5) est imm\'ediate.\ste  On proc\`ede de
m\^eme pour la r\'eciproque.\qed
\begin{lem}
Soit  $t$, ${\cal BD}$ v\'erifiant les
conditions 1) \`a 5) du Lemme pr\'ec\'edent. Soit $t'\in J_{1}$ tel
que :
$$t'^{A'\alpha}=t'^{\alpha},\>\> \alpha \in \Gamma'_{+}$$
Alors, on a :  ${\cal T}_{{\cal BD}, (1,t)({\cal W})}=
{\cal T}_{{\cal BD}, (t',t't)({\cal W})}$, et ${\cal T}_{{\cal BD},
(1,t)({\cal W})}$ est conjugu\'e par $(t',t'^{\sigma_{1}})\in
\eta_{1}(G_{1})$ \`a ${\cal T}_{{\cal BD},
(1,(t't'^{-\sigma_{1}})t)({\cal W})}$.
\end{lem}
\dem
Notons ${\cal T}_{{\cal BD}, (1,t){\cal W}}=(B,\i,\i')$ et
utilisons les notations du Th\'eor\`eme 1.  La stabilit\'e de $\i$ par
$(t',t')$ est claire. Par ailleurs,  $(t',t')\in J_{0}$ laisse
stable  $\p'$, donc $\n'$ et laisse fixe point par point les
\'el\'ements de $\a'$. Il reste \`a voir que $\h'$ est invariant. Mais
cette alg\`ebre est engendr\'ee par :
$$(X_{\alpha},t^{A'\alpha}X_{A'\alpha}),\>\>
(Y_{\alpha},t^{-A'\alpha}Y_{A'\alpha}), \>\> \alpha\in
\Gamma'_{+}$$
Le Lemme en r\'esulte imm\'ediatement. \qed
Si $\alpha\in \Gamma_{0}:=\Gamma'_{+}\cup \Gamma_{-}$, on note
${\cal C}(\alpha)$, l'ensemble des $\beta \in \Gamma_{0}$ tels qu' il
existe $n\in \N$ tel que  :
$$\alpha, A'\alpha, \dots, A'^{n-1}\alpha \in \Gamma'_{+}, \>\>
et \>\>\beta=A'^{n}\alpha $$
ou  :
$$\beta, A'\beta , \dots, A'^{n-1}\beta \in \Gamma'_{+}, \>\> et\>\>
\alpha=A'^{n}\beta $$
Les ${\cal C}(\alpha)$ sont soit distincts soit confondus. De plus la
deuxi\`eme \'egalit\'e de la condition 4) du Lemme 25, montre
facilement :
\beq \theta({\cal C}(\alpha)={\cal C}(\theta{\alpha})\eeq
(cf [P1], Lemme 4.11, pour un Lemme analogue)
Par ailleurs, gr\^ace \`a la d\'efinition de  de
$\varepsilon_{\alpha}$, on a :
\beq \varepsilon_{\theta(\alpha)}=\varepsilon_{\alpha}, \>\>
\alpha\in \Sigma_{1}\eeq
Le Th\'eor\`eme  suivant est du \`a A. Panov (cf. [P], Th\'eor\`emes 4.5,
4.13). \begin{theo}\ste
(i) Le complexifi\'e d'un $\h_{1}$-triple est conjugu\'e par un
\'el\'ement de $\eta_{1}(G_{1})$ \`a un triple
${\cal T}_{{\cal BD}, (1,u)({\cal W})}$, o\`u $u$, ${\cal BD}$
v\'erifient les conditions 1) \`a 5) du Lemme 25 (avec $t$ remplac\'e
par $u$) et $u$ v\'erifie de plus :\ste
1) $u^{\sigma_{1}}=u=u^{-1}$ \ste
2) $u^{\alpha}=1 \>\> ou\>\> -1, \>\>  si \>\> \alpha\in
\Sigma_{1}$
\ste 3) $u^{\alpha}=1$, si $\alpha\in\Sigma_{1}\setminus
\Gamma_{0}$ et $\theta(\alpha)\not= \alpha$ , ou si $\alpha \in
\Gamma_{0}$ et $\theta({\cal C}(\alpha)\not ={\cal C}(\alpha)$
\ste (ii)
R\'eciproquement si $u$ et ${\cal BD}$ v\'erifient les propri\'et\'es
ci-dessus (i.e. 1) et 2) de (i) et les conditions 1) \`a 5) du
Lemme 25),
${\cal T}_{{\cal BD}, (1,u){\cal W}}$ est conjugu\'e au complexifi\'e d'un
$\h_{1}$-triple, unique modulo la conjugaison de
$G_{1}^{\sigma_{1}}$.
\ste (iii) Dans le cas o\`u $\f_{1}$ est l'alg\`ebre
de Lie d'un tore maximal compact dans $G_{1}$, $\theta$ est triviale,
$\Gamma_{+}$ est vide et les conditions sur $u$ ci-dessus se
r\'eduisent aux deux premi\`eres, outre celles du Lemme 25.  \end{theo}
\dem
On sait que le complexifi\'e d'un $\h_{1}$-triple est conjugu\'e par
un \'el\'ement de $\eta_{1}(G_{1})$ \`a un triple
${\cal T}_{{\cal BD}, (1,t){\cal W}}$, o\`u $t$, ${\cal
BD}$ v\'erifient les conditions 1) \`a 5) du Lemme 26 (avec $t$
remplac\'e par $u$. \ste L'id\'ee est d'appliquer le Lemme 25, avec
$t'$ bien choisi.Il  suffit de d\'efinir $t'^{\alpha}, \alpha
\in \Sigma_{1}$.\ste
On \'etablit d'abord quelques r\'esultats auxiliaires.\ste
On rappelle que d'apr\`es (4.26) et la condition 1) du Lemme 25,
on a : \beq {\overline
{t^{\theta(\alpha)}}}=t^{\alpha},\>\> \alpha \in
\Sigma_{1}\eeq
Par ailleurs, d'apr\`es la troisi\`eme \'egalit\'e de la condition 4) du
Lemme 25 , on a :  $$
t^{A'\alpha}=\varepsilon_{A'\alpha}\varepsilon_{\alpha}t^{\alpha},
\>\> \alpha \in \Gamma'$$
d'o\`u l'on d\'eduit :
\beq t^{\beta}=\varepsilon _{\beta}\varepsilon_{\beta
'}t^{\beta '}, \>\> \beta, \beta'\in {\cal C}(\alpha), \>\> \alpha \in
\Gamma_{0} \eeq
Enfin la condition 4) du
Lemme 25    montre que :
\beq \varepsilon_{\alpha}= \varepsilon_{\theta(\alpha)}, \>\>
\alpha \in \Sigma_{1}\eeq Pour d\'efinir $t'^{\alpha}$, on distingue
plusieurs cas. On note $\Sigma_{1*}$, un sous-ensemble de
$\Sigma_{1}$ tel que  tout $\beta \in \Sigma_{1}$ soit \'el\'ement
 d'un  ${\cal C}(\alpha)$ pour un unique $\beta$ appartenant \`a
$\Sigma_{1*}\cup\theta( \Sigma_{1*})$, et tel que les \'el\'ements de
l'intersection de $\Sigma_{1*}$ et $\theta( \Sigma_{1*})$ soient
fix\'es par $\theta$.\ste
 Soit
$\alpha\in \Sigma_{1*}$  :
\ste 1) Si
$\alpha\notin
\Gamma_{0}$, et
$\theta(\alpha)=\alpha$,  (4.31) implique que  $t^{\alpha}$ est
r\'eel et on pose :
$$t'^{\alpha}=\vert t^{-\alpha}\vert^{1/2}$$
 2) Si $\alpha\notin \Gamma_{0}$,
et
$\theta(\alpha)\not=\alpha$, on note $z$ une racine carr\'ee de
$t^{-\alpha}$, et on pose :
$$j^{\alpha}=z,\>\>j^{\theta(\alpha)}={\overline z}$$
 3) Si $\alpha\in \Gamma_{0}$ et $\theta ({\cal
C}(\alpha))={\cal
C}(\alpha)$, on a,  d' apr\`es (4.32) :
$$t^{\theta(\alpha)}=
\varepsilon_{\theta(\alpha)}\varepsilon_{\alpha}t^{\alpha}$$
Tenant compte de (4.31),  cela implique que $t^{\alpha}$ est
r\'eel. On note $z$ une racine carr\'ee de $\vert t^{-\alpha}\vert $,
et on pose : $$t'^{\beta}=z, \>\> \beta\in {\cal
C}(\alpha)$$
4)  Si $\alpha\in \Gamma_{0}$ et $\theta (C(\alpha))\cap
C(\alpha)=\emptyset $, on note $z$ une racine carr\'ee de
$t^{-\alpha}$, et on pose :
$$t'^{\beta}=z, \>\> t'^{\theta(\beta)}={\overline z}, \>\>\beta
\in C(\alpha)$$
On voit que les relations pr\'ec\'edentes d\'efinissent $t'$,
qui v\'erifie :
$$t'^{A'\alpha}=t'^{\alpha}, \>\> \alpha \in \Gamma'_{+}$$
Le Lemme 25 et les conditions impos\'ees \`a $t'$ montrent que
$u:=t't'^{-\sigma_{1}}t$ v\'erifie les conditions voulues.\ste
(ii) est un cas particulier du Lemme 25 (ii).
\ste
Traitons le cas o\`u $\f_{1}$ est l'alg\`ebre de Lie d'un tore maximal
compact de $\g_{1}$. Alors $\f_{1}$ est l'ensemble des \'el\'ements de
$\j_{1}$ sur lesquels toutes les racines sont imaginaires pures.
\ste Tout
$\h_{1}$-triple  est conjugu\'e sous
$G_{1}$ \`a un triple r\'eel fortement standard $(B_{1},
\i_{1},\i'_{1})$. Alors
$\i_{1}\cap\j_{1}$, qui est conjugu\'e par $G_{1}$ \`a $\f_{1}$, est
l'alg\`ebre de Lie d'un tore maximal compact de $G_{1}$, donc est \'egal
\`a
$\f_{1}$. Par ailleurs $\h'_{1}\cap\j_{1}$ est une \sac fondamentale
de $\h'_{1}$. Son intersection avec  $\f_{1}$ est donc non
r\'eduite \`a z\'ero sauf si $\h'_{1}$ est r\'eduite \`a
z\'ero. On en d\'eduit que la sous-alg\`ebre Lagrangienne $\i'_{1}$ est
sous
une alg\`ebre de Borel. Par complexification, il en r\'esulte que dans
(i), on doit avoir
$\Gamma'_{+}=\emptyset$
 La d\'efinition de $\theta$ et le fait que les
racines soient imaginaires pures sur $\f_{1}$ montrent que $\theta$
est l'identit\'e. Ceci ach\`eve la preuve du Th\'eor\`eme.\qed

\ste {\bf R\'ef\'erences }\ste

\noindent[BD], BELAVIN A., DRINFELD G., {\em Triangle equations and
simple Lie algebras}, Mathematical Physics Reviews, vol. 4,
 93-165

\noindent[Bor], BOREL A., {\em Linear algebraic groups, Second
Enlarged Edition}, Graduate Text in Math.126, 1991, Springer Verlag,
 New York, Berlin, Heidelberg.

\noindent[Bou], BOURBAKI N., {Groupes et Alg\`ebres de Lie, Chapitre
I, Chapitres IV, V, VI, Chapitres VII, VIII}, Actualit\'es
Scientifiques et Industrielles 1285, 1337, 1364, Hermann, Paris,
1960, 1968, 1975.

\noindent[De], DELORME P., {Sur les triples de Manin pour une
alg\`ebre r\'eductive complexe}, Preprint 1999.

\noindent[G],
GANTMACHER F.,  {\em Canonical representation of automorphism of
a semisimple Lie group}, Math Sb., 47, (1939), 101-144.

\noindent[K1], KAROLINSKY E., {\em A classification of Poisson
homogeneous spaces of a compact Poisson Lie group}, Math.
Phys., Anal. and Geom., 3 (1996), 545-563.

\noindent[K2], KAROLINSKY E., {\em A classification of Poisson
homogeneous spaces of a compact Poisson Lie group},
Dokl. Ak. Nauk, 359 (1998), 13-15.

\noindent[K3], KAROLINSKY E.,{\em A classification of Poisson
homogeneous spaces of a reductive complex Poisson Lie group,}
 Preprint, 1999

\noindent[M1], MATSUKI T., {\em The orbits of affine symmetric
spaces under the action of minimal parabolic subgroups}, J. Math.
Soc. Japan, 31 (1979), 331-357.

\noindent[M2], MATSUKI T., {\em Orbits of affine symmetric
spaces under the action of  parabolic subgroups}, Hiroshima J.
Math., 12 (1982), 307-320.

\noindent[P1], PANOV A.., {\em Manin triples of real simple Lie
algebras, Part 1}, Preprint.

\noindent[P2], PANOV A., {\em Manin triples of real simple Lie
algebras, Part 2}, Preprint, \ste QA 9905028.

\noindent[W], WARNER G., {\em Harmonic Analysis on semi-simple Lie
groups}, Grundleh
ren der math. Wis. in Einz., Vol 188, Springer
Verlag, Berlin-Heidelberg-New York 1972

\ste {\em Institut de Math\'ematiques de Luminy, U.P.R. 9016 du
C.N.R.S. \ste
Universit\'e de la M\'editerrann\'ee, 163 Avenue de Luminy, Case
907,\ste
13288, Marseille Cedex 09, France \ste
e-mail : delorme@iml.univ-mrs.fr}
\end{document}